\documentclass[a4,10pt]{article}
\usepackage{graphicx}
\usepackage{appendix}
\usepackage{amsmath,amssymb}
\usepackage{ulem}

\usepackage[nameinlink,noabbrev]{cleveref}%

\usepackage{color}
\usepackage{here}

\makeatletter
\def\@seccntformat#1{\@ifundefined{#1@cntformat}%
   {\csname the#1\endcsname\quad}%      default
   {\csname #1@cntformat\endcsname}%    enable individual control
}
\makeatother

\begin{document}
\setlength{\baselineskip}{15pt} %20pt
\def\pot{\mathaccent"7017} %={}^o 
\makeatletter
 \renewcommand{\theequation}{%
  \thesection.\arabic{equation}}
 \@addtoreset{equation}{section}
\makeatother

\newtheorem{remark}{Remark}[section]
\newtheorem{theo}{Theorem}[section]
\newtheorem{lemma}{Lemma}[section]
\newtheorem{prop}{Proposition}[section]
\newtheorem{assume}{Assumption}[section]
\newtheorem{cor}{Corollary}[section]
\newtheorem{example}{Example}[section]

\def\ep{\varepsilon}
\def\ov{\overline}
\def\un{\underline}
\def\del{\partial}
\def\norm{\parallel}
\def\no{\noindent}
\def\lam{\lambda}
\def\dis{\displaystyle}
\def\bhat{\widehat}
\def\lap{\bigtriangleup}
\def\R{\mbox{\bf R}}
\def\C{\mbox{\bf C}}
\def\lg{\langle}
\def\rg{\rangle}

\newcommand{\Qed}{ $\square$ }

\title{
Stability of Single Transition Layer in Mass-Conserving
Reaction-Diffusion Systems with Bistable Nonlinearity}

\author{
Hideo Ikeda\thanks{
Department of Mathematics, University of Toyama, Toyama,
930-8555, Japan,
email: hideoikeda5@gmail.com}\;
and 
Masataka Kuwamura\thanks{
Graduate School of Human Development and Environment,
Kobe University, Kobe 
657-8501, Japan,
email: kuwamura@main.h.kobe-u.ac.jp }
}

\date{}

\maketitle
\begin{abstract}
Mass-conserving reaction-diffusion systems with bistable nonlinearity are considered 
under general assumptions. The existence of stationary solutions with a
single internal transition layer in such reaction-diffusion systems is shown using the analytical singular perturbation theory.
Moreover, a stability criterion for the stationary solutions is provided by calculating the Evans function. 
\end{abstract}

\vspace{1ex}
%\noindent
%Abbreviated title:  Stability of Transition Layer in Reaction-Diffusion Systems \\
%
\hspace{1.75cm} stability, Evans function, analytical singular perturbation method 
Corresponding author: Hideo Ikeda (hideoikeda5@gmail.com)

\section{Introduction }
Reaction-diffusion systems provide a theoretical framework for understanding pattern formation in various fields of science and technology. Here, we consider the following reaction-diffusion system: 
\begin{equation}\label{a1}
\left \{
\begin{array}{rcl}
\begin{array}{l}
u_t = \ep^2 u_{xx} +  f(u,v), \\[1ex]
v_t = Dv_{xx}  - f(u, v), 
\end{array} \quad  (t,x) \in (0,\infty) \times  (0,1) \\[0.4cm]
(u_x, v_x)(t,0) = (0, 0) = (u_x, v_x)(t,1), \ 0 < t < \infty, 
\end{array} 
\right.
\end{equation}
where $\ep$ and $D$ are positive constants satisfying $0 < \ep \ll D$. 
We note that \eqref{a1} is a mass-conserving reaction-diffusion system because 
\begin{equation}\label{a2}
\xi := \int_0^1 \left\{ u(0, x) + v(0, x) \right\} dx =  
\int_0^1 \left\{ u(t, x) + v(t, x) \right\} dx 
\end{equation}
holds for any (smooth) solutions.
It was shown in \cite{MJE1} that \eqref{a1}
with an appropriate nonlinear term $f$ is a useful model for 
understanding the wave-pinning phenomena in cell division and differentiation. 
In simple terms, this biological phenomena will be mathematically
interpreted as the dynamics of  \eqref{a1} as follows: a propagating front solution such that the initial position of the front exists near the boundary converges to a stationary solution with a single internal 
transition layer.  In fact, \cite{MJE1, MJE2} concluded that the system \eqref{a1}  with an appropriate nonlinear term $f$ 
has a stable stationary solution with a single internal transition layer under certain conditions using a formal analysis, numerical computations, and a perturbative argument against the background of cell biology.  

In this paper, we consider the existence and stability of a single internal transition layer solution of the system \eqref{a1} with more general 
nonlinear term $f$ including the specific one given in \cite{MJE1, MJE2}. 
According to our previous paper \cite{KTI},
we assume that
the nonlinear term $f$ has the following bistability: 

\noindent
{\bf Assumption}. 
{\rm (A1)} 
The ODE $\dot{u} = f(u, v)$
is bistable in $u$ for each fixed $v \in I = (\un{v}, \ov{v})$.
That is, $f(u, v)=0$ has exactly three roots 
$h^-(v) < h^0(v) < h^+(v) $ for each $v \in I$ satisfying
\begin{equation*}%\label{a3}
f_u ( h^{\pm}(v), v ) < 0
\ \ \
\text{and}
\ \ \ 
f_u ( h^{0}(v), v ) > 0. 
\end{equation*}
{\rm (A2)} The function 
\begin{equation}\label{a4}
J(v) := \int_{h^-(v)}^{h^+(v)} f(u, v)du \ \ \ ( v \in I)
\end{equation}
has an isolated zero at $v = v^* \in I$ such that 
\begin{equation}\label{a4x}
J'(v^*) = \int_{h^-(v^*)}^{h^+(v^*)} f_v(u, v^*)du \neq 0.
\end{equation}
\no
{\rm (A3)} \ 
\begin{equation*}%\label{a5x}
f_u ( h^{\pm}(v), v ) < f_v ( h^{\pm}(v), v )  \ \ \ ( v \in I). 
\end{equation*}
\no
{\rm (A4)} 
The conserved mass $\xi$ defined by \eqref{a2} satisfies the following inequality:
\begin{equation}\label{a6}
h^-(v^*) + v^* < \xi < h^+(v^*) + v^*.
\end{equation}

To justify the formal and numerical result given in \cite{MJE1, MJE2},  
we showed in \cite{KTI} that under the assumptions (A1),(A2) and (A4), 
for any given $\xi$ satisfying \eqref{a6},
the system \eqref{a1} has a family of single transition layer solutions $(u,v)(x;\ep)$ 
satisfying 
\begin{equation}\label{b1_1}
\xi = \int_0^1 \left\{ u(x; \ep) + v(x; \ep) \right\} dx,
\end{equation}
for sufficiently small $\ep > 0$. These solutions satisfy
$
\lim_{\ep \to 0} u(x; \ep) = U^*(x)
$
compact uniformly on $[0,x^*) \cup (x^*,1]$ and 
$
\lim_{\ep \to 0} v(x; \ep) = v^*
$
uniformly on $[0, 1]$, where 
$$
U^*(x) 
= 
\left \{
\begin{array}{l}
h^-(v^*)  \ \ \ (0 \leq x \leq x^*) \\[1ex]
h^+(v^*)  \ \ \ (x^* < x \leq 1) 
\end{array} 
\right.
$$
and
$
x^* = (h^+(v^*) + v^* - \xi)/(h^+(v^*) - h^-(v^*)) 
$
is derived from \eqref{b1_1} in the limit of $\ep \to 0$.
We note that $(u,v)(x;\ep)$ are called jump-up transition layer
solutions (see, Figure~\ref{fig.1}).
Moreover, we showed that they are stable if $J'(v^*) > 0$ under the assumptions (A1)-(A4). 
Furthermore, \cite{KTI} investigated examples for both $J'(v^*) > 0$ and $J'(v^*) < 0$ with the aid of numerical simulations, 
and proposed a conjecture
that they are unstable if $J'(v^*) < 0$ under the assumptions (A1)-(A4), i.e., the sign of $J'(v^*)$ gives a stability criterion for the single transition layer solutions.

Generally, there are two different methods for solving the singularly perturbed 
eigenvalue problem concerning the stability of transition 
layer solutions of reaction-diffusion systems: 
an analytical approach called the 
Singular Limit Eigenvalue Problem (SLEP) method \cite{NF, NMIF}, and a geometrical 
approach called the Nonlocal Eigenvalue Problem (NLEP) method \cite{DGK}. 
These methods are based on the linearized stability principle; 
the former method solves the linearized eigenvalue problem directly and 
derives a well-defined singular limit equation called the SLEP equation 
as $\varepsilon \to 0$, while the latter method defines the Evans function \cite{AGJ}
for the linearized equations and subsequently applies a topological method to it. 
Although the NLEP method can apply to a broader class of reaction-diffusion systems,
we used the SLEP method in \cite{KTI} because
it provides detailed information on the behavior of 
critical eigenvalues, which essentially determine the stability of transition 
layer solutions. 
Consequently, we obtained the precise characterization of the 
critical eigenvalues, which leads to the stability result for 
the transition layer solutions $(u,v)(x;\ep)$.
However, we could not show the existence of the critical eigenvalues, 
and hence we could not obtain the instability result for them.
In fact, the Lax-Milgram theorem cannot apply to the solvability of the SLEP equation. This shortcoming seems to be common to singular perturbation 
problems for mass-conserving reaction-diffusion systems such as \eqref{a1}. 

\par The aim of this paper is 
to prove that the sign of $J'(v^*)$ gives a stability 
criterion for the family of single transition layer solutions $(u,v)(x;\ep)$ satisfying \eqref{b1_1} for 
sufficiently small $\ep > 0$. Here, we 
calculate the Evans function $g(\ep;\lambda)$ (see Section \ref{S3}) for the linearized eigenvalue problem in the same spirit of \cite{I2}. 
Although its definition is different from that in \cite{AGJ}, it was shown in \cite{INS}
that the principal parts of both functions are equivalent up to a constant multiple. 
Our main result concerning the existence and stability of single transition layer solutions is summarized as follows: 
\begin{theo}\label{th1}
Under the assumptions (A1) - (A4), for any given $\xi$ satisfying \eqref{a6}, 
the mass-conserving reaction-diffusion system \eqref{a1} has a family of single transition layer solutions $(u,v)(x;\ep)$ satisfying 
\eqref{b1_1} 
for sufficiently small $\ep > 0$. Moreover,
$(u,v)(x;\ep)$ are stable if $J'(v^*) > 0 $, while they are unstable if $J'(v^*) < 0 $.
\end{theo}

\begin{remark} \rm
The solutions $(u,v)(x;\ep)$ are called jump-down transition layer
solutions when they satisfy
$
\lim_{\ep \to 0} u(x; \ep) = U^*(x)
$
compact uniformly on $[0,x^*) \cup (x^*,1]$ and 
$
\lim_{\ep \to 0} v(x; \ep) = v^*
$
uniformly on $[0, 1]$, where 
$$
U^*(x) 
= 
\left \{
\begin{array}{l}
h^+(v^*)  \ \ \ (0 \leq x \leq x^*) \\[1ex]
h^-(v^*)  \ \ \ (x^* < x \leq 1) 
\end{array} 
\right.
$$
and $
x^* = (h^-(v^*) + v^* - \xi)/(h^-(v^*) - h^+(v^*))
$
is derived from \eqref{b1_1} in the limit of $\ep \to 0$.
In this case, we can also obtain the same result as shown in 
the above theorem.
Moreover, as mentioned in \cite{KTI}, 
we can easily find $\ep^2 u(x; \ep) + D v(x; \ep) \equiv C(\ep)$,
where $C(\ep)$ denotes a constant independent of $x$. 
On the other hand, $u(x; \ep)$ exhibits
a single internal transition layer with $O(1)$-amplitude at $x = x^*$ when $\ep$ is sufficiently small. 
Therefore, noting  $v(x; \ep) = (C(\ep) - \ep^2 u(x; \ep))/D$, 
$v(x; \ep)$ exhibits
a single internal transition layer with $O(\ep^2)$-amplitude at $x = x^*$ when $\ep$ is sufficiently small.
However, our stability analysis does not require the information about
this small amplitude transition layer.
\end{remark}

The precise version of this theorem is given by Theorem \ref{th2} and Theorem \ref{l3-8}. Theorem \ref{th1} justifies the conjecture in \cite{KTI}, that is, 
the stability of the single transition layer solutions is determined by the sign of $J'(v^*)$. 
The remainder of this paper is organized as follows. In the next section, we construct a family of stationary solutions with a single internal transition layer
by using the analytical singular perturbation method \cite{F, I1, I2, KTI, MTH}, 
which gives $O(\ep)$ uniform approximations for the stationary 
solutions.
They are different from those in \cite{KTI}, and more suitable for calculating the Evans function.
In Section~\ref{S3}, we show the stability of the stationary solutions by calculating the Evans function $g(\ep;\lambda)$ for the linearized eigenvalue problem around them. 
The calculations can be performed in the same spirit of \cite{I2}. However, they require delicate analysis and advanced technique because they involve four-dimensional 
linear ODE systems concerning the linearized eigenvalue problem;
the difficulty of them would reflect the fact that 
the SLEP method with no justification by the Lax-Milgram theorem 
cannot show the existence of the critical eigenvalues \cite[Remark 3.1]{KTI}.
The proof of Lemma \ref{l3-6}, which plays a key role in the 
calculations of the Evans function, is found in the appendix.
We emphasize that the present paper and \cite{KTI} 
give a comparison between the SLEP method and the approach based on the Evans function;
the Evans function approach can provide a necessary and sufficient
condition for the stability of 
transition layer solutions of mass-conserving reaction-diffusion systems, 
whereas the SLEP method can provide only a sufficient condition for their stability.

\section{Existence of single transition layer solutions } % Section 2 

\begin{figure}[!htb]
\centering
\includegraphics[width=8cm]{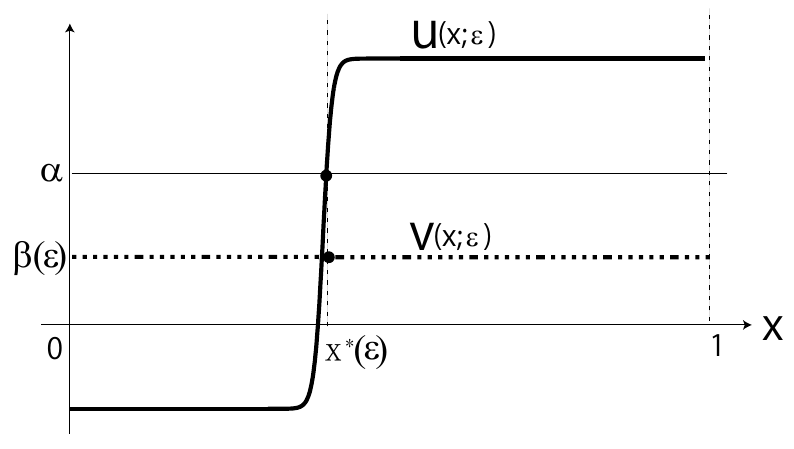} 
\caption{Schematic profile of a single jump-up transition layer solution. This profile does not represent the $O(\ep^2)$-amplitude transition layer of $v$-component 
because it is not required in our analysis.}
\label{fig.1} 
\end{figure}

In this section, we consider a single transition layer solution of 
\begin{equation}\label{b1}
\left \{
\begin{array}{l}
\begin{array}{l}
\ep^2 u_{xx} +  f(u,v) = 0, \\[1ex]
Dv_{xx}  - f(u, v) = 0, 
\end{array} \quad  x \in (0,1) \\[0.4cm]
(u_x, v_x)(0)  = (0, 0) = (u_x, v_x)(1)
\end{array} 
\right.
\end{equation}
satisfying \eqref{b1_1} 
for a given constant $\xi$ in (A4) under the assumptions (A1) and (A2). 
We assume that the $u$-component of a solution $(u,v)(x;\ep)$ of \eqref{b1} exhibits a sharp jump-up transition layer 
with $O(1)$-amplitude 
at $x= x^*(\ep) \in (0,1)$ (see Figure \ref{fig.1}).\par

To solve this problem, in Section \ref{S2.1} we divide the interval $[0, 1]$ into two subintervals $[0, x^*(\ep)]$ and $[x^*(\ep), 1]$, and consider the following two boundary value problems: 
\begin{equation}\label{b5}
\left \{
\begin{array}{l}
\begin{array}{l}
\ep^2 u_{xx} + f(u,v) = 0, \\[1ex]
Dv_{xx} - f(u,v) = 0,
\end{array} \ x \in  (0, x^*(\ep))   \\[0.4cm]
(u_x, v_x)(0) = (0,0), \ (u,v)(x^*(\ep)) = (\alpha, \beta(\ep))
\end{array} 
\right.
\end{equation}
and
\begin{equation}\label{b6}
\left \{
\begin{array}{l}
\begin{array}{l}
\ep^2 u_{xx} + f(u,v) = 0, \\[1ex]
Dv_{xx} - f(u, v) = 0,
\end{array} \ x \in  (x^*(\ep), 1) \\[0.4cm]
(u,v)(x^*(\ep)) = (\alpha,\beta(\ep)), \ (u_x,v_x)(1) = (0,0), 
\end{array} 
\right.
\end{equation}
where $\alpha$ is a constant satisfying 
$h^-(\beta(0)) < \alpha < h^+(\beta(0))$
and $\beta(\ep)$ is determined by $v(x^*(\ep)) = \beta(\ep)$. Put 
\begin{equation}\label{b8}
x^*(\ep) = x_0 +  \ep x_1
\end{equation}
and
\begin{equation}\label{b9}
\beta(\ep) = \beta_0 + \ep \beta_1.
\end{equation}
By using the singular perturbation method used in \cite{F, I1, I2, KTI, MTH}, 
we show the existence of solutions satisfying \eqref{b5} and \eqref{b6}. In Section \ref{S2.2}, we match these solutions in $C^1$-sense at $x = x^*(\ep)$, from which we find an approximate solution of \eqref{b1} up to $O(\ep)$ by determining the relations between $x_j$ and $\beta_j$ for $j=0,1$. Similarly, we use the equation \eqref{b1_1}, and obtain the other relations between $x_j$ and $\beta_j$ for $j=0,1$ in Section \ref{S2.3}. 
Finally, in Section \ref{S2.4}, using the result in Sections \ref{S2.2} and \ref{S2.3}, we determine the unknown constants $x^*(\ep)$ and $\beta(\ep)$ uniquely, and obtain the desired result about the existence of a single transition layer solution. Here we emphasize that $O(\ep)$ uniformly approximate solutions are required for the stability analysis (see Section \ref{S3}). \par 
In \cite{KTI}, we solved a single equation for only $u$, which comes from 
$\ep^2 u + Dv = Const.$ obtained by adding the first and second equations of \eqref{b1},
under the constrained condition \eqref{b1_1}.
But here, we 
will solve the full problem \eqref{b1} with \eqref{b1_1} because the solving procedure is very important for the stability analysis.

We use the following function spaces with a positive $\nu$ in this paper: 
\begin{equation*}
\begin{array}{lcl}
C^2_\nu[0,1] & := & \left\{ u \in C^2[0,1] \ | \ \dis \sum_{j=0}^2 \max_{0 \leq x \leq 1} \left| \left( \nu \frac {d}{dx}\right)^{\! j} \! u(x) \right| < \infty \right\}, \\[3.5ex]
\pot{C}^2_\nu[0,1] & := & \left\{ u \in C^2_\nu[0,1] \ | \ u_x(0) = 0, \ u_x(1) = 0 \right\},  \\[2ex]
C^2_{\nu,0}[0,1] & := & \left\{ u \in C^2_\nu[0,1] \ | \ u_x(0) = 0, \ u(1) = 0 \right\},  \\[2ex]
C^2_{\nu,1}[0,1] & := & \left\{ u \in C^2_\nu[0,1] \ | \ u(0) = 0, \ u_x(1) = 0 \right\}.
\end{array} 
\end{equation*}

\subsection{Solutions of \eqref{b5} and \eqref{b6}}\label{S2.1} % Subsection 2.1 

First, we consider the approximation of the solution of \eqref{b5} up to $O(\ep)$. Applying the change of variables $x = x^*(\ep)y$, we have
\begin{equation}\label{b7}
\left \{
\begin{array}{l}
\begin{array}{l}
\ep^2 u_{yy} + (x^*(\ep))^2 f(u,v) = 0, \\[1ex]
Dv_{yy} - (x^*(\ep))^2 f(u,v) = 0,
\end{array} \ y \in  (0, 1)   \\[0.4cm]
(u_y, v_y)(0) = (0,0), \ (u,v)(1) = (\alpha, \beta(\ep)).
\end{array} 
\right.
\end{equation}
To construct the outer approximation of the solution of \eqref{b7},
substituting 
\begin{equation*}%\label{b10} 
\left \{
\begin{array}{l}
u(y) = U^-_0(y) + \ep U^-_1(y) + O(\ep^2), \\[1ex]
v(y) = V^-_0(y) + \ep V^-_1(y) + O(\ep^2)
\end{array} \right.
\end{equation*}
into \eqref{b7}, and comparing the coefficients of powers of $\ep$, we have 
\begin{equation}\label{b11} \left\{
\begin{array}{l}
f( U^-_0, V^-_0 ) = 0,\\[1ex]
D V^-_{0,yy} = 0, \ y \in (0,1) \\[1ex]
V^-_{0,y}(0) = 0, \ V^-_{0}(1) = \beta_0 
\end{array} \right.
\end{equation}
and
\begin{equation}\label{b12} \left\{
\begin{array}{l}
f_u^- U^-_1 + f_v^- V^-_1 = 0, \\[1ex]
D V^-_{1,yy} = 0, \ y \in (0,1) \\[1ex]
V^-_{1,y}(0) = 0, \ V^-_{1}(1) = \beta_1, 
\end{array} \right.
\end{equation}
where $f_u^- := f_u( U^-_0, V^-_0)$ and $f_v^- := f_v( U^-_0, V^-_0)$. 
Since we consider a jump-up solution at $y=1$, it follows from \eqref{b11} and
\eqref{b12} that
\begin{equation}\label{b14} 
\begin{array}{l}
V^-_0(y) = \beta_0, \ 
U^-_0(y) = h^{-}(\beta_0)
\end{array}
\end{equation}
and
\begin{equation}\label{b15}
\begin{array}{l}
V^-_1(y;\beta_1) = \beta_1, \ 
U^-_1(y;\beta_1) = - f_v^- \beta_1 / f_u^- = h^-_v(\beta_0)\beta_1,
\end{array}
\end{equation}
where we used the relation $f_u(h^-(v), v) h^-_v(v) + f_v(h^-(v), v) = 0$
obtained by the differentiation of $f( h^-(v), v) = 0$ in $v$.
It should be noted that $(U^-_j,V^-_j)$ $(j=0, 1)$ are constants independent of $y$. 
Since these outer approximations do not satisfy the boundary condition at $y = 1$, 
we must consider the correction of the above approximation to the solution of \eqref{b7} in a neighborhood of $y = 1$ 
with the aid of the inner approximation given by
\begin{equation}\label{b17} \left\{ 
\begin{array}{l}
u(y)  =  U^-_0(y) + \ep U^-_1(y;\beta_1)
+ \phi^-_0( \frac{y-1}{\ep}) + \ep \phi^-_1(  \frac{y-1}{\ep} ) + O(\ep^2), \\[1ex]
v(y)  =  V^-_0(y) + \ep V^-_1(y;\beta_1) 
+ \ep^2 \psi^-_0( \frac{y-1}{\ep}) 
 + \ep^3 \psi^-_1(  \frac{y-1}{\ep} ) + O(\ep^4). 
\end{array} \right. 
\end{equation}
Introducing the stretched coordinate $z = (y-1)/\ep$, and substituting 
\eqref{b17} into \eqref{b7}, and 
comparing the coefficients of powers of $\ep$, we have 
\begin{equation}\label{b18}
\left \{
\begin{array}{l}
\begin{array}{l}
\ddot{\phi}_0^- + x_0^2 \tilde{f}^- = 0, \\[1ex]
D \ddot{\psi}_0^- + \ddot{\phi}_0^- = 0, 
\end{array} \ z \in (-\infty, 0)  \\[0.4cm]
\phi^-_0(-\infty) = 0, \ \phi^-_0(0) = \alpha - U^-_0(1), \\[1ex]
\psi^-_0(-\infty) = 0, \ \dot{\psi}^-_0(-\infty) = 0 
\end{array} \right.
\end{equation}
and 
\begin{equation}\label{b19}
\left \{
\begin{array}{l}
\begin{array}{l}
\ddot{\phi}_1^- + x_0^2 \tilde{f}^-_u \phi^-_1 = F_1^-(z;\beta_1,x_1), \\[1ex]
D \ddot{\psi}_1^- + \ddot{\phi}_1^- = 0, 
\end{array}  \ z \in (-\infty, 0)   \\[0.4cm]
\phi^-_1(-\infty) = 0, \ \phi^-_1(0) = -U^-_1(1;\beta_1), \\[1ex]
\psi^-_1(-\infty) = 0, \ \dot{\psi}^-_1(-\infty) = 0,
\end{array} \right.
\end{equation}
where the dot notation denotes $d/dz$, and 
$$
F_1^-(z;\beta_1,x_1) :=  -2x_0 x_1 \tilde{f}^- - x_0^2 \tilde{f}^-_u U^-_1(1;\beta_1) -x_0^2 \tilde{f}^-_vV^-_1(1;\beta_1)  ,
$$
$$ \tilde{f}^- := f( h^-(\beta_0) + \phi^-_0, \beta_0 ), \ 
\tilde{f}^-_u := f_u( h^-(\beta_0) + \phi^-_0, \beta_0 ),
$$
and $\tilde{f}^-_v$ is similarly defined. 
From Assumption (A1), we find that \eqref{b18} has a unique monotone
increasing solution $\phi^-_0(z)$ if $h^-(\beta_0) < \alpha \le \bar{\alpha}(\beta_0)$, where $\bar{\alpha}(\beta_0)$ is defined by 
\begin{equation}\label{Ham1}
\displaystyle \int_{h^-(\beta_0)}^{\bar{\alpha}(\beta_0)} f(u,\beta_0) du = 0.
\end{equation}
We then have $\psi^-_0(z) = - \phi^-_0(z)/D$. Moreover, noting $(\dot{\phi}_0^-)_{zz} + x_0^2 \tilde{f}^-_u \dot{\phi}_0^- = 0$ obtained by differentiating the first equation of \eqref{b18} in $z$, 
we see that the solutions of \eqref{b19} are explicitly given by
\begin{equation}\label{b21} \left\{ 
\begin{array}{l}
\phi^-_1(z;\beta_1,x_1) = -U^-_1(1;\beta_1) \dis\frac{ \dot{\phi}_0^-(z) }{ \dot{\phi}_0^-(0) } \\[0.3cm]
\hspace*{2cm} - \ \dot{\phi}_0^-(z)  \dis\int_z^0 \frac{1}{ (\dot{\phi}_0^-(\eta))^2} 
\dis\int_{-\infty}^{\eta}  \dot{\phi}_0^-(\zeta) F_1^-(\zeta;\beta_1,x_1) d\zeta d\eta, \\[0.4cm]
\psi^-_1(z;\beta_1,x_1) = - \phi^-_1(z;\beta_1,x_1)/D. 
\end{array} \right.
\end{equation}
\par
Now, we put 
\begin{equation*} \left\{ 
\begin{array}{lcl}
U^-(y;\ep;\beta_1,x_1) & = & U^-_0(y) + \ep U^-_1(y;\beta_1) \\[0.2cm]
 &  & + \ \theta(y)\{ \phi^-_0( \frac{y-1}{\ep}) 
+  \ep \phi^-_1(  \frac{y-1}{\ep};\beta_1,x_1 )\} , \\[0.3cm]
V^-(y;\ep;\beta_1,x_1) & = & V^-_0(y) + \ep V^-_1(y;\beta_1) \\[0.2cm]
 &  & + \  \theta(y)\{ \ep^2 \psi^-_0( \frac{y-1}{\ep}) 
+  \ep^3 \psi^-_1(  \frac{y-1}{\ep};\beta_1,x_1 ) \\[0.2cm]
 & & - \ \ep^2 \psi^-_0(0) - \ \ep^3 \psi^-_1(0;\beta_1,x_1) \} ,
\end{array} \right.
\end{equation*}
where $\theta(y) \in C^{\infty}[0,1]$ is a cut-off function satisfying 
\begin{equation*}
\begin{array}{c}
\theta(y) = 0, \quad y \in [0, 1/2]; \quad \theta(y) = 1, \quad y \in [3/4, 1]; \\[1ex]
 0 \leq \theta(y) \leq 1, \quad y \in (1/2, 3/4). 
\end{array}
\end{equation*}
Here, we note that we are correcting $V^-$ to satisfy the boundary condition at $y=1$ of  \eqref{b7}.
Moreover, noting that when $y \in (1/2,3/4)$, 
$$
\left\{ 
\begin{array}{lcl}
U^-(y;\ep;\beta_1,x_1) = U^-_0(y) + \ep U^-_1(y;\beta_1) + O(e^{-\kappa/\ep}),
\\[1ex]
V^-(y;\ep;\beta_1,x_1) = V^-_0(y) + \ep V^-_1(y;\beta_1) + O(e^{-\kappa/\ep}) 
\end{array}
\right.
$$ 
for some $\kappa >0$, we see that $(U^-(y;\ep;\beta_1,x_1),V^-(y;\ep;\beta_1,x_1))$ is an $O(\ep)$ approximation to a solution of \eqref{b7}. We then find a solution of \eqref{b7} in the following form:
\begin{equation*} \left\{ 
\begin{array}{l}
\tilde{u}^-(y;\ep;\beta_1,x_1) = U^-(y;\ep;\beta_1,x_1) + \ep r^-(y;\ep;\beta_1,x_1) 
+ \ep h^-_v(V^-_0(y)) s^-(y;\ep;\beta_1,x_1), \\[1ex]
\tilde{v}^-(y;\ep;\beta_1,x_1) = V^-(y;\ep;\beta_1,x_1) + \ep s^-(y;\ep;\beta_1,x_1).
\end{array} \right. 
\end{equation*}
This type of formulation concerning a remainder term $(r^-,s^-)$ is crucial to justify our argument about the existence of a solution satisfying \eqref{b7} (for more information, see \cite{It} and the proof of Lemma \ref{lem1b}).
Consequently, we rewrite \eqref{b7} as the following form with respect to $(r^-,s^-)$: 
\begin{equation}\label{b22_1}
\left \{
\begin{array}{l}
\begin{array}{l}
\ep^2 r^-_{yy} + \ep^2 (h^-_v(V^-_0)s^-)_{yy} + \ep U^-_{yy}(y;\ep;\beta_1,x_1) \\[1ex]
+ \ (x^*(\ep))^2 f(U^- +\ep r^- + \ep h^-_v(V^-_0)s^-, V^- +\ep s^-) /\ep 
 \ = \ 0,   \\[1ex]
D s^-_{yy} + D V^-_{yy}(y;\ep;\beta_1,x_1)/\ep \\[1ex]
- \ (x^*(\ep))^2 f(U^- +\ep r^- + \ep h^-_v(V^-_0)s^-, V^- +\ep s^-) /\ep
\ = \ 0,\\[1ex]
\end{array}  \ y \in (0, 1) \\[1ex]
(r^-,s^-)_y(0) = (0,0), \ (r^-,s^-)(1) = (0,0).
\end{array} 
\right.
\end{equation}
When we simply write \eqref{b22_1} as 
\begin{equation}\label{b22_2}
\begin{array}{c}
T(r^-,s^-;\ep;\beta_1,x_1) \ = \ 0, 
\end{array} 
\end{equation}
$T$ is a smooth mapping from $X_\ep^- := C^2_{\ep,0}[0,1] \times C^2_{1,0}[0,1]$ to $Y := C[0,1] \times C[0,1]$, and then we have the following lemma:
\begin{lemma}\label{lem1b} 
For any given constants $\beta_1^*$ and $x_1^*$, there exist $\ep_0 > 0$, $\rho_0 > 0$, and $K >0$ such that for any 
$\ep \in (0, \ep_0)$ and $(\beta_1,x_1) \in \Delta_{\rho_0} := \{(\beta_1,x_1) \in {\bf R}^2 \ | \ |(\beta_1,x_1) - (\beta^*_1,x^*_1)| \leq \rho_0 \}$, 
\begin{description}
\item[(i)] \  $||T(0,0;\ep;\beta_1,x_1) ||_{Y} = o(1)$ uniformly in $(\beta_1,x_1) \in \Delta_{\rho_0}$ as $\ep \to 0$; 
\item[(ii)] \ for any $(r_1,s_1), (r_2,s_2) \in X_\ep^-$, 
\begin{equation*}
\begin{array}{l}
\displaystyle{ \left|\left| \left( \frac {\partial T}{\partial r^-}, \frac {\partial T}{\partial s^-} \right) 
(r_1,s_1;\ep;\beta_1,x_1) -  \left( \frac {\partial T}{\partial r^-}, \frac {\partial T}{\partial s^-} \right) 
(r_2,s_2;\ep;\beta_1,x_1) \right|\right|_{X_\ep^- \to Y}} \\[0.4cm]
\displaystyle{\hspace*{3.5cm} \ \leq K || (r_1,s_1) - (r_2,s_2)||_{X_\ep^-} }; 
\end{array} 
\end{equation*}
\item[(iii)] \ \hspace*{0.5cm}$ \displaystyle{ \left|\left| \left( \frac {\partial T}{\partial r^-}, \frac {\partial T}{\partial s^-} \right)^{ \! \! -1}  \! \! (0,0;\ep;\beta_1,x_1) \right|\right|_{Y \to X_\ep^-}  \leq K.  }$ 
\end{description}
Moreover, the results (i)-(iii) hold also for $\partial T/\partial \beta_1$ and $\partial T/\partial x_1$ in place of $T$.
\end{lemma}

{\bf Proof}.
Since (i) and (ii) can be proved by the argument similar to that of \cite[Lemma 4.3]{MTH}, we give a proof of (iii). Note that 
\begin{equation*}
\begin{array}{l}
\displaystyle 
\left( \frac {\partial T}{\partial r^-}, \frac {\partial T}{\partial s^-} \right)(0,0;\ep;\beta_1,x_1) \ =: \ \left[
\begin{array}{ll}
T_{11} & T_{12} \\
T_{21} & T_{22} 
\end{array} \right] 
\end{array} 
\end{equation*} 
is represented as 
\begin{equation*}
\begin{array}{lll}
T_{11} & = & \ep^2 \frac {d^2}{d y^2} + (x^*(\ep))^2f_u(U^-,V^-) \\[1ex]
 & = & \ep^2 \frac {d^2}{d y^2} + x_0^2f_u(U_0^- + \phi_0^-,V_0^-) + O(\ep) \\[1ex]
 & =: & T^0_{11} + O(\ep), \\[1ex]
T_{12} & = & \ep^2 \{h_v^-(V_0^-)\frac {d^2}{d y^2} + 2(h_v^-(V_0^-))_y \frac {d}{d y} + (h_v^-(V_0^-))_{yy}\} \\[1ex]
 & & + \ (x^*(\ep))^2\{f_u(U^-,V^-)h^-_v(V_0^-)  + f_v(U^-,V^-)\} \\[1ex]
 & = &  \ep^2 h_v^-(V_0^-)\frac {d^2}{d y^2} + 2\ep^2(h_v^-(V_0^-))_y \frac {d}{d y} \\[1ex] 
 & & + \ x_0^2\{f_u(U_0^- + \phi_0^-,V_0^-)h^-_v(V_0^-) + f_v(U_0^- + \phi_0^-,V_0^-)\} + O(\ep) \\[1ex] 
  & =: & T^0_{12} + O(\ep), \\[1ex]
T_{21} & = & - \ (x^*(\ep))^2f_u(U^-,V^-) 
 =  -  x_0^2f_u(U_0^- + \phi_0^-,V_0^-) + O(\ep) \\[1ex]
  & =: & T^0_{21} + O(\ep), \\[1ex]
T_{22} & = & D \frac {d^2}{d y^2} - (x^*(\ep))^2\{ f_u(U^-,V^-) h^-_v(V_0^-)  + f_v(U^-,V^-)\} \\[1ex]
 & = & D \frac {d^2}{d y^2} - x_0^2\{f_u(U_0^- + \phi_0^-,V_0^-)h^-_v(V_0^-) + f_v(U_0^- + \phi_0^-,V_0^-)\} + O(\ep) \\[1ex]
  & =: & T^0_{22} + O(\ep).
\end{array} 
\end{equation*}
To prove (iii), it suffices to show that for any $F = {}^t(F_1, F_2)\in Y$, there uniquely exists $w \in X_{\ep}^-$ satisfying $T w = F$ such that 
\begin{equation} \label{sss0}
\begin{array}{l}
 ||w||_{X_{\ep}^-} \ \leq \ K ||F||_Y  
\end{array} 
\end{equation}
holds for some constant $K$. By \cite{F}, we easily find that $T^0_{11}$ has a uniformly bounded inverse in $\ep,\beta_1$ and $x_1$. To show the invertibility of $T^0_{22}$, we divide $T^0_{22} = T^0_1 + T^0_2$, where 
$ T^0_1 := D \frac {d^2}{d y^2}$ and $T^0_2 := - x_0^2\{f_u(U_0^- + \phi_0^-,V_0^-)h^-_v(V_0^-) + f_v(U_0^- + \phi_0^-,V_0^-)\}$. 
We see that a solution $g \in C^2_{1,0}[0,1]$ satisfying $T^0_1 g = G$ for any $G \in C[0,1]$ is uniquely represented by 
$$  g(y) = - \ \frac 1D \int_y^1 \int _0^x G(z)dzdx,  $$
which implies that $||g||_{C^2_{1,0}[0,1]} \leq K_1 ||G||_{C[0,1]}$ for some $K_1 > 0$. We then obtain $||(T^0_{1})^{-1}||_{C[0,1] \to C^2_{1,0}[0,1]} \leq K_1$.
On the other hand, for $g \in C^2_{1,0}[0,1]$, we have 
$$ 
|| T^0_2 g ||_{C[0,1]} \leq K_2 ||\phi^-_0 g||_{C[0,1]} 
$$
for some $K_2 >0$.
By \cite[Lemma 4.3]{It}, we find that there exists $K_3 > 0$ satisfying 
\begin{equation} \label{sss1}
\begin{array}{l}
||\phi^-_0 g||_{C[0,1]} \leq \ep K_3||g||_{C^2_{1,0}[0,1]}, 
\end{array} 
\end{equation}
from which we obtain $||T^0_2||_{C^2_{1,0}[0,1] \to C[0,1]} \leq \ep K_4$ for some $K_4>0$. Thus, $T^0_{22}$ is also invertible uniformly in $\ep,\beta_1$ and $x_1$. Therefore $T w = F$ is reduced to 
\begin{equation*} \left\{ 
\begin{array}{l}
  r = - (T^0_{11})^{-1} T^0_{12} s + (T^0_{11})^{-1} F_1,  \\[1ex]
  s = - (T^0_{22})^{-1} T^0_{21} r + (T^0_{22})^{-1} F_2
\end{array} \right. 
\end{equation*}
or 
$$ r =  (T^0_{11})^{-1} T^0_{12}(T^0_{22})^{-1} T^0_{21} r - (T^0_{11})^{-1}(
T^0_{12}(T^0_{22})^{-1} F_2 - F_1). $$
It is easy to see that $||T^0_{21}||_{C^2_{\ep,0}[0,1] \to C[0,1]} \leq K_5$ and $||T^0_{12}||_{C^2_{1,0}[0,1] \to C[0,1]} \leq \ep K_5$, where we used the inequality \eqref{sss1} for some $K_5 > 0$. If we choose $\ep_1 > 0$ to be sufficiently small, it holds that for any $\ep \in (0,\ep_1)$, 
$$    || (T^0_{11})^{-1} T^0_{12}(T^0_{22})^{-1} T^0_{21}||_{C^2_{\ep,0}[0,1] \to C^2_{\ep,0}[0,1]} < 1. $$ 
Therefore \eqref{sss0} holds, and the proof of (iii) is completed. \Qed

Owing to Lemma \ref{lem1b}, we can apply the implicit function theorem to \eqref{b22_2}, and hence we obtain the following: 
\begin{prop}\label{prop1b} 
There exist $\ep_1 > 0$ and $\rho_1 > 0$ such that for any $\ep \in (0,\ep_1)$ and $\rho \in \Delta_{\rho_1}$, there exists $(r^-,s^-)(y;\ep;\beta_1,x_1) \in X_\ep^-$ satisfying 
$$   
T((r^-,s^-)(y;\ep;\beta_1,x_1)) \ = \ 0. 
$$
$(r^-,s^-)(y;\ep;\beta_1,x_1)$, $(\partial (r^-,s^-) \! / \partial \beta_1)(y;\ep;\beta_1,x_1)$ and $(\partial (r^-,s^-) \! / \partial x_1)(y;\ep;\beta_1,x_1)$ are uniformly continuous with respect to $(\ep,\beta_1,x_1) \in (0, \ep_1) \times \Delta_{\rho_1}$ in $X_\ep^-$-topology and satisfy 
\begin{equation*}
\left. 
\begin{array}{l}
 ||(r^-, s^-)(y;\ep;\beta_1,x_1) ||_{X_\ep^-} \\[2ex]
\displaystyle{ \left|\left| \frac {\partial (r^-,s^-)}{\partial \beta_1}(y;\ep;\beta_1,x_1) \right|\right|_{X_\ep^-}} \\[3ex]
\displaystyle{ \left|\left| \frac {\partial (r^-,s^-)}{\partial x_1}(y;\ep;\beta_1,x_1) \right|\right|_{X_\ep^-}}
\end{array} 
\right\}
= o(1) \  \mbox{uniformly in} \ (\beta_1,x_1) \in \Delta_{\rho_1} \ \mbox{as} \ \ep \to 0. 
\end{equation*}
\end{prop}

Thus, we obtain the solution of \eqref{b5} on $[0, x^*(\ep)]$, which takes the form
\begin{equation}\label{b34} \left\{ 
\begin{array}{l}
u^-(x ; \ep;\beta_1,x_1)  \ :=  \tilde{u}^-(\frac{x}{x^*(\ep)};\beta_1,x_1) \\[1ex]
\hspace*{1cm}  = U^-_0( \frac{x}{x^*(\ep)}) + \ep U^-_1(\frac{x}{x^*(\ep)};\beta_1) \\[1ex]
\hspace{1.2cm} + \ \theta(\frac{x}{x^*(\ep)}) \{\phi^-_0( \frac{x-x^*(\ep)}{\ep x^*(\ep)}) 
+  \ep \phi^-_1( \frac{x-x^*(\ep)}{\ep x^*(\ep)};\beta_1,x_1)\} \\[1ex]
\hspace{1.2cm} + \ \ep r^-(\frac{x}{x^*(\ep)};\ep;\beta_1,x_1 ) \\[1ex]
\hspace{1.2cm} + \ \ep h^-_v(V^-_0(\frac{x}{x^*(\ep)})) s^-(\frac{x}{x^*(\ep)};\ep;\beta_1,x_1 ), \\[1ex]
v^-(x ; \ep;\beta_1,x_1)  \ :=  \tilde{v}^-(\frac{x}{x^*(\ep)};\beta_1,x_1) \\[1ex]
\hspace*{1cm}  = V^-_0( \frac{x}{x^*(\ep)}) + \ep V^-_1(\frac{x}{x^*(\ep)};\beta_1) \\[1ex]
\hspace{1.2cm} + \ \theta(\frac{x}{x^*(\ep)}) \{\ep^2 \psi^-_0( \frac{x-x^*(\ep)}{\ep x^*(\ep)}) 
+  \ep^3 \psi^-_1(   \frac{x-x^*(\ep)}{\ep x^*(\ep)};\beta_1,x_1) \\[1ex]
\hspace{1.2cm} - \ep^2 \psi^-_0(0) - \ep^3 \psi^-_1(0;\beta_1,x_1 )  \} \\[1ex]
\hspace{1.2cm} + \ \ep s^-(\frac{x}{x^*(\ep)};\ep;\beta_1,x_1 ), \\[1ex]
\hspace*{6cm} \ x \in [0, x^*(\ep)]. 
\end{array} \right. 
\end{equation}

Next, we consider the solution of \eqref{b6}.
Applying the change of variables $x = x^*(\ep) + (1-x^*(\ep))y$, we have
\begin{equation}\label{b23}
\left \{
\begin{array}{l}
\begin{array}{l}
\ep^2 u_{yy} + (1- x^*(\ep))^2 f(u,v) = 0, \\[1ex]
D v_{yy} - (1- x^*(\ep))^2 f(u,v) = 0,
\end{array}
 \ y \in (0, 1)   \\[0.4cm]
(u,v)(0) = (\alpha,\beta(\ep)), \ (u_y,v_y)(1) = (0,0), \ 
\end{array} 
\right.
\end{equation}
where $\beta(\ep)$ and $x^*(\ep)$ are given by \eqref{b8} and \eqref{b9}, respectively.
Applying the same lines of argument as applied to \eqref{b7}, 
we can obtain the outer approximation of \eqref{b23} 
\begin{equation*}%\label{b24} 
\left\{ 
\begin{array}{l}
u(y)  =  U^+_0(y) + \ep U^+_1(y) + O(\ep^2),\\[1ex]
v(y)  =  V^+_0(y) + \ep V^+_1(y) + O(\ep^2),
\end{array} 
\right.
\end{equation*}
where
\begin{equation}\label{b25}
V^+_0(y) = \beta_0, \ U^+_0(y) = h^{+}(\beta_0),
\end{equation}

\begin{equation}\label{b26}
V^+_1(y;\beta_1) = \beta_1, \ U^+_1(y;\beta_1) = - (f_v^+ / f_u^+) \beta_1 = h^+_v(\beta_0) \beta_1,
\end{equation}
and $f_u^+ := f_u( U^+_0, V^+_0)$ and $f_v^+ := f_v( U^+_0, V^+_0)$.
We note that $(U^+_j,V^+_j )$ $(j=0, 1)$ are constants independent of $y$.

Similarly to the problem \eqref{b7}, these outer approximations do not satisfy the boundary condition at $y = 0$. Hence, 
we must consider the correction of the above approximation in a neighborhood of $y = 0$ 
with the aid of the inner approximation given by 
\begin{equation}\label{b28} \left\{ 
\begin{array}{l}
u(y)  =  U^+_0(y) + \ep U^+_1(y;\beta_1) + \ \phi^+_0( \frac{y}{\ep}) 
+  \ep \phi^+_1(  \frac{y}{\ep} ) + O(\ep^2),\\[1ex]
v(y)  =  V^+_0(y) + \ep V^+_1(y;\beta_1) + \  \ep^2 \psi^+_0( \frac{y}{\ep}) 
+  \ep^3 \psi^+_1(  \frac{y}{\ep} ) + O(\ep^4).
\end{array} \right.
\end{equation}
Introducing the stretched coordinate $z = y/\ep$, and substituting 
\eqref{b28} into \eqref{b23}, and 
comparing the coefficients of powers of $\ep$, we have 
\begin{equation}\label{b29}
\left \{
\begin{array}{l}
\begin{array}{l}
\ddot{\phi}_0^+ + (1-x_0)^2 \tilde{f}^+ = 0, \\[1ex]
D \ddot{\psi}_0^+ + \ddot{\phi}_0^+ = 0,
\end{array}
 \ z \in (0, \infty) \\[0.4cm]
\phi^+_0(0) = \alpha - U^+_0(0), \  \phi^+_0(\infty) = 0, \\[1ex]
\psi^+_0(\infty) = 0, \  \dot{\psi}^+_0(\infty) = 0,
\end{array} 
\right.
\end{equation}
and 
\begin{equation}\label{b30}
\left \{
\begin{array}{l}
\begin{array}{l}
\ddot{\phi}_1^+ + (1-x_0)^2 \tilde{f}^+_u \phi^+_1 = F_1^+(z;\beta_1,x_1), \\[1ex]
D \ddot{\psi}_1^+ + \ddot{\phi}_1^+ = 0, 
\end{array}
 \ z \in (0, \infty) \\[0.4cm]
 \phi^+_1(0) = -U^+_1(0;\beta_1), \ \phi^+_1(\infty) = 0, \\[1ex]
 \psi^+_1(\infty) = 0, \ \dot{\psi}^+_1(\infty) = 0,
\end{array} 
\right.
\end{equation}
where the dot notation denotes $d/dz$, and 
$$
F_1^+(z;\beta_1,x_1) :=  2(1-x_0) x_1 \tilde{f}^+ - (1- x_0)^2 \tilde{f}^+_u U^+_1(0;\beta_1) 
-(1-x_0)^2 \tilde{f}^+_v V_1^+(0;\beta_1) ,
$$
$$ \tilde{f}^+ := f( h^+(\beta_0) + \phi^+_0, \beta_0), \ 
\tilde{f}^+_u := f_u( h^+(\beta_0) + \phi^+_0, \beta_0 ),
$$
and $\tilde{f}^+_v$ is similarly defined.
From Assumption (A1), we find that \eqref{b29} has a unique monotone 
increasing solution $\phi^+_0(z)$ if $ \underline{\alpha}(\beta_0) \leq \alpha < h^+(\beta_0)$, 
where $\underline{\alpha}(\beta_0)$ is defined by 
\begin{equation}\label{Ham2}
\displaystyle \int^{h^+(\beta_0)}_{\underline{\alpha}(\beta_0)} f(u,\beta_0) du = 0.
\end{equation}
We then have $\psi^+_0(z) = - \phi^+_0(z)/D$. 
Moreover, we see that the solution of \eqref{b30} is explicitly given by
\begin{equation}\label{b32} \left\{ 
\begin{array}{l}
\phi^+_1(z;\beta_1,x_1) = -U^+_1(0;\beta_1) \dis\frac{ \dot{\phi}_0^+(z) }{ \dot{\phi}_0^+(0) } \\[0.3cm]
\hspace*{2cm} - \ \dot{\phi}_0^+(z)  \dis\int_0^z \frac{1}{ (\dot{\phi}_0^+(\eta))^2} 
\dis\int_{\eta}^{\infty}  \dot{\phi}_0^+(\zeta) 
F_1^+(\zeta;\beta_1,x_1) d\zeta d\eta, \\[3ex]
\psi^+_1(z;\beta_1,x_1) = - \phi^+_1(z;\beta_1,x_1)/D. 
\end{array} \right. 
\end{equation}

Applying a similar argument to obtain the solution 
$(\tilde{u}^-, \tilde{v}^-)(y;\ep;\beta_1,x_1)$ of \eqref{b7}, 
we can find the solution 
$(\tilde{u}^+, \tilde{v}^+)(y;\ep;\beta_1,x_1)$ of \eqref{b23} as follows:
\begin{equation*} \left\{
\begin{array}{l}
 \tilde{u}^+(y;\ep;\beta_1,x_1)  =  U^+(y;\ep;\beta_1,x_1) + \ep r^+(y;\ep;\beta_1,x_1) + \ep h^+_v(V^+_0(y))s^+(y;\ep;\beta_1,x_1), \\[1ex]
 \tilde{v}^+(y;\ep;\beta_1,x_1)  =  V^+(y;\ep;\beta_1,x_1) + \ep s^+(y;\ep;\beta_1,x_1), 
\end{array} \right. 
\end{equation*}
where
\begin{equation*} \left\{
\begin{array}{l}
U^+(y;\ep;\beta_1,x_1) =  U^+_0(y) + \ep U^+_1(y;\beta_1)  \hspace*{1cm} \\[0.2cm]
\hspace*{1cm}+ \ \theta(1-y)\{\phi^+_0( \frac{y}{\ep}) 
+  \ep \phi^+_1(  \frac{y}{\ep};\beta_1,x_1 ) \}, \\[0.2cm]
V^+(y;\ep;\beta_1,x_1) =  V^+_0(y) + \ep V^+_1(y;\beta_1)  \hspace*{1cm} \\[0.2cm]
\hspace*{1cm}+ \  \theta(1-y)\{\ep^2 \psi^+_0( \frac{y}{\ep}) 
+  \ep^3 \psi^+_1(  \frac{y}{\ep};\beta_1,x_1 ) \\[0.2cm]
\hspace*{1cm}- \ \ep^2 \psi^+_0(0) -  \ep^3\psi^+_1(0;\beta_1,x_1 )\}.
\end{array} \right.
\end{equation*}
Here $(r^+,s^+)(y;\ep;\beta_1,x_1) \in X_\ep^+$ satisfies 
\begin{equation*}
\left. 
\begin{array}{l}
 ||(r^+,s^+)(y;\ep;\beta_1,x_1) ||_{X_\ep^+} \\[2ex]
\displaystyle{ \left|\left| \frac {\partial (r^+,s^+)}{\partial \beta_1}(y;\ep;\beta_1,x_1) \right|\right|_{X_\ep^+}} \\[3ex]
\displaystyle{ \left|\left| \frac {\partial (r^+,s^+)}{\partial x_1}(y;\ep;\beta_1,x_1) \right|\right|_{X_\ep^+ }}
\end{array} 
\right\}
= o(1) \  \mbox{uniformly in} \ (\beta_1,x_1) \in \Delta_{\rho_2} \ \mbox{as} \ \ep \to 0. 
\end{equation*}
$(r^+,s^+)(y;\ep;\beta_1,x_1)$, $(\partial (r^+,s^+) \! /\partial \beta_1)(y;\ep;\beta_1,x_1)$ and $(\partial (r^+,s^+) \! /\partial x_1)(y;\ep;\beta_1,x_1)$ are uniformly continuous with respect to $(\ep,\beta_1,x_1) \in (0, \ep_2) \times \Delta_{\rho_2}$
in $X_\ep^+$-topology, 
where $X_\ep^+ := C^2_{\ep,1}[0,1] \times C^2_{1,1}[0,1]$ and, $\ep_2$ and $\rho_2$ are positive constants. 
Thus, we obtain the solutions of \eqref{b6} on $[x^*(\ep), 1]$ which takes the form 
\begin{equation}\label{b35} \left\{ 
\begin{array}{l}
\! \! \! u^+(x ; \ep;\beta_1,x_1)  \ :=  \tilde{u}^+(\frac{x-x^*(\ep)}{1- x^*(\ep)};\beta_1,x_1) \\[1ex]
\hspace*{1cm}  = 
U^+_0( \frac{x-x^*(\ep)}{1- x^*(\ep)}) + \ep U^+_1(\frac{x-x^*(\ep)}{1- x^*(\ep)};\beta_1) \\[1ex]
\hspace{1.2cm} + \ \theta(\frac{1-x}{1- x^*(\ep)}) \{\phi^+_0( \frac{x-x^*(\ep)}{\ep(1- x^*(\ep))}) 
+  \ep \phi^+_1(  \frac{x-x^*(\ep)}{\ep(1- x^*(\ep))};\beta_1,x_1)\} \\[1ex]
\hspace{1.2cm} + \ \ep r^+( \frac{x-x^*(\ep)}{1- x^*(\ep)};\ep,\beta_1,x_1) \\[1ex]
\hspace{1.2cm} + \ \ep h^+_v(V^+_0( \frac{x-x^*(\ep)}{1- x^*(\ep)}))s^+( \frac{x-x^*(\ep)}{1- x^*(\ep)};\ep,\beta_1,x_1), \\[1ex]
\! \! \! v^+(x ; \ep;\beta_1,x_1)  \ :=  \tilde{v}^+(\frac{x-x^*(\ep)}{1- x^*(\ep)};\beta_1,x_1) \\[1ex]
\hspace*{1cm}  = 
V^+_0( \frac{x-x^*(\ep)}{1- x^*(\ep)}) + \ep V^+_1(\frac{x-x^*(\ep)}{1- x^*(\ep)};\beta_1) \\[1ex]
\hspace{1.2cm} + \ \theta(\frac{1-x}{1- x^*(\ep)})\{\ep^2  \psi^+_0( \frac{x-x^*(\ep)}{\ep(1- x^*(\ep))}) 
+  \ep^3 \psi^+_1(  \frac{x-x^*(\ep)}{\ep(1- x^*(\ep))};\beta_1,x_1) \\[1ex]
\hspace{1.2cm}  - \ \ep^2 \psi^+_0(0) - \ep^3 \psi^+_1(0;\beta_1,x_1 )
\}  \\[1ex]
\hspace{1.2cm} + \ \ep s^+( \frac{x-x^*(\ep)}{1- x^*(\ep)};\ep,\beta_1,x_1), \\[1ex]
\hspace*{6cm} x \in [x^*(\ep), 1]. 
\end{array} \right.
\end{equation}

\subsection{$C^1$-matching at $x = x^*(\ep)$}\label{S2.2} %Subsection 2.2 

We now consider the $C^1$-matching of $(u^-,v^-)(x ; \ep;\beta_1,x_1)$ and $(u^+,v^+)(x ; \ep;\beta_1,x_1)$ at $x = x^*(\ep)$ to obtain the approximation of the solution of \eqref{b1} up to $O(\ep)$. Since these two solutions are already continuous at $x = x^*(\ep)$ and we recall that $\ep^2 u_x + D v_x = 0$ holds by \eqref{b1}, we can impose this condition only on the $u$-component. Then 
we determine the values of $\beta_j$ and $x_j$ $(j=0, 1)$ in such a way that
\[
\Phi(\ep) := \ep x^*(\ep) ( 1 - x^*(\ep)) \{
\frac{d}{dx}u^-(x^*(\ep) ; \ep;\beta_1,x_1) - \frac{d}{dx}u^+(x^*(\ep) ; \ep;\beta_1,x_1) \} = o(\ep) 
\]
holds for small $\ep >0$. Noting that $U^-_j$ and $U^+_j$ are constants,  
it follows from \eqref{b34} and \eqref{b35} that
\[
\begin{array}{lcl}
\Phi(\ep) & = & \ep  ( 1 - x^*(\ep)) \{
\dot{\phi}^-_0(0)/\ep + \dot{\phi}^-_1(0;\beta_1,x_1) + o(1) \} \\[1ex]  
 & & - \ \ep  x^*(\ep) \{
\dot{\phi}^+_0(0)/\ep + \dot{\phi}^+_1(0;\beta_1,x_1) + o(1) \} \\[1ex] 
& = & \{ (1-x_0) \dot{\phi}^-_0(0) - x_0  \dot{\phi}^+_0(0)\} + \ep \{ (1-x_0) \dot{\phi}^-_1(0;\beta_1,x_1)  \\[1ex] 
& &  - \ x_1  \dot{\phi}^-_0(0) - x_0  \dot{\phi}^+_1(0;\beta_1,x_1)
- x_1  \dot{\phi}^+_0(0) \} + o(\ep)\\[1ex] 
&  =: & \Phi_0 + \ep \Phi_1(\beta_1,x_1) + o(\ep).
\end{array}
\]
First, we consider 
\begin{equation} \label{b00009}
\begin{array}{l}
\Phi_0 \ = \ (1-x_0) \dot{\phi}^-_0(0) - x_0  \dot{\phi}^+_0(0) = 0.
\end{array}
\end{equation}
It follows from \eqref{b18} that 
\[
\begin{array}{rcl}
0 & = & \dis\int_{-\infty}^0 \{ \ddot{\phi}^-_0 \dot{\phi}^-_0
+ x_0^2 f( h^-(\beta_0) + \phi^-_0, \beta_0 ) \dot{\phi}^-_0 \} dz 
\\[2ex]
& = & 
\dis\frac{(\dot{\phi}^-_0(0))^2}{2} + x_0^2
\dis\int_{h^-(\beta_0)}^\alpha f(u, \beta_0) du,
\end{array}
\]
which implies 
\[
\dot{\phi}^-_0(0) = x_0 \sqrt{-2 \dis\int_{h^-(\beta_0)}^\alpha f(u,\beta_0) du}. 
\]
Similarly, we see from \eqref{b29} that 
\[
\dot{\phi}^+_0(0) = (1-x_0) \sqrt{ 2 \dis\int_\alpha^{h^+(\beta_0)} f(u,\beta_0) du}. 
\]
Therefore, we have
\[
\Phi_0 = 
-\dis\frac{ 2 x_0(1-x_0) J(\beta_0)}
{ \sqrt{- 2 \dis\int_{h^-(\beta_0)}^\alpha f(u,\beta_0) du}
+ \sqrt{  2 \dis\int_\alpha^{h^+(\beta_0)} f(u,\beta_0) du} },
\]
where $J=J(v)$ is given by \eqref{a4}. Hence, noting Assumption (A2),
it follows from \eqref{b14}, \eqref{b25} and $\Phi_0 = 0$ that
\begin{equation}\label{b35x}
V^{\pm}_0(y) = \beta_0 = v^* \ \ \ \text{and} \ \ \ U^{\pm}_0(y) = h^{\pm}(v^*).
\end{equation}
Moreover, we can take $\alpha$ satisfying $h^-(v^*) < \alpha < h^+(v^*)$
because $\underline{\alpha}(v^*) = h^-(v^*)$
and  $\bar{\alpha}(v^*) = h^+(v^*) $ by \eqref{Ham1}, \eqref{Ham2}
and Assumption (A2).
We note that though $\Phi_0$ depends on both $\beta_0$ and $x_0$, the solution satisfying $\Phi_0 = 0$ is determined by only $\beta_0 = v^*$ for any $x_0 \in (0,1)$.
Moreover, we have
\begin{equation}\label{b36}
\dot{\phi}^-_0(0) = x_0 \sqrt{-2 \dis\int_{h^-(v^*)}^\alpha f(u, v^* ) du}
\end{equation}
and
\begin{equation}\label{b37}
\dot{\phi}^+_0(0) = (1-x_0) \sqrt{ 2 \dis\int_\alpha^{h^+(v^*)} f(u, v^* ) du}. 
\end{equation}

Next, we consider
$\Phi_1(\beta_1,x_1) = (1-x_0) \dot{\phi}^-_1(0;\beta_1,x_1) - x_1  \dot{\phi}^-_0(0) - x_0  \dot{\phi}^+_1(0;\beta_1,x_1) - x_1  \dot{\phi}^+_0(0) = 0$.
Note the following relation: 
\[
\dis\int_{-\infty}^0 x_0^2 \tilde{f}^-_u  \dot{\phi}^-_0 dz
=
 - \dis\int_{-\infty}^0 \dddot{\phi}^-_0  dz
= 
-  \ddot{\phi}^-_0 (0), 
\]
where we used the relation $\dddot{\phi}^-_0 + x_0^2 \tilde{f}_u^- \dot{\phi}_0^- = 0$ obtained by
the differentiation of the first equation of \eqref{b18} in $z$. 
Then, we have from \eqref{b21} that
\begin{equation}\label{b38}
\begin{array}{l}
 \dot{\phi}^-_1(0;\beta_1,x_1) 
= \dis  -U^-_1(1;\beta_1)\frac {\ddot{\phi}^-_0 (0)}{\dot{\phi}_0^-(0)} + 
 \frac{1}{\dot{\phi}_0^-(0)}\dis\int_{-\infty}^0 F_1^-(z;\beta_1,x_1) \dot{\phi}^-_0 dz
\\[3ex]
\ \ =
\dis\frac{1}{\dot{\phi}_0^-(0)} \left(   -U^-_1(1;\beta_1)\ddot{\phi}^-_0 (0) 
-2x_0 x_1 \dis\int_{-\infty}^0  \tilde{f}^- \dot{\phi}^-_0 dz \right. \\[3ex]
\left. \hspace*{2.5cm} - x_0^2 U^-_1(1;\beta_1) \dis\int_{-\infty}^0  \tilde{f}^-_u \dot{\phi}^-_0 dz
 - x_0^2 \beta_1 \dis\int_{-\infty}^0 \tilde{f}^-_v \dot{\phi}^-_0 dz \right)
\\[3ex]
\ \ = 
\dis\frac{1}{\dot{\phi}_0^-(0)} \left(
-2x_0 x_1 \dis\int_{h^-(v^*)}^\alpha  f(u, v^*) du
- \dis x_0^2 \beta_1 \dis\int_{h^-(v^*)}^\alpha  f_v(u, v^*) du \right). 
\end{array}
\end{equation}

Similarly, we can see from \eqref{b29} and \eqref{b32} that 
\begin{equation}\label{b39}
\begin{array}{l}
\dot{\phi}^+_1(0;\beta_1,x_1)  
= 
\dis\frac{1}{\dot{\phi}_0^+(0)} \left(
-2(1-x_0 ) x_1 \dis\int^{h^+(v^*)}_\alpha  f(u, v^*) du \right.
\\[2ex]
\ \ \ \ \ \ \ \ \ \ \ \ \ \ \ \ \ \ \ \ \ \ \ \ \ \left.
 + \dis (1-x_0)^2 \beta_1 \dis\int^{h^+(v^*)}_\alpha  f_v(u, v^*) du \right).
\end{array}
\end{equation}
Hence, it follows from \eqref{b38} and \eqref{b39} that
\[
\begin{array}{rcl}
\Phi_1(\beta_1,x_1) & = & (1-x_0) \dot{\phi}^-_1(0;\beta_1,x_1) - x_1  \dot{\phi}^-_0(0) - x_0  \dot{\phi}^+_1(0;\beta_1,x_1) - x_1  \dot{\phi}^+_0(0)
\\[2ex]
& =: & K(x_0)x_1 + M(x_0)\beta_1 + R(x_0),
\end{array}
\]
where
\begin{equation}\label{b40}
\begin{array}{l}
K(x_0) = \dis\frac{x_0(1-x_0)}{\dot{\phi}^-_0(0)}
\left( -2 \dis\int_{h^-(v^*)}^\alpha f(u, v^* ) du \right) - \dot{\phi}^-_0(0)
\\[3ex]
\ \ \ \ \ \ \ \ \ \ \ 
+ \dis\frac{x_0(1-x_0)}{\dot{\phi}^+_0(0)}
\left( 2 \dis\int_\alpha^{h^+(v^*)} f(u, v^* ) du \right) - \dot{\phi}^+_0(0),
\end{array}
\end{equation}
\begin{equation}\label{b41}
\begin{array}{l}
M(x_0) = - \dis\frac{x_0^2(1-x_0)}{\dot{\phi}^-_0(0)}
\left(  \dis\int_{h^-(v^*)}^\alpha f_v(u, v^* ) du \right) 
\\[3ex]
\ \ \ \  \ \ \ \ \ \ \ \ \ \ 
- \dis\frac{x_0(1-x_0)^2}{\dot{\phi}^+_0(0)}
\left( \dis\int_\alpha^{h^+(v^*)} f_v(u, v^* ) du \right),
\end{array}
\end{equation}
and $R(x_0) = 0$. Moreover, it follows from \eqref{b36}, \eqref{b37} and 
\eqref{b40} that
\[
\begin{array}{rcl}
K(x_0) & = & (1-2x_0) \left( \sqrt{ -2 \dis\int_{h^-(v^*)}^\alpha f(u, v^* ) du } 
- \sqrt{ 2 \dis\int_\alpha^{h^+(v^*)} f(u, v^* ) du } \right)
\\[4ex]
& = & \dis\frac{ -2(1-2x_0) J(v^*)}{ \sqrt{ -2 \dis\int_{h^-(v^*)}^\alpha f(u, v^* ) du }
+ \sqrt{ 2 \dis\int_\alpha^{h^+(v^*)} f(u, v^* ) du } } = 0,
\end{array}
\]
where $J(v)$ is given by \eqref{a4}. 
This implies that $\Phi_1$ does not depend on $x_1$.
Noting
\[
\dis\int_\alpha^{h^+(v^*)} f(u, v^* ) du = - \dis\int_{h^-(v^*)}^\alpha f(u, v^* ) du
\]
by $J(v^*) = 0$, 
it follows from \eqref{b36}, \eqref{b37} and 
\eqref{b41} that
\[
M(x_0)  = 
- \dis\frac{x_0(1-x_0)}{ \sqrt{ -2 \dis\int_{h^-(v^*)}^\alpha f(u, v^* ) du } }
\dis\int_{h^-(v^*)}^{h^+(v^*)} f_v(u, v^* ) du \neq 0 
\]
by \eqref{a4x}. Therefore, from \eqref{b15}, \eqref{b26} and $\Phi_1 = 0$, we have 
\begin{equation}\label{b43x}
V^{\pm}_1(y ; \beta_1) = \beta_1 = -R(x_0)/M(x_0) = 0 \ \ \ \ \text{and} \ \ \ \  U^{\pm}_1(y ; \beta_1) = 0.
\end{equation}

\subsection{Computation of \eqref{b1_1}}\label{S2.3} %Subsection 2.3 

To complete the construction of the approximate solution of \eqref{b1} with \eqref{b1_1}, 
we determine the values of $\beta_j$ and $x_j$ $(j=0, 1)$ by the conservation law \eqref{b1_1}. 
Although we know that $\beta_1 = 0$ by \eqref{b43x}, we purposely give expressions including
$\beta_1$ in the following calculations, which are helpful to understand an argument
in the next subsection.

When we put 
\begin{equation}\label{b44}
\Psi(\ep) := \dis \int_0^1 \{u(x;\ep;\beta_1,x_1) + v(x;\ep;\beta_1,x_1)\}dx - \xi,
\end{equation}
\eqref{b1_1} is equivalent to $\Psi(\ep) = 0$. 
Using \eqref{b14}, \eqref{b15}, \eqref{b34}, \eqref{b25}, \eqref{b26}, \eqref{b35} and \eqref{b35x}, 
we have
\begin{equation}\label{b45} 
\begin{array}{l}
\dis \int_0^1 u(x;\ep;\beta_1,x_1)dx  = \int_0^{x^*(\ep)} u^-(x;\ep;\beta_1,x_1)dx  
+ \int_{x^*(\ep)}^1 u^+(x;\ep;\beta_1,x_1)dx, \\[3ex]
\dis \int_0^1 v(x;\ep;\beta_1,x_1)dx  = \int_0^{x^*(\ep)} v^-(x;\ep;\beta_1,x_1)dx  
+ \int_{x^*(\ep)}^1 v^+(x;\ep;\beta_1,x_1)dx, 
\end{array} 
\end{equation}
where
\[
\begin{array}{l}
\dis\int_0^{x^*(\ep)} u^-(x;\ep;\beta_1,x_1)dx  = x^*(\ep) \{
\dis\int_0^1 ( U^-_0(y) +\ep  U^-_1(y;\beta_1)) dy 
\\[3ex]
\ \ \ \ \ \ \ \ \ \ \ \ \ \ \ 
+ \, \ep \dis\int_{-\infty}^0  {\phi}_0^-(z)dz + o(\ep)
\, \} 
\\[3ex]
\ \ \ \ \ 
= ( x_0 + \ep x_1 + o(\ep) )\{ h^-(v^*) + \ep h^-_v(v^*) \beta_1 \\[3ex]
\ \ \ \ \ \ \ \ \ \ \ \ \ \ 
+ \, \ep \dis\int_{-\infty}^0   {\phi}_0^-(z) dz + o(\ep) \}, 
\end{array}
\]
\[
\begin{array}{l}
\dis\int_{x^*(\ep)}^1 u^+(x;\ep;\beta_1,x_1)dx  = (1-x^*(\ep)) \{
\dis\int_0^1 ( U^+_0(y) +\ep  U^+_1(y;\beta_1)) dy 
\\[3ex]
\ \ \ \ \ \ \ \ \ \ \ \ \ \ \ 
+ \, \ep \dis\int^{\infty}_0  {\phi}_0^+(z) dz + o(\ep)
\, \} 
\\[3ex]
\ \ \ \ \ 
= (1- x_0 - \ep x_1+o(\ep))\{ h^+(v^*) + \ep h^+_v(v^*) \beta_1 \\[3ex]
\ \ \ \ \ \ \ \ \ \ \ \ \ \
+ \, \ep \dis\int^{\infty}_0   {\phi}_0^+(z) dz + o(\ep) \}, 
\end{array}
\]
\[
\begin{array}{l}
\dis\int_0^{x^*(\ep)} v^-(x;\ep;\beta_1,x_1)dx  = x^*(\ep) \{
\dis\int_0^1 ( V^-_0(y) +\ep  V^-_1(y;\beta_1)) dy + o(\ep)  
\, \} 
\\[3ex]
\ \ \ \ \ 
= ( x_0 + \ep x_1 + o(\ep))\{v^* + \ep \beta_1   + o(\ep) \}, 
\end{array}
\]
and
\[
\begin{array}{l}
\dis\int_{x^*(\ep)}^1 v^+(x;\ep;\beta_1,x_1)dx  = (1-x^*(\ep)) \{
\dis\int_0^1 ( V^+_0(y) +\ep  V^+_1(y;\beta_1)) dy + o(\ep) 
\, \} 
\\[3ex]
\ \ \ \ \ 
= (1- x_0 - \ep x_1 + o(\ep) )\{v^* + \ep \beta_1 + o(\ep) \}. 
\end{array}
\]
Substituting \eqref{b45} into \eqref{b44}, we have  
\begin{equation}\label{b45_1}
\begin{array}{l}
0 = \Psi(\ep) = \dis \{v^* + x_0 h^-(v^*) + (1- x_0) h^+(v^*) - \xi \} \\[1ex]
+ \ \ep \left\{ \beta_1 + x_0 h_v^-(v^*) \beta_1 + (1-x_0) h_v^+(v^*) \beta_1 + x_0 \dis\int_{-\infty}^0   {\phi}_0^-(z) dz
\right. \\[2ex]
+ \ \left. 
  x_1 h^-(v^*) +  (1-x_0) \dis\int^{\infty}_0   {\phi}_0^+(z) dz   - x_1 h^+(v^*) \right\} 
+ o(\ep) \\[3ex]
 =: \Psi_0 + \ep \Psi_1(\beta_1,x_1)  + o(\ep). 
\end{array}
\end{equation}
Comparing the each coefficients of powers of $\ep$ in \eqref{b45_1}, we have 
$ \Psi_i = 0 \ (i=0,1)$. 
Noting assumption (A4) and $h^-(v^*) < h^+(v^*)$ by Assumption (A1), 
it follows from $\Psi_0 = 0$ that
\begin{equation}\label{b48_0}
\begin{array}{l}
x_0 = \dis\frac{v^* + h^+(v^*) - \xi}{h^+(v^*) - h^-(v^*) }
\ \ \ \text{and} \ \ \ 
0 < x_0 < 1.
\end{array}
\end{equation}
Moreover, since $\beta_1 = 0$ by \eqref{b43x}, it follows from $\Psi_1(\beta_1,x_1) = 0$ that
\begin{equation}\label{b48_1}
\begin{array}{l}
x_1 = \dis\frac{I_1(x_0)}{h^+(v^*) - h^-(v^*) },
\end{array}
\end{equation}
where 
\[
I_1(x_0) =  x_0 \dis\int_{-\infty}^0   {\phi}_0^-(z) dz 
+ (1-x_0) \dis\int^{\infty}_0   {\phi}_0^+(z) dz
\]
is a function of $x_0$. 
Thus, we see that \eqref{b34} and \eqref{b35} give 
the approximation of the solution of \eqref{b1} with \eqref{b1_1}. 

\subsection{Determination of $\beta(\ep)$ and $x^*(\ep)$}\label{S2.4} % Subsection 2.4 
Finally, we determine  $\beta(\ep)$ and $x^*(\ep)$ uniquely such that \eqref{b1} with \eqref{b1_1} have a 
single transition layer solution $(u, v)(x;\ep)$ at the layer position $ x = x^*(\ep)$.

First, the coefficients $\beta_i$ and $x_i$ $(i=0,1)$ are determined step by step as follows:  $\beta_0$ (by \eqref{b35x}) $\longrightarrow$ 
$x_0$ (by \eqref{b48_0})  $\longrightarrow$ $\beta_1$ (by \eqref{b43x})  $\longrightarrow$ $x_1$ (by \eqref{b48_1}). We note that $\Phi(\ep) = o(\ep)$ and $\Psi(\ep) = o(\ep)$; $\Phi$ and $\Psi$ are not identically zero for these $\beta_i$ and $x_i$ $(i=0,1)$. 

Next, we set $(\beta_1^*, x_1^*) = (\beta_1, x_1)$ in Lemma \ref{lem1b}, and 
consider $\beta(\ep) = \beta_0 + \ep \bar{\beta}_1$ and 
$x(\ep) = x_0 + \ep \bar{x}_1$. 
We can take $(\bar{\beta}_1, \bar{x}_1)$ around $(\beta_1^*, x_1^*)$ 
so as to satisfy $\Phi(\ep) = 0$ and $\Psi(\ep) = 0$ as follows:
Let us define $\Phi^*(\bar{\beta}_1, \bar{x}_1;\ep)$ and  
$\Psi^*(\bar{\beta}_1, \bar{x}_1;\ep)$ by $\Phi(\ep) = \ep \Phi^*(\bar{\beta}_1, \bar{x}_1;\ep)$ and $\Psi(\ep) = \ep \Psi^*(\bar{\beta}_1, \bar{x}_1;\ep)$, respectively. We easily find that there exist two positive constants $\delta$ and $\ep_3 (< \min \{\ep_1, \ep_2\})$ 
such that $\Phi^*(\bar{\beta}_1, \bar{x}_1;\ep)$ and $\Psi^*(\bar{\beta}_1, \bar{x}_1;\ep)$ are continuous in $\bar{\beta}_1 \in (\beta_1^*-\delta, \beta_1^*+\delta), \bar{x}_1 \in (x_1^*-\delta, x_1^*+\delta)$ and $\ep \in [0,\ep_3)$, and are $C^1$-functions of $\beta_1$ and $x_1$. Moreover, we can easily find that 
\[ \left\{ 
\begin{array}{l}
\Phi^*(\beta_1^*, x_1^*;0) = 0, \
\frac {\partial \Phi^*}{\partial \bar{\beta}_1}(\beta_1^*, x_1^*;0) = M(x_0), \
\frac {\partial \Phi^*}{\partial \bar{x}_1}(\beta_1^*, x_1^*;0) = 0, \\[2ex]
\Psi^*(\beta_1^*, x_1^*;0) = 0, \ 
\frac {\partial \Psi^*}{\partial \bar{\beta}_1}(\beta_1^*, x_1^*;0) 
= (1 + x_0 h_v^-(v^*) +  (1 - x_0) h_v^+(v^*))/D,  \\[2ex]
\frac {\partial \Psi^*}{\partial \bar{x}_1}(\beta_1^*, x_1^*;0) = h^+(v^*) - h^-(v^*), 
\end{array} \right.
\]
which implies that 
\[
\begin{array}{l}
 \dfrac {\partial (\Phi^*, \Psi^*)}{\partial (\bar{\beta}_1, \bar{x}_1)}  (\beta_1^*, x_1^* ;0) 
 = M(x_0) (h^+(v^*) - h^-(v^*)) \  \ne \ 0. 
\end{array}
\]
Then, we can apply the implicit function theorem to $\Phi^*(\bar{\beta}_1, \bar{x}_1;\ep) = 0$ and  
$\Psi^*(\bar{\beta}_1, \bar{x}_1;\ep) = 0$, and find that there uniquely exist $\bar{\beta}_1 = \bar{\beta}_1(\ep)$
and $\bar{x}_1 = 
\bar{x}_1(\ep)$ for $\ep \in [0,\ep_3)$ satisfying $\bar{\beta}_1(0) = \beta_1^*, \ \bar{x}_1(0) = 
 x_1^*$,  
$$\Phi^*(\bar{\beta}_1, \bar{x}_1;\ep) = 0, \ \text{and} \ \Psi^*(\bar{\beta}_1, \bar{x}_1;\ep) = 0.$$

Substituting $\beta_1 = \bar{\beta}_1(\ep)$ and $x_1 = \bar{x}_1(\ep)$ into 
\eqref{b8}, \eqref{b9}, \eqref{b34} and \eqref{b35}, we obtain the following existence result:

\begin{theo}\label{th2}
Assume that (A1), (A2) and (A4). 
For any $\ep \in (0, \ep_3)$, there exists a family of single jump-up transition layer solutions $(u,v)(x;\ep) \in \pot{C}^2_\ep[0,1] \times \pot{C}^2_1[0,1]$ of \eqref{b1} with \eqref{b1_1}. Furthermore, the following estimate holds:
\begin{equation*}
\begin{array}{c}
\dis \left|\left| u(x;\ep) - U^-(\frac x{x^*(\ep)};\ep;  \bar{\beta}_1(\ep), \bar{x}_1(\ep)) \right|\right|_
{C^2_{\ep,0}[0,x^*(\ep)]} \\[4ex]
 + \ \dis \left|\left| v(x;\ep) - V^-(\frac x{x^*(\ep)};\ep;  \bar{\beta}_1(\ep), \bar{x}_1(\ep)) \right|\right|_
{C^2_{1,0}[0,x^*(\ep)]} \\[4ex]
+ \ \dis \left|\left| u(x;\ep) - U^+(\frac {x-x^*(\ep)}{1-x^*(\ep)};\ep; \bar{\beta}_1(\ep), \bar{x}_1(\ep)) \right|\right|_{C^2_{\ep,1}[x^*(\ep),1]} \\[4ex]
+ \ \dis \left|\left| v(x;\ep) - V^+(\frac {x-x^*(\ep)}{1-x^*(\ep)};\ep; \bar{\beta}_1(\ep), \bar{x}_1(\ep)) \right|\right|_{C^2_{1,1}[x^*(\ep),1]} = o(\ep)
\end{array}
\end{equation*}
as $\ep \to 0$. 
\end{theo}

\begin{remark}
In Section \ref{S2.1}, we took $U_0^-(y) = h^-(\beta_0)$ in \eqref{b14} and $U_0^+(y) = h^+(\beta_0)$ in \eqref{b25} as a jump-up solution. If we take $U_0^-(y) = h^+(\beta_0)$ in \eqref{b14} and $U_0^+(y) = h^-(\beta_0)$ in \eqref{b25}, 
we can obtain a family of single jump-down transition layer solutions $(u,v)(x;\ep) \in \pot{C}^2_\ep[0,1] \times \pot{C}^2_1[0,1]$ of \eqref{b1} with \eqref{b1_1} such that $u(x;\ep)$ has a jump-down layer at $x = x^*(\ep)$. 
\end{remark}

\begin{remark} \label{rem2.2}
In the next section, the leading terms 
$(U^{\pm}, V^{\pm})(\, \cdot \, ;  \ep; \beta_1, x_1)$ given in Theorem \ref{th2} are very crucial to calculate the eigenvalue problems (\ref{4c3}) and (\ref{4c4}). We then give these explicit forms here. 
\[ \left\{ 
\begin{array}{lcl}
U^-(\frac x{x^*(\ep)};\ep; \beta_1, x_1) & = & h^-(v^*) + \phi^-_0(\frac {x-x^*(\ep)}{\ep x^*(\ep)}) + \ep \phi^-_1(\frac {x-x^*(\ep)}{\ep x^*(\ep)};\beta_1,x_1), \\[0.2cm]
V^-(\frac x{x^*(\ep)};\ep; \beta_1, x_1) & = & v^* + \ep^2 \psi^-_0(\frac {x-x^*(\ep)}{\ep x^*(\ep)}) + \ep^3 \psi^-_1(\frac {x-x^*(\ep)}{\ep x^*(\ep)};\beta_1,x_1), \\[0.2cm]
 & & \hspace{3cm} x \in [0,x^*(\ep)]
\end{array}  \right.
\]
and 
\[ \left\{ 
\begin{array}{lcl}
U^+(\frac {x-x^*(\ep)}{1-x^*(\ep)};\ep; \beta_1, x_1) & = & h^+(v^*) + \phi^+_0(\frac {x-x^*(\ep)}{\ep(1-x^*(\ep))}) + \ep \phi^+_1(\frac {x-x^*(\ep)}{\ep(1-x^*(\ep))};\beta_1,x_1), \\[0.2cm]
V^+(\frac {x-x^*(\ep)}{1-x^*(\ep)};\ep; \beta_1, x_1) & = & v^* + \ep^2 \psi^+_0(\frac {x-x^*(\ep)}{\ep(1-x^*(\ep))}) + \ep^3 \psi^+_1(\frac {x-x^*(\ep)}{\ep(1-x^*(\ep))};\beta_1,x_1), \\[0.2cm]
 & & \hspace{2.7cm} x \in [x^*(\ep),1], 
\end{array}  \right.
\]
where we note that $U^{\pm}_1 = 0$ and $V^{\pm}_1 = 0$ by \eqref{b43x}. 

\end{remark}

\section{Stability analysis of the transition layer solutions} \label{S3} % Section 3

In this section, we perform the stability analysis of the linearization of \eqref{a1}
around $(u,v)(x;\ep)$ given by
Theorem~\ref{th2} under a natural constraint derived from 
the conservation law \eqref{a2}.

We consider the linearized eigenvalue problem of \eqref{a1}
\begin{equation}\label{4c1}
{\cal L}^\ep 
\left[
\begin{array}{c}
p\\
q
\end{array}
\right] :=  
\left[
\begin{array}{cc}
\dis \varepsilon^2\frac{d^2}{dx^2}+f_u^\ep& f_v^\ep\\
 -f_u^\ep&D \dis \frac{d^2}{dx^2}- f_v^\ep
\end{array}
\right] 
\left[
\begin{array}{c}
p\\
q
\end{array}
\right]
=
\lambda \ 
\left[
\begin{array}{c}
p\\
q
\end{array}
\right],
\end{equation}
under the Neumann boundary condition, where 
$f_u^\ep := f_u ( u(x;\ep), v(x;\ep) )$, $f_v^\ep := f_v (u(x;\ep), v(x;\ep))$ and 
$\lambda \in {\mathbb C}$.
The underlying space for \eqref{4c1} can be taken as $BC[0,1]\times BC[0,1]$ with 
$${\cal D}({\cal L}^\ep) := \{ (p, q) \in \pot{C}^2_\ep[0,1] \times \pot{C}^2_1[0,1] \ | \  \int_0^1 (p+q) dx = 0 \} $$
by virtue of \eqref{a2}. 
We note that for $(p, q ) \in \pot{C}^2_\ep[0,1] \times \pot{C}^2_1[0,1]$ satisfying \eqref{4c1},
the condition
\begin{equation*}%\label{4cc1}
\begin{array}{c}
\displaystyle \lambda \int_0^1 (p+q) dx = 0
\end{array}
\end{equation*}
always holds
by integrating the equations with respect to $p$ and $q$ in \eqref{4c1} on the interval $[0,1]$ under the Neumann boundary conditions. This fact implies that $(p, q) \in {\cal D}({\cal L}^\ep)$ 
if $(p, q ) \in \pot{C}^2_\ep[0,1] \times \pot{C}^2_1[0,1]$ satisfies \eqref{4c1}
for $\lambda \neq 0$.
The equation (\ref{4c1}) can be rewritten equivalently as 
\begin{equation} \label{4c2}\left\{ 
\begin{array}{l}
\dis \frac{d}{dx}\bar{V}=A(x;\varepsilon;\lambda)\bar{V}, \ \ x \in (0,1) \\[2ex]
(p_x, q_x)(0) = (0,0), \ \ (p_x,p_x)(1) = (0,0)
\end{array} \right. 
\end{equation}
for $\bar{V}=\bar{V}(x;\varepsilon;\lambda) := \displaystyle{ \left(p,\ep p_x,q,q_x\right)(x;\varepsilon;\lambda)}$, 
where $A(x;\varepsilon;\lambda)$ is defined by 
\begin{eqnarray*}
A(x;\varepsilon;\lambda)
:= 
\left[
\begin{array}{cccc}
0&1/\varepsilon&0&0 \\
(\lambda - f_u^\ep)/\varepsilon&0&- f_v^\ep/\varepsilon &0\\
0&0&0&1\\
f_u^\ep/D&0&(\lambda + f_v^\ep)/D&0\\
\end{array}
\right].
\end{eqnarray*}

In a similar manner to the construction of a transition layer solution, let us consider the following problems with suitable boundary conditions: 
\begin{equation}\label{4c3}
\left \{
\begin{array}{l}
\begin{array}{l}
\ep^2 p_{xx} +  f_u^\ep p + f_v^\ep q = \lambda p, \\[1ex]
Dq_{xx}  - f_u^\ep p - f_v^\ep q = \lambda q, 
\end{array} \quad  x \in (0,x^*(\ep)) \\[0.4cm]
(p_x, q_x)(0) = (0, 0), \ (p, q)(x^*(\ep)) = (a,b) 
\end{array} 
\right.
\end{equation}
and 
\begin{equation}\label{4c4}
\left \{
\begin{array}{l}
\begin{array}{l}
\ep^2 p_{xx} +  f_u^\ep p + f_v^\ep q = \lambda p, \\[1ex]
Dq_{xx}  - f_u^\ep p - f_v^\ep q = \lambda q, 
\end{array} \quad  x \in (x^*(\ep),1) \\[0.4cm]
(p, q)(x^*(\ep)) = (a,b), \ (p_x, q_x)(1) = (0,0),
\end{array} 
\right.
\end{equation}
where $a,b$ are given real numbers. For any $\lambda \in \mathbb{C}$, let 
$(p^-,q^-)(x;\varepsilon;\lambda;a,b)$ and 
$(p^+,q^+)(x;\varepsilon;\lambda;a,b)$ be 
solutions of (\ref{4c3}) and (\ref{4c4}), respectively. Then, any solution 
$\bar{V}(x;\varepsilon;\lambda)$ of  (\ref{4c2}) on $[0,x^*(\ep)]$ is 
represented as a linear combination of two independent solutions 
\begin{equation}\label{v12}
\begin{array}{l}
\bar{V}_1(x;\varepsilon;\lambda)
:=
\left[
\begin{array}{c}
p^-(x;\varepsilon;\lambda;1,0)\\
\varepsilon p_x^-(x;\varepsilon;\lambda;1,0)\\
q^-(x;\varepsilon;\lambda;1,0)\\
q_x^-(x;\varepsilon;\lambda;1,0)
\end{array}
\right], \ 
\bar{V}_2(x;\varepsilon;\lambda)
:=
\left[
\begin{array}{c}
p^-(x;\varepsilon;\lambda;0,1)\\
\varepsilon p_x^-(x;\varepsilon;\lambda;0,1)\\
q^-(x;\varepsilon;\lambda;0,1)\\
q_x^-(x;\varepsilon;\lambda;0,1)
\end{array}
\right].
\end{array}
\end{equation}
Similarly, any solution of  (\ref{4c2}) on $[x^*(\ep),1]$ is 
represented as a linear combination of two independent solutions 
\begin{equation}\label{v34}
\begin{array}{l}
\bar{V}_3(x;\varepsilon;\lambda)
:=
\left[
\begin{array}{c}
p^+(x;\varepsilon;\lambda;1,0)\\
\varepsilon p_x^+(x;\varepsilon;\lambda;1,0)\\
q^+(x;\varepsilon;\lambda;1,0)\\
q_x^+(x;\varepsilon;\lambda;1,0)
\end{array}
\right], \ 
\bar{V}_4(x;\varepsilon;\lambda)
:=
\left[
\begin{array}{c}
p^+(x;\varepsilon;\lambda;0,1)\\
\varepsilon p_x^+(x;\varepsilon;\lambda;0,1)\\
q^+(x;\varepsilon;\lambda;0,1)\\
q_x^+(x;\varepsilon;\lambda;0,1)
\end{array}
\right].
\end{array}
\end{equation}
Since the coefficient matrix $A(x;\varepsilon;\lambda)$ of (\ref{4c2}) depends analytically on $\lambda$, we can consider, without loss of generality, that 
$\bar{V}_i(x;\varepsilon;\lambda)\ (i=1,2,3,4)$ also depend analytically on $\lambda$. 

Let $\bar{V}(x;\varepsilon;\lambda)$ be a nontrivial solutions of (\ref{4c2}) for some $\lambda \in \mathbb{C}$. Then, there exist constants 
$\alpha_i\ (i=1,2,3,4)$ satisfying ${\sum}_{i=1}^4|\alpha_i|\neq 0$ such that 
$\bar{V}(x;\varepsilon;\lambda)$ must be represented as 
\begin{equation}\label{eq002}
\bar{V}(x;\varepsilon;\lambda)=
\left\{
\begin{array}{ll}
\alpha_1\bar{V}_1(x;\varepsilon;\lambda)+\alpha_2\bar{V}_2(x;\varepsilon;\lambda), \ x\in [0,x^*(\ep)], \\[1ex]
\alpha_3\bar{V}_3(x;\varepsilon;\lambda)+\alpha_4\bar{V}_4(x;\varepsilon;\lambda), \ x\in [x^*(\ep),1], 
\end{array}
\right.
\end{equation}
which implies that the relation 
\begin{equation}
\begin{array}{c} \label{eq001}
\alpha_1\bar{V}_1(x^*(\varepsilon);\varepsilon;\lambda)+\alpha_2\bar{V}_2(x^*(\varepsilon);
\varepsilon;\lambda) = \ \alpha_3\bar{V}_3(x^*(\varepsilon);\varepsilon;\lambda)+\alpha_4\bar{V}_4(x^*(\varepsilon);\varepsilon;\lambda).
\end{array}
\end{equation}
holds at $x=x^*(\varepsilon)$; four vectors 
$\bar{V}_i(x^*(\varepsilon);\varepsilon;\lambda)$ $ (i=1,2,3,4)$ are linearly dependent. 
Defining 
\begin{eqnarray*}
\begin{array}{c}
g(\varepsilon;\lambda) := {\rm det}[\bar{V}_1(x^*(\varepsilon);\varepsilon;\lambda), \ \bar{V}_2(x^*(\varepsilon);\varepsilon;\lambda), 
\bar{V}_3(x^*(\varepsilon);\varepsilon;\lambda), \bar{V}_4(x^*(\varepsilon);\varepsilon;\lambda)], 
\end{array} 
\end{eqnarray*}
we find that $g(\varepsilon;\lambda)$ is an analytic function of $\lambda \in 
{\mathbb C}$ and have the next lemma: 
\begin{lemma}\label{Evans} \ 
Let $\lambda \neq 0$. 
Then,  $\lambda \in \mathbb{C}$ is an eigenvalue of 
$(\ref{4c1})$ if and only if $g(\varepsilon;\lambda)=0$. 
\end{lemma} 

We call $g(\varepsilon;\lambda)$ the {\it Evans function} of the 
transition layer solution, which enables us to investigate the distribution of eigenvalues of $(\ref{4c1})$ in $\mathbb{C}$. 
We note that 
$\lambda = 0$ is an eigenvalue of $(\ref{4c1})$ if $g(\varepsilon; \lambda)\big{|}_{\lambda = 0} =0$ and $(0, 0) \neq (p, q ) \in \pot{C}^2_1[0,1] \times \pot{C}^2_1[0,1]$ satisfying \eqref{4c1} for $\lambda = 0$ satisfies $\int_0^1 ( p + q ) dx = 0$.
To calculate the Evans function, we have to construct functions $\bar{V}_i(x;\varepsilon;\lambda) \ 
(i=1,2,3,4)$ as we constructed a transition layer solution in the previous section. According 
to the dependency of $\lambda \in \mathbb{C}$ on $\varepsilon$, 
we must divide our argument into the following three cases: \\
$(\mathrm{I}) \quad \lambda = \lambda(\varepsilon) = O(\ep) \ in \ \mathbb{C}$
as $\varepsilon \to 0$. 

For the other two cases, we have $\lambda(\ep)/\ep \to \infty$ as $\ep \to 0$. 
From \cite[Lemma 1.1.1]{Ec}, we find that there exists a positive  and continuous real function $\omega(\ep) \to \infty$ as $\ep \to 0$ such that $\lambda(\ep)$ is represented as 
$$ \lambda(\ep) = \ep \omega(\ep) \hat{\lambda}(\ep), $$
where $\hat{\lambda}(\ep)$ satisfies $\hat{\lambda}(0) \ne 0$. Then, we 
consider two cases according to the magnitude of $\ep \omega(\ep)$ as follows:\\
$(\mathrm{II}) \quad \ep \omega(\ep) \to 0$ and $\omega(\ep) \to \infty$ as $\ep \to 0$; \\
$(\mathrm{III}) \quad \ep \omega(\ep) \to \omega_0$ as $\ep \to 0$ for some positive constant $\omega_0$. 

\vspace{3mm}
In what follows, we construct $\bar{V}_i(x; \ep ; \lambda) \, (i=1,2,3,4)$
by using the same lines of argument in Section~\ref{S2.1}.

%%%%%%%%%%%%%%%%%%%%%%%%%%%%%%%%%%%%%%%%%%%%%%%%%%%%%%%%%%%%%%%%%%%%%%%%%%% 3-1
\subsection{Case $(\mathrm{I}) \ \lambda = \lambda(\varepsilon) = O(\ep)$ as $\ep \to 0$} \label{S3.1} %Section 3.1 

\subsubsection{Construction of $\bar{V}_1$ and $\bar{V}_2$} \label{S3.1.1}
Putting $\lambda(\ep) = \ep \kappa$ for $\kappa \in {\mathbb C}$,
we consider the problem (\ref{4c3}). Using the transformation $x = x^*(\ep)y$ in (\ref{4c3}), we have 
\begin{equation} \label{4c5}
\left\{
\begin{array}{l}
\begin{array}{l}
\varepsilon^2p_{yy}+x^*(\varepsilon)^2(f_u^\ep - \varepsilon\kappa)p 
+ x^*(\varepsilon)^2 f_v^\ep q = 0, \\[1ex]
D q_{yy}- x^*(\ep)^2f_u^\ep p - x^*(\ep)^2(f_v^\ep + \varepsilon\kappa)q = 0,
\end{array} y \in (0,1)\\[0.4cm]
(p_y,q_y)(0)=(0,0), \ (p,q)(1)= (a, b). 
\end{array}
\right.
\end{equation}
We first consider the outer approximation of the form 
\begin{eqnarray*}
\left\{
\begin{array}{ll}
p^-(y) = P_0^-(y) + {\varepsilon}P_1^-(y)  + O(\ep^2),\\ 
q^-(y) = Q_0^-(y) + {\varepsilon}Q_1^-(y) + O(\ep^2).
\end{array}
\right.
\end{eqnarray*}
Substituting this into (\ref{4c5}),
we equate the coefficients of the same powers of 
$\varepsilon$. \\[0.2cm]
$O({\varepsilon}^0)$:
\begin{eqnarray*}
\left\{
\begin{array}{l}
\begin{array}{l}
f_u^- P_0^- + f_v^- Q_0^- = 0, \\[1ex]
D Q_{0,yy}^- = 0,  
\end{array} y \in (0,1) \\[0.4cm]
Q_{0,y}^-(0)=0,\  Q_0^-(1)=b.
\end{array}
\right.
\end{eqnarray*}
$Q_0^-(y;b) = b$ and then $P_0^-(y;b) = -b f_v^- /f_u^-$, where $f_u^-= f_u(U^-_0,V^-_0)$
and $f_v^-= f_v(U^-_0,V^-_0)$. \\[0.2cm]
$O({\varepsilon}^1)$:
\begin{eqnarray*}
\left\{
\begin{array}{l}
\begin{array}{l}
f_u^-P_1^- + f_v^- Q_1^- - \kappa P_0^-(y;b) = 0, \\[1ex]
D Q_{1,yy}^- + x_0^2 \kappa \{ f_v^-/f_u^- - 1\}Q_0^-(y;b) = 0, \\[1ex]
\end{array} y \in (0,1)\\[0.4cm]
Q_{1,y}^-(0)=0,\  Q_1^-(1)=0. 
\end{array}
\right.
\end{eqnarray*}
We have 
\begin{eqnarray*}
\left\{
\begin{array}{l}
Q_1^-(y;\kappa;b)= b x_0^2 \kappa \{ -f_v^-/f_u^- + 1\}(y^2-1)/(2D), \\[1ex]
P_1^-(y;\kappa;b) = - \ f_v^- Q_1^-(y;\kappa;b)/f_u^- + \kappa P_0^-(y;b)/f_u^-.  \\[1ex]
\end{array}
\right.
\end{eqnarray*}
Since the $p$ component does not satisfy the boundary condition at $y=1$, 
we have to modify the above approximation by adding correction terms $\rho_i^-, \pi_i^- \ (i=0,1)$ of the 
form 
\begin{eqnarray*}
\left\{
\begin{array}{lcl}
p^-(y) &=& P_0^-(y;b) + {\varepsilon}P_1^-(y;\kappa;b) \\[1ex]
 & & + \ \rho_{0}^-(\frac{y-1}{\varepsilon})+ \ep \rho_1^-(\frac{y-1}{\varepsilon}) + O(\ep^2),\\[1ex] 
q^-(y) &=& Q_0^-(y;b) + \varepsilon Q_1^-(y;\kappa;b) \\[1ex]
 & & + \ {\varepsilon}^2 \pi_0^-(\frac{y-1}{\varepsilon}) + {\varepsilon}^3 \pi_1^-(\frac{y-1}{\varepsilon})
 + O(\ep^4).
\end{array}
\right.
\end{eqnarray*}
Substituting this into (\ref{4c5}) and using $z=(y-1)/\varepsilon$, we equate the coefficients of 
the same powers of $\varepsilon $. \\[0.2cm]
$O({\varepsilon}^{0}):$
\begin{eqnarray*}
\left\{
\begin{array}{l}
\begin{array}{l}
\ddot{\rho}_0^- + x_0^2\tilde{f}_u^-\rho_0^- = b x_0^2 \{\tilde{f}_u^- f_v^-/f_u^- - \tilde{f}_v^-\}, \\[1ex]
D \ddot{\pi}_0^- + \ddot{\rho}_0^-  =0, \\[1ex]
\end{array} \qquad z \in (-\infty ,0) \\[0.4cm]
\rho_0^-(-\infty)=0,\ \rho_0^-(0) = a+ b f_v^-/f_u^- , \\[1ex]
\pi_0^-(-\infty)=0,\ \dot{\pi}_0^-(-\infty)=0, 
\end{array}
\right.
\end{eqnarray*}
where $\tilde{f}_u^- = f_u(h^-(v^*)+\phi^-_0,v^*), \tilde{f}_v^- = f_v(h^-(v^*)+\phi^-_0,v^*)$,
and the dot notation denotes $d/dz$.
In the same way as the solution of \eqref{b19} is represented by \eqref{b21}, we have 
\begin{eqnarray}\label{com1}
\begin{array}{l}
\rho_{0}^-(z;a,b) = (a+ b f_v^- /f_u^- ) \dot{\phi}_0^-(z)/\dot{\phi}_0^-(0)\\[0.2cm]
\hspace{1cm}- \ b x_0^2  \dot{\phi}_0^-(z)
\dis\int_z^0 (\dot{\phi}_0^-(\eta))^{-2} \dis\int_{-\infty}^\eta \dot{\phi}_0^-(\zeta)\{\tilde{f}_u^- f_v^-/f_u^- - \tilde{f}_v^-\}d\zeta d\eta
\end{array}
\end{eqnarray}
and $\pi_0^-(z;a,b)=- \rho_0^-(z;a,b)/D$.  \\[0.2cm]
$O({\varepsilon}^1)$: \\
Noting that $U^-_1= V^-_1= 0$ (see Remark \ref{rem2.2}) and $Q_1^-(1; \kappa; b)=0$, we have 
\begin{equation*}
\left\{
\begin{array}{l}
\begin{array}{l}
\ddot{\rho}_1^- + x_0^2\tilde{f}_u^-\rho_1^- = -\{ 2x_0x_1\tilde{f}_u^- +  
x_0^2 \tilde{f}_{uu}^-\phi_1^-  - x_0^2 \kappa \}\rho_0^-(z;a,b) \\[1ex]
\hspace*{3.5cm} - \ x_0^2 \tilde{f}_u^-P_1^-(1;\kappa;b)
+ b R_1^- , \\[1ex]
D \ddot{\pi}_1^- + \ddot{\rho}_1^- = x_0^2 \kappa \rho_0^-(z;a,b), \\[1ex]
\end{array} z \in (-\infty ,0)\\[0.4cm]
\rho_1^-(-\infty)=0,\ \rho_1^-(0)= -P_1^-(1;\kappa;b), \\[1ex]
\pi_1^-(-\infty)=0,\ \dot{\pi}_1^-(-\infty)=0,
\end{array}
\right.
%\label{eqn:B4}
\end{equation*}
where 
$R_1^- :=  \{x_0^2(\tilde{f}_{uu}^-\phi_1^- -\kappa)+2x_0x_1\tilde{f}_u^-\}f_v^-/f_u^-
 - (x_0^2\tilde{f}_{vu}^-\phi_1^- +2x_0x_1\tilde{f}_v^-)$,
$\tilde{f}_{uu}^- := f_{uu}(h^-(v^*)+\phi^-_0,v^*)$, 
and 
$\tilde{f}_{vu}^- := f_{vu}(h^-(v^*)+\phi^-_0,v^*)$. 
We obtain 
\begin{equation} 
\left\{
\begin{array}{l}
\rho_1^-(z;\kappa;a,b)= -P_1^-(1;\kappa;b)\dot{\phi}_0^-(z)/\dot{\phi}_0^-(0) 
\\[1ex]
\hspace*{0.5cm}+ \ \dot{\phi}_0^-(z)\int_{z}^{0}(\dot{\phi}_0^-(\eta))^{-2}\int_{-\infty}^{\eta}
\dot{\phi}_0^-(\zeta) 
\\[1ex]
\hspace*{1.5cm} \times \{ (2x_0x_1\tilde{f}_u^- + x_0^2 \tilde{f}_{uu}^-\phi_1^-  - x_0^2 \kappa )
\rho_0^-(\zeta;a,b) 
\\[1ex]
\hspace*{1.5cm}+ \ x_0^2 \tilde{f}_u^-P_1^-(1;\kappa;b)  - b R_1^-  \} d{\zeta}d{\eta}, 
\\[1ex]
\pi_1^-(z;\kappa;a,b)= \{ -\rho_1^-(z;\kappa;a,b) 
\\[1ex]
\hspace*{2.5cm}+ \ x_0^2 \kappa\int_{-\infty}^{z}\int_{-\infty}^{\eta} \rho_0^-(\zeta;a,b) d{\zeta}d{\eta}\}/D.
\end{array} \right.
\label{eqn:C4}
\end{equation}

%
%%%%%%%%%%%%%%%%%%%%%%%%%%%%%%%%%%%%%%%%%%%%%%%%%%%%%%%%%%%%%%%%% 
\indent For any fixed $\kappa^* \in {\mathbb C}$, let us define 
$\Delta_{\nu} := \{ \kappa \in {\mathbb C} \ | 
\ |\kappa - \kappa^*| \leq \nu \}$ for some $\nu > 0$.
Using the above approximate solutions, we can construct uniform approximations up to 
$O({\varepsilon})$ of (\ref{4c5}) for any $\kappa \in \Delta_{\nu}$, which take the form
\begin{eqnarray*}
\left\{
\begin{array}{lcl}
P^-(y;{\varepsilon};\kappa;a,b) &=& P_0^-(y;b) + {\varepsilon}P_1^-(y;\kappa;b) \\[1ex]
   &+&{\theta}(y)\{\rho_{0}^-(\frac{y-1}{\varepsilon};a,b) + \ep 
           \rho_1^-(\frac{y-1}{\varepsilon};\kappa;a,b)\},\\[1ex] 
Q^-(y;{\varepsilon};\kappa;a,b) &=& Q_0^-(y;b) + {\varepsilon}Q_1^-(y;\kappa;b) \\[1ex]
 &+&{\theta}(y)\{{\varepsilon}^2\pi_0^-(\frac{y-1}{\varepsilon};a,b) + 
           {\varepsilon}^3\pi_1^-(\frac{y-1}{\varepsilon};\kappa;a,b) \\[1ex]
 & - & \varepsilon^2\pi_0^-(0;a,b) - \ep^3 \pi_1^-(0;\kappa;a,b)\},
\end{array}
\right.
\end{eqnarray*}
where ${\theta}(y)$ is the same function as is defined in Section~2.1. 
It is clear that $(P^-, Q^-)(y;{\varepsilon};\kappa;a,b)$ 
satisfies the both boundary conditions at $y=0$ and $y=1$. 
Hence, we look for exact solutions of (\ref{4c5}) of the form
\begin{eqnarray*}
\left\{
\begin{array}{lcl}
p^-(y;{\varepsilon};\kappa;a,b) &=& P^-(y;{\varepsilon};\kappa;a,b)
+{\varepsilon}\tilde{p}^-(y;{\varepsilon};\kappa;a,b) \\[1ex]
 & & + \ \ep h_v^-(V_0^-(y))\tilde{q}^-(y;{\varepsilon};\kappa;a,b),\\[1ex]
q^-(y;{\varepsilon};\kappa;a,b) &=& 
Q^-(y;{\varepsilon};\kappa;a,b)
    +{\varepsilon}\tilde{q}^-(y;{\varepsilon};\kappa;a,b). 
\end{array}
\right.
\end{eqnarray*}
Substituting this into (\ref{4c5}), we have 
\begin{equation}\label{b222_1}
\left \{
\begin{array}{l}
\begin{array}{l}
\ep^2 \tilde{p}^-_{yy} + \ep^2 (h^-_v(V^-_0)\tilde{q}^-)_{yy} \\[1ex]
\hspace*{1cm} + \ (x^*(\ep))^2 \{(f_u^\ep - \ep \kappa)(\tilde{p}^- + h^-_v(V^-_0)\tilde{q}^-) + f_v^\ep \tilde{q}^-\} \\[1ex]
\hspace*{1cm} + \ \ep P^-_{yy} + \ (x^*(\ep))^2 \{(f_u^\ep - \ep \kappa)P^- + f_v^\ep Q^-\}/\ep \ = \ 0,   \\[1ex]
D \tilde{q}^-_{yy} - \ (x^*(\ep))^2\{ f_u^\ep (\tilde{p}^- + h^-_v(V^-_0)\tilde{q}^-) + (f_v^\ep + \ep \kappa) \tilde{q}^-\} \\[1ex]
\hspace*{1cm}+ \ D Q^-_{yy}/\ep  - \ (x^*(\ep))^2 \{f_u^\ep P^- + (f_v^\ep + \ep \kappa)Q^-\}/\ep \ = \ 0,   \\[1ex]
\end{array}  \ y \in (0, 1) \\[1ex]
(\tilde{p}^-,\tilde{q}^-)_y(0) = (0,0), \ (\tilde{p}^-,\tilde{q}^-)(1) = (0,0).
\end{array} 
\right.
\end{equation}
To solve (\ref{b222_1}), for $w = (\tilde{p}^-,\tilde{q}^-)$, we define the following operator 
$T(w; {\varepsilon};\kappa;a,b) = (T_1,T_2)$ $(w; {\varepsilon};\kappa;a,b)$:
\begin{eqnarray*}
\left\{
\begin{array}{ll}
T_1(w; {\varepsilon};\hat{\lambda};a,b,c) := 
\ep^2 \tilde{p}^-_{yy} + \ep^2 (h^-_v(V^-_0)\tilde{q}^-)_{yy} \\[1ex]
\hspace*{0.5cm} + \ (x^*(\ep))^2 \{(f_u^\ep - \ep \kappa)(\tilde{p}^- + h^-_v(V^-_0)\tilde{q}^-) + f_v^\ep \tilde{q}^-\} \\[1ex]
\hspace*{0.5cm} + \ \ep P^-_{yy} + \ (x^*(\ep))^2 \{(f_u^\ep - \ep \kappa)P^- + f_v^\ep Q^-\}/\ep, \\[1ex]
T_2(w; {\varepsilon};\hat{\lambda};a,b,c) := 
D \tilde{q}^-_{yy} - \ (x^*(\ep))^2\{ f_u^\ep (\tilde{p}^- + h^-_v(V^-_0)\tilde{q}^-) \\[1ex]
\hspace*{0.5cm}+ \  (f_v^\ep + \ep \kappa) \tilde{q}^-\} + \ D Q^-_{yy}/\ep  -  
(x^*(\ep))^2\{f_u^\ep P^- + (f_v^\ep + \ep \kappa)Q^-\}/\ep.
\end{array}
\right.
\end{eqnarray*}
We find that (\ref{b222_1}) is rewritten as $T(w; {\varepsilon};
\kappa;a,b)=0$, and that
$T(w; \varepsilon;\kappa;a,b)$ is a 
continuously differentiable mapping from $X_\ep^- \times (0,\varepsilon_0) \times \Delta_\nu$ to $Y$, where 
$X_\ep^-$ and $Y$ are defined in Section~2.1. 
\begin{lemma} For any given $\kappa^* \in {\mathbb C}$, there exist positive constants 
$\varepsilon_0, \nu_0$ and $K$ such that for any $\varepsilon \in (0,\varepsilon_0)$ and $\kappa \in \Delta_{\nu_0}$, 
%\vspace*{-0.4cm}
\begin{eqnarray*}
&\mathrm{(i)}&\|T(0;\varepsilon;\kappa;a,b)\|_Y = o(1) \ \mbox{uniformly in } \kappa \in  \Delta_{\nu_0} \ \mbox{as} \ \varepsilon \to 0;\\
&\mathrm{(ii)}&\|T_w(w_1; \varepsilon;\kappa;a,b)-T_w(w_2; \varepsilon;\kappa;a,b)\|_
{X_\ep^- \to Y}
\leqq K\|w_1-w_2\|_{X_\ep^-} \hspace*{2cm} \\
& & \hspace*{6.5cm} \mbox{for any} \ w_1, w_2 \in X_\ep^-; \\
&\mathrm{(iii)}&\|T_w^{-1}(0;\varepsilon;\kappa;a,b)\|_{Y\to {X_\ep^-}}\leq K.
\end{eqnarray*}
Moreover, the results (i) - (iii) hold for $\partial T/\partial \kappa$ in place of $T$.
\label{L3.5}
\end{lemma}

By this lemma, we can apply the Implicit Function Theorem to $T(w; \varepsilon;
\kappa;a,b)=0$, and hence there exists 
$w(\varepsilon;\kappa;a,b) \in X_\ep^-$ satisfying $T(w; \varepsilon;\kappa;a,b)=0$
under the same assumption of Lemma \ref{L3.5}. Moreover, 
$w(\varepsilon;\kappa;a,b)$ and $\partial w/\partial \kappa(\varepsilon;\kappa;a,b)$ 
are uniformly continuous with respect to $(\varepsilon,\kappa) \in (0,\varepsilon_0)\times \Delta_{\nu_0}$ in the $X_\ep^-$-topology, and satisfy 
$$ 
\begin{array}{l}
\|w(\varepsilon;\kappa;a,b) \|_{X_\ep^-}, \ \|\partial w/\partial 
\kappa(\varepsilon;\kappa;a,b) \|_{X_\ep^-} =  o(1) \hspace{2cm} \\ 
\hspace*{4cm} \mbox{as} \ \varepsilon \to 0 
\ \mbox{uniformly in} \ \kappa \in \Delta_{\nu_0}.
\end{array}
$$
Thus, we have exact solutions of (\ref{4c5}) on $[0,1]$ of the form 
\begin{eqnarray*}
\left\{
\begin{array}{lcl}
p^-(y;{\varepsilon};\kappa;a,b) &=& P^-(y;{\varepsilon};\kappa;a,b)
+{\varepsilon}\tilde{p}^-(y;{\varepsilon};\kappa;a,b) \\[1ex]
 & & + \ \ep h_v^-(V_0^-(y))\tilde{q}^-(y;{\varepsilon};\kappa;a,b),\\[1ex]
q^-(y;{\varepsilon};\kappa;a,b) &=& 
Q^-(y;{\varepsilon};\kappa;a,b)
    +{\varepsilon}\tilde{q}^-(y;{\varepsilon};\kappa;a,b). 
\end{array}
\right.
\end{eqnarray*} 
and then we obtain solutions of \eqref{4c3} on $[0,x^*({\varepsilon})]$ 
\begin{equation}\label{b312xx}
\left\{
\begin{array}{lcl}
p^-(x;{\varepsilon};\kappa;a,b) &=& P^-(\frac{x}{x^*(\varepsilon)};{\varepsilon};\kappa;a,b)
+{\varepsilon}\tilde{p}^-(\frac{x}{x^*(\varepsilon)};{\varepsilon};\kappa;a,b) \\[1ex]
 & & + \ \ep h_v^-(V_0^-(\frac{x}{x^*(\varepsilon)}))\tilde{q}^-(\frac{x}{x^*(\varepsilon)};{\varepsilon};\kappa;a,b),\\[1ex]
q^-(x;{\varepsilon};\kappa;a,b) &=& 
Q^-(\frac{x}{x^*(\varepsilon)};{\varepsilon};\kappa;a,b)
    +{\varepsilon}\tilde{q}^-(\frac x{x^*(\varepsilon)};{\varepsilon};\kappa;a,b). 
\end{array}
\right.
%\label{abc1}
\end{equation}

%
%%%%%%%%%%%%%%%%%%%%%%%%%%%%%%%%%%%%%%%%%%%%%%%%%%%%%%%%%%%%%%%%%%%3-1-2
\subsubsection{Construction of $\bar{V}_3$ and $\bar{V}_4$} \label{S3.1.2}
Next, we consider the problem (\ref{4c4}). By using the transformation 
$x = x^*({\varepsilon}) +(1-x^*({\varepsilon}))y$, 
we have 
\begin{equation} \label{eqn:n2-1}
\left\{
\begin{array}{l}
\begin{array}{l}
{\varepsilon}^2p_{yy}+(1-x^*({\varepsilon}))^2(f_u^\ep -{\varepsilon}\kappa)p + 
(1-x^*({\varepsilon}))^2f_v^\ep q = 0, \\[1ex]
Dq_{yy}-(1-x^*({\varepsilon}))^2f_u^\ep p 
-(1-x^*({\varepsilon}))^2(f_v^\ep+{\varepsilon}\kappa)q = 0, 
\end{array}
 y \in (0,1)\\[3ex]
(p,q)(0) = (a,b), \ (p_y,q_y)(1)=(0,0).
\end{array}
\right.
\end{equation}
First, we consider outer approximations of the form 
\begin{eqnarray*}
\left\{
\begin{array}{l}
p^+(y) = P_0^+(y) + \varepsilon P_1^+(y) + O({\varepsilon}^2),\\[1ex]
q^+(y) = Q^+_0(y) + \varepsilon Q_1^+(y) + O({\varepsilon}^2).
\end{array}
\right.
\end{eqnarray*}
Substituting this into (\ref{eqn:n2-1}), we equate the coefficients of the same powers of 
$\varepsilon$ \\[0.2cm]
$O(\varepsilon^0)$:
\begin{eqnarray*}
\left\{
\begin{array}{l}
\begin{array}{l}
f_u^+ P_0^+ + f_v^+ Q_0^+ = 0, \\[0.1cm]
D Q_{0,yy}^+ = 0,
\end{array} y \in (0,1 ) \\[0.3cm]
Q_0^+(0)=b, \ Q_{0,y}^+(1)=0.
\end{array} \right.
\end{eqnarray*}
We easily find that $Q_0^+(y;b) = b $ and $P_0^+(y;b) = - b f_v^+ /f_u^+$, where $f_u^+ = f_u(U_0^+, V_0^+)$ and $f_v^+ = f_v(U_0^+, V_0^+)$.\\[0.2cm]
$O(\varepsilon^1)$:
\begin{eqnarray*}
\left\{
\begin{array}{l}
\begin{array}{l}
f_u^+ P_1^+ + f_v^+ Q_1^+ - \kappa P_0^+(y;b) = 0, \\[0.1cm]
D Q_{1,yy}^+ + (1-x_0)^2 \kappa \{ f_v^+/f_u^+ - 1\}Q_0^+(y;b)=0,
\end{array} y \in (0,1) \\[0.3cm]
Q_1^+(0)=0,\  Q_{1,y}^+(1)=0.
\end{array}
\right.
\end{eqnarray*}
We have 
\begin{eqnarray*}
\left\{
\begin{array}{l}
Q_1^+(y;\kappa;b)= b (1-x_0)^2 \kappa \{ -f_v^+/f_u^+ + 1\}(y^2-2y)/(2D), \\[1ex]
P_1^+(y;\kappa;b) = - \ f_v^+ Q_1^+(y;\kappa;b)/f_u^+ + \kappa P_0^+(y;b)/f_u^+. \\[1ex]
\end{array}
\right.
\end{eqnarray*}
%
%
%%%%%%%%%%%%%%%%%%%%%%%%%%%%%%%%%%%%%%%%%%%%%%%%%%%%%%%%%%%%%%%%%%%%%%%%%%%%%%%%%%%%%%%%%%%%%%%  3-1-2
\indent Since the $p$ component does not satisfy the boundary condition at $y=0$, we have to modify the above approximation by 
adding correction terms $\rho_i^+, \pi_i^+ \ (i=0,1)$ of the form 
\begin{eqnarray*}
\left\{
\begin{array}{lcl}
p^+(y) &=& P_0^+(y;b) + \varepsilon P_1^+(y;\kappa;b)
 +  \rho_{0}^+(\frac{y}{\varepsilon}) + \ep \rho_1^+(\frac{y}{\varepsilon})+O(\ep^2),\\[0.1cm] 
q^+(y) &=& Q_0^+(y;b) + \varepsilon Q_1^+(y;\kappa;b) 
 +  {\varepsilon}^2\pi_0^+(\frac{y}{\varepsilon}) + {\varepsilon}^3\pi_1^+(\frac{y}{\varepsilon})+ O(\ep^4).
\end{array}
\right.
\end{eqnarray*}
Substituting this into (\ref{eqn:n2-1}) and using $z=y/{\varepsilon}$, we equate the coefficient 
of the same powers of $\varepsilon$. \\[0.2cm]
$O({\varepsilon}^0)$:
\begin{eqnarray*}
\left\{
\begin{array}{l}
\begin{array}{l}
\ddot{\rho}_0^+ + (1-x_0)^2 \tilde{f}_u^+\rho_0^+= b(1-x_0)^2\{
\tilde{f}_u^+ f_v^+/f_u^+ - \tilde{f}_v^+\}, \\[0.1cm]
D \ddot{\pi}_0^+ + \ddot{\rho}_{0}^+ = 0,
\end{array} z \in (0 ,\infty)\\[0.3cm]
\rho_0^+(0)=a + b f_v^+/f_u^+, \ \rho_0^+(\infty)=0, \\[0.1cm]
\pi_0^+(\infty)=0, \ \dot{\pi}_0^+(\infty)=0, 
\end{array}
\right.
\end{eqnarray*}
where $\tilde{f}_u^+ = f_u(h^+(v^*)+\phi^+_0,v^*), \tilde{f}_v^+ = f_v(h^+(v^*)+\phi^+_0,v^*)$
and the dot notation denotes $d/dz$. 
In the same way as the solution of \eqref{b30} is represented by \eqref{b32}, we have
\begin{eqnarray}\label{c22}
\begin{array}{l}
\hspace*{0.5cm}\rho_{0}^+(z;a,b)=(a + b f_v^+/f_u^+)\dot{\phi}_0^+(z)/\dot{\phi}_0^+(0) \\[0.2cm]
\hspace*{1cm}- \ b (1-x_0)^2  \dot{\phi}_0^+(z) \dis\int_0^z (\dot{\phi}_0^+(\eta))^{-2} \int_\eta^\infty \dot{\phi}_0^+(\zeta)\{\tilde{f}_u^+ f_v^+/f_u^+ - \tilde{f}_v^+\}d\zeta d\eta,
\end{array}
\end{eqnarray}
and  
$\pi_0^+(z;a,b)=- \rho_0^+(z;a,b)/D$. \\[0.2cm]
$O({\varepsilon}^1)$: \\
Noting $U^+_1 = V^+_1 = 0$ and $Q_1^+(0; \kappa, b) =0$, we have
\begin{eqnarray*}
\left\{
\begin{array}{l}
\begin{array}{l}
 \ddot{\rho}_1^+ + (1-x_0)^2 \tilde{f}_u^+ \rho_1^+ = \{ 2(1-x_0)x_1\tilde{f}_u^+ -(1-x_0)^2(\tilde{f}_{uu}^+\phi_1^+   \\[1ex]
\hspace*{0.5cm}  - \ \kappa) \}\rho_0^+(z;a,b) - \ (1-x_0)^2 \tilde{f}_u^+P_1^+(0;\kappa;b) 
+ \ b R_1^+,  \\[1ex]
D \ddot{\pi}_1^+ + \ddot{\rho}_1^+ = (1-x_0)^2 \kappa \rho_0^+(z;a,b), 
\end{array}  z \in (0 ,\infty) \\[0.7cm]
 \rho_1^+(0) = -P_1^+(0;\kappa;b), \ \rho_1^+(\infty)=0, \\[0.1cm]
 \pi_1^+(\infty)=0,\ \dot{\pi}_1^+(\infty)=0,
\end{array}
\right.
\end{eqnarray*}
where 
$R_1^+ :=  \{(1-x_0)^2(\tilde{f}_{uu}^+\phi_1^+ -\kappa)-2(1-x_0)x_1\tilde{f}_u^+\}f_v^+ /f_u^+ 
- \{(1-x_0)^2\tilde{f}_{vu}^+\phi_1^+ -2(1-x_0)x_1\tilde{f}_v^+\}$,
$\tilde{f}_{uu}^+ := f_{uu}(h^+(v^*)+\phi^+_0,v^*)$, 
and 
$\tilde{f}_{vu}^+ := f_{vu}(h^+(v^*)+\phi^+_0,v^*)$. 
We obtain 
\begin{equation} 
\left\{ 
\begin{array}{l}
\rho_1^+(z;\kappa;a,b) = -P_1^+(0;\kappa;b)\dot{\phi}_0^+(z)/\dot{\phi}_0^+(0) + \dot{\phi}_0^+(z)\int_{0}^{z}(\dot{\phi}_0^+(\eta))^{-2}  \\[1ex]
\hspace*{1cm} \times \ \int^{\infty}_{\eta}\dot{\phi}_0^+(\zeta) \{ ( -2(1-x_0)x_1\tilde{f}_u^+ + (1-x_0)^2 \tilde{f}_{uu}^+\phi_1^+ \\[1ex]
\hspace*{1cm} - \ (1- x_0)^2 \kappa ) \rho_0^+(\zeta;a,b) + \ (1-x_0)^2 \tilde{f}_u^-P_1^-(0;\kappa;b)   - b R_1^+ \} d{\zeta}d{\eta}, \\[0.2cm]
\pi_1^+(z;\kappa;a,b) = \{-\rho_1^+(z;\kappa;a,b) \\[1ex]
\hspace*{1cm} + \ (1-x_0)^2 \kappa 
\int_{z}^{\infty}\int_{\eta}^{\infty} \rho_0^+(\zeta;a,b)d{\zeta}d{\eta}\}/D. 
\end{array}
\right.
\label{C6}
\end{equation}

%
%%%%%%%%%%%%%%%%%%%%%%%%%%%%%%%%%%%%%%%%%%%%%%%%%%%%%%%%%%%%%%%%%%%%%%%%%%%%%  3-1-2
\indent Using the above approximate solutions, we can construct uniform approximations 
up to $O({\varepsilon})$ of (\ref{eqn:n2-1}) for any fixed $\kappa \in \Delta_\nu$, 
which takes the form 
\begin{eqnarray*}
\left\{
\begin{array}{lcl}
P^+(y;{\varepsilon};\kappa;a,b) &=& P_0^+(y;b) + {\varepsilon}P_1^+(y;\kappa;b) \\[1ex]
           &+& \theta(1-y)\{\rho_0^+(\frac{y}{\varepsilon};a,b) + \ep \rho_1^+(\frac{y}{\varepsilon};\kappa;a,b)\},\\[1ex] 
Q^+(y;{\varepsilon};\kappa;a,b) &=& Q_0^+(y;b) + {\varepsilon}Q_1^+(y;\kappa;b)  \\[1ex]
	   &+&\theta(1-y)\{{\varepsilon}^2\pi_0^+(\frac{y}{\varepsilon};a,b) + {\varepsilon}^3\pi_1^+(\frac{y}{\varepsilon};\kappa;a,b) \\[1ex]
 &	- & {\varepsilon}^2\pi_0^+(0;a,b) - \ep^3\pi_1^+(0;\kappa;a,b) \}. 
\end{array}
\right.
\end{eqnarray*}
It is clear that $(P^+,Q^+)(y;{\varepsilon};\kappa;a,b)$ satisfies the both 
boundary conditions at $y=0$ and $y=1$.
Hence, we look for 
exact solutions of (\ref{eqn:n2-1}), which take the form 
\begin{eqnarray*}
\left\{
\begin{array}{lcl}
p^+(y;{\varepsilon};\kappa;a,b) &=& P^+(y;{\varepsilon};\kappa;a,b)
+{\varepsilon}\tilde{p}^+(y;{\varepsilon};\kappa;a,b) \\[1ex]
 & & + \ \ep h_v^+(V_0^+(y))\tilde{q}^+(y;{\varepsilon};\kappa;a,b),\\[1ex]
q^+(y;{\varepsilon};\kappa;a,b) &=& 
Q^+(y;{\varepsilon};\kappa;a,b)
    +{\varepsilon}\tilde{q}^+(y;{\varepsilon};\kappa;a,b). 
\end{array}
\right.
\end{eqnarray*}
Substituting this into (\ref{eqn:n2-1}), we have 
\begin{equation}\label{eqn:d4}
\left \{
\begin{array}{l}
\begin{array}{l}
\ep^2 \tilde{p}^+_{yy} + \ep^2 (h^+_v(V^+_0)\tilde{q}^+)_{yy} \\[1ex]
\hspace*{1cm} + \ (1-x^*(\ep))^2 \{(f_u^\ep - \ep \kappa)(\tilde{p}^+ + h^+_v(V^+_0)\tilde{q}^+) + f_v^\ep \tilde{q}^+ \} \\[1ex]
\hspace*{1cm} + \ \ep P^+_{yy} + \ (1-x^*(\ep))^2 \{(f_u^\ep - \ep \kappa)P^+ + f_v^\ep Q^+\}/\ep \ = \ 0,   \\[1ex]
D \tilde{q}^+_{yy} - \ (1-x^*(\ep))^2 \{f_u^\ep (\tilde{p}^+ + h^+_v(V^+_0)\tilde{q}^+) + (f_v^\ep + \ep \kappa) \tilde{q}^+ \} \\[1ex]
\hspace*{1cm}+ \ D Q^+_{yy}/\ep  - \ (1-x^*(\ep))^2 \{f_u^\ep P^+ + (f_v^\ep + \ep \kappa)Q^+\}/\ep \ = \ 0,   \\[1ex]
\end{array}  \ y \in (0, 1) \\[1ex]
(\tilde{p}^+,\tilde{q}^+)(0) = (0,0), \ (\tilde{p}^+,\tilde{q}^+)_y(1) = (0,0).
\end{array} 
\right.
\end{equation}
Applying a similar argument as applied to \eqref{b222_1}, we obtain 
solutions of \eqref{eqn:d4},
which has the same properties as those of \eqref{b222_1}. 
Thus, we have exact solutions of (\ref{eqn:n2-1}) on $[0,1]$ of the form 
\begin{eqnarray*}
\left\{
\begin{array}{lcl}
p^+(y;{\varepsilon};\kappa;a,b) &=& P^+(y;{\varepsilon};\kappa;a,b)
+{\varepsilon}\tilde{p}^+(y;{\varepsilon};\kappa;a,b) \\[1ex]
 & & + \ \ep h_v^+(V_0^+(y))\tilde{q}^+(y;{\varepsilon};\kappa;a,b),\\[1ex]
q^+(y;{\varepsilon};\kappa;a,b) &=& 
Q^+(y;{\varepsilon};\kappa;a,b)
    +{\varepsilon}\tilde{q}^+(y;{\varepsilon};\kappa;a,b),
\end{array}
\right.
\end{eqnarray*} 
where $||(\tilde{p}^+, \tilde{q}^+)(y;{\varepsilon};\kappa;a,b)||_{X^+_\ep}$, $||{\partial}(\tilde{p}^+, \tilde{q}^+)/{\partial \kappa}(y;{\varepsilon};\kappa;a,b)||_{X^+_\ep} = o(1)$ as $\ep \to 0$ uniformly in $\kappa \in \Delta_{\nu_0}$, 
and then we obtain solutions of \eqref{4c4} on $[x^*({\varepsilon}),1]$ 
\begin{equation}\label{b317xx} 
\left\{
\begin{array}{lcl}
p^+(x;{\varepsilon};\kappa;a,b) &=& P^+(\frac{x-x^*(\varepsilon)}{1-x^*(\varepsilon)};{\varepsilon};\kappa;a,b)
+{\varepsilon}\tilde{p}^+(\frac{x-x^*(\varepsilon)}{1-x^*(\varepsilon)};{\varepsilon};\kappa;a,b) \\[1ex]
 & & + \ \ep h_v^+(V_0^+(\frac{x-x^*(\varepsilon)}{1-x^*(\varepsilon)}))\tilde{q}^+(\frac{x-x^*(\varepsilon)}{1-x^*(\varepsilon)};{\varepsilon};\kappa;a,b),\\[1ex]
q^+(x;{\varepsilon};\kappa;a,b) &=& 
Q^+(\frac{x-x^*(\varepsilon)}{1-x^*(\varepsilon)};{\varepsilon};\kappa;a,b)
    +{\varepsilon}\tilde{q}^+(\frac {x-x^*(\varepsilon)}{1-x^*(\varepsilon)};{\varepsilon};\kappa;a,b). 
\end{array}
\right.
%\label{abc2}
\end{equation}
%

%%%%%%%%%%%%%%%%%%%%%%%%%%%%%%%%%%%%%%%%%%%%%%%%%%%%%%%%%%%%%%%%%%%%%%%%%%%%%%%%%%%%%%%%%%%%%%%%%%%%%%
% Section 3.1.3
\subsubsection{Evans function $g(\varepsilon;\ep \kappa)$ corresponding to (\ref{4c1}) in  the case (I)} \label{S3.1.3}
First, we note that 
$$
%\begin{array}{l}
\dis \frac 1{x^*(\ep)} = \dis \frac 1{x_0} - \ep \frac {x_1}{x_0^2} + o(\ep), \quad 
\dis \frac 1{1 - x^*(\ep)} = \dis \frac 1{1 - x_0} + \ep \frac {x_1}{(1 - x_0)^2} + o(\ep). 
%\end{array}
$$
By using $(p^-,q^-)(x;\varepsilon;\kappa;a,b)$ and 
$(p^+,q^+)(x;\varepsilon;\kappa;a,b)$ given by \eqref{b312xx} and \eqref{b317xx}, respectively,
we can calculate 
$\bar{V}^{(i)}(\varepsilon;\kappa) := \bar{V}_i(x^*(\varepsilon);\varepsilon;\ep \kappa)\ 
(i=1,2,3,4)$ (see \eqref{v12} and \eqref{v34}) as follows: 
$$
\begin{array}{l}
\bar{V}^{(1)}(\varepsilon;\kappa)
=
\left[
\begin{array}{c}
1 \\
a_{10}+\varepsilon a_{11}+o(\varepsilon)\\
0\\
\varepsilon a_{21}+\varepsilon^2 a_{22}+o(\varepsilon^2)\\
\end{array}
\right], 
\bar{V}^{(2)}(\varepsilon;\kappa)
=
\left[
\begin{array}{c}
0\\
b_{10} + \ep b_{11} + o(\varepsilon)\\
1\\
\varepsilon b_{21} + \ep^2 b_{22} +o(\varepsilon^2)\\
\end{array}
\right],
\end{array}
$$
$$
\begin{array}{l}
\bar{V}^{(3)}(\varepsilon;\kappa)
=
\left[
\begin{array}{c}
1 \\
c_{10} + \varepsilon c_{11}+o(\varepsilon)\\
0\\
\varepsilon c_{21} + \varepsilon^2 c_{22}+o(\varepsilon^2)\\
\end{array}
\right],
\bar{V}^{(4)}(\varepsilon;\kappa)
=
\left[
\begin{array}{c}
0\\
d_{10} + \ep d_{11}+ o(\varepsilon)\\
1\\
\varepsilon d_{21} + \ep^2 d_{22}+ o(\varepsilon^2)\\
\end{array}
\right], 
\end{array} 
$$
where 
\begin{eqnarray*}
\begin{array}{l}
a_{10} = \dot{\rho}^-_0(0;1,0)/x_0, \ a_{11} = \dot{\rho}^-_1(0;\kappa;1,0)/x_0 - x_1 \dot{\rho}^-_0(0;1,0)/(x_0)^2, \\[0.1cm]
a_{21} = \dot{\pi}^-_0(0;1,0)/x_0, \ a_{22} = \dot{\pi}^-_1(0;\kappa;1,0)/x_0 -  x_1 \dot{\pi}^-_0(0;1,0)/(x_0)^2, \\[0.1cm]
b_{10} = \dot{\rho}^-_0(0;0,1)/x_0, \ b_{21} = (Q^-_{1,y}(1;\kappa;1) + \dot{\pi}^-_0(0;0,1))/x_0, \\[0.1cm]
c_{10} = \dot{\rho}^+_0(0;1,0)/(1-x_0), \ c_{11} = \dot{\rho}^+_1(0;\kappa;1,0)/(1-x_0) + x_1 \dot{\rho}^+_0(0;1,0)/(1-x_0)^2, \\[0.1cm]
c_{21} = \dot{\pi}^+_0(0;1,0)/(1-x_0), \ c_{22} = \dot{\pi}^+_1(0;\kappa;1,0)/(1-x_0) +  x_1 \dot{\pi}^+_0(0;1,0)/(1-x_0)^2, \\[0.1cm]
d_{10} = \dot{\rho}^+_0(0;0,1)/(1-x_0), \ d_{21} = (Q^+_{1,y}(0;\kappa;1) + \dot{\pi}^+_0(0;0,1))/(1-x_0), 
\end{array}
\end{eqnarray*}
and $b_{11}, b_{22}, d_{11}, d_{22}$ can be similarly calculated, but 
they are not necessary in our calculation.

By the relations $\rho_0^{\pm}(z;1,0) = \dot{\phi}_0^{\pm}(z)/\dot{\phi}_0^{\pm}(0)$, $\pi_0^{\pm}(z;1,0) = - \rho_0^{\pm}(z;1,0)/D$ and \eqref{b18}, \eqref{b29} and \eqref{b00009}, we have 
\begin{eqnarray*}
\hspace{0.5cm} a_{10} \ = \ c_{10} \ \ \mbox{and} \ \ a_{21} \ = \ c_{21}, \hspace*{2cm}
\end{eqnarray*}
and hence we find that $g(\varepsilon;\varepsilon \kappa)$ is represented as 
follows:
\begin{eqnarray*}
\begin{array}{lcl}
g(\varepsilon;\varepsilon \kappa) & = &
\det [\bar{V}^{(1)}(\varepsilon;\kappa), \bar{V}^{(2)}(\varepsilon;\kappa), \bar{V}^{(3)}(\varepsilon;\kappa), \bar{V}^{(4)}(\varepsilon;\kappa) ] \\[1ex]
 & = & \left\{(a_{22}-c_{22})(b_{10}-d_{10}) 
 -(a_{11}-c_{11})(b_{21}-d_{21}) \right\}\varepsilon^2 + o(\varepsilon^2). 
\end{array}
\end{eqnarray*}
To determine $\kappa = \kappa(\varepsilon)$ satisfying $g(\varepsilon;\varepsilon \kappa) = 0$, we define $\tilde{g}(\varepsilon;\kappa) := g(\varepsilon;\varepsilon \kappa)/\varepsilon^2$. That is,
$$
\tilde{g}(\varepsilon; \kappa)  = \left\{(a_{22}-c_{22})(b_{10}-d_{10}) 
 -(a_{11}-c_{11})(b_{21}-d_{21}) \right\} + o(1)
$$
for small $\varepsilon > 0$.
The next lemma is a key for our purpose, which is proved in the appendix. 
\begin{lemma}\ \label{l3-6}
Define 
\begin{eqnarray*} 
\begin{array}{l}
f^*_u := \left\{
\begin{array}{l}
f_u(h^-(v^*),v^*), \ x \in [0,x_0], \\
f_u(h^+(v^*),v^*), \ x \in (x_0,1], 
\end{array} \right. \
f^*_v := \left\{
\begin{array}{l}
f_v(h^-(v^*),v^*), \ x \in [0,x_0], \\
f_v(h^+(v^*),v^*), \ x \in (x_0,1]. 
\end{array} \right. 
\end{array} 
\end{eqnarray*}
We have 
\begin{eqnarray*}  \left\{ 
\begin{array}{l}
a_{11}-c_{11} = \displaystyle \kappa \int^\infty_{-\infty} (\dot{ W}(z))^2 dz/
(\dot{ W}(0))^2 , \\[0.3cm]
a_{22}-c_{22}  = \left[ -(a_{11}-c_{11}) + \kappa (h^+(v^*) - h^-(v^*))/\dot{ W}(0)
\right]/D, \\[0.3cm]
b_{10}-d_{10} =  - \displaystyle \int^{h^+(v^*)}_{h^-(v^*)} f_v(u,v^*) du / \dot{ W}(0), \\[0.3cm]
b_{21}-d_{21} = \displaystyle  \left[ -(b_{10}-d_{10}) + \kappa \int^1_0 \left( 
\frac {f_u^* - f_v^*}{f_u^*} \right) dx \right]/D, \\[0.2cm]
\end{array} \right.
%\label{eqn:xax1}
\end{eqnarray*}
where $ W(z)$ is a solution of 
\begin{equation*} \left\{ 
\begin{array}{l}
\ddot{ W}(z) + f( W(z), v^*) = 0, \ z \in {\bf R}, \\[1ex]
 W(\pm \infty) = h^{\pm}(v^*), \  W(0) = \alpha; 
\end{array} \right. 
\end{equation*}
that is, 
$$
 W(z) =
 \left\{ 
\begin{array}{l}
\phi^-_0(z/x_0)+h^-(v^*), \ z \in (-\infty,0], \\[1ex]
\phi^+_0(z/(1-x_0))+h^+(v^*), \ z \in [0,\infty) 
\end{array}
\right.
$$ and $\alpha$ is a constant satisfying $h^-(v^*) < \alpha < h^+(v^*)$.
\end{lemma}

By this lemma, we can calculate
\begin{equation*} 
\begin{array}{lll}
\tilde{g}(0; \kappa^*) & = & \left\{(a_{22}-c_{22})(b_{10}-d_{10}) 
 -(a_{11}-c_{11})(b_{21}-d_{21}) \right\} \\[0.3cm]
 & = & \displaystyle -\kappa^* \left[ \kappa^* \int^\infty_{-\infty} (\dot{ W}(z))^2 dz \int^1_0 \left( \frac {f_u^* - f_v^*}{f_u^*} \right) dx \right. \\[0.5cm]
  & & \left. \displaystyle 
 + \ (h^+(v^*) - h^-(v^*))
 \int^{h^+(v^*)}_{h^-(v^*)} f_v(u,v^*) du \right] /\{ D(\dot{ W}(0))^2 \}. 
\end{array}
\end{equation*}
Therefore, we find two solutions of $\tilde{g}(0; \kappa^*) = 0$ such that (i) $\kappa^* = 0$ or  
$$
\mbox{(ii)} \hspace*{1cm}
\kappa^* = - \frac { \displaystyle (h^+(v^*) - h^-(v^*))
 \int^{h^+(v^*)}_{h^-(v^*)} f_v(u,v^*) du} {\displaystyle \int^\infty_{-\infty} (\dot{ W}(z))^2 dz \int^1_0 \left( \frac {f_u^* - f_v^*}{f_u^*} \right) dx} \ne 0. 
 \hspace*{4cm}
$$

For the case (ii), we can look for a single solution $\kappa(\ep) = \kappa^* \ep + o(\ep)$ satisfying $\tilde{g}(\varepsilon;\kappa(\ep)) = g(\ep;\ep \kappa(\ep)) = 0$. Hence, we see that \eqref{4c1} has
an eigenvalue $\lambda = \lambda(\ep) = \ep \kappa(\ep) \ne 0$ by Lemma~\ref{Evans}.

The transition layer solution can be denoted by $(u, v)(x;\ep;\xi)$ because
we can regard $\xi$ satisfying \eqref{a6} as a free parameter including it.
Then, we find that $({\partial u}/{\partial \xi},  {\partial v}/{\partial \xi})(x;\ep;\xi)$ is a solution of \eqref{4c1} with $\lambda = 0$, 
which corresponds to the case (i). Moreover, noting
$\int_0^1 \{{\partial u}/{\partial \xi}(x;\ep;\xi) + {\partial v}/{\partial \xi}(x;\ep;\xi)\} dx = 1 \neq 0 $ by \eqref{b1_1}, we find that $ \lambda = 0$ is not an eigenvalue of \eqref{4c1}.
 
\begin{remark}\ \label{rem3.1}
We can consider the case (i) as follows: 
When $\kappa^* =0$, we first construct $(p(x;\ep), q(x;\ep))$, the solution of \eqref{4c1} for $\lambda = \ep \kappa(\ep) = o(\ep)$. Solving the equation \eqref{eq001} with respect to $\alpha_i \ (i = 1,2,3,4)$, we obtain 
\begin{equation} \label{eq003} \left\{ 
\begin{array}{l}
\alpha_2 = \displaystyle - \alpha_1 \frac {\{p^-_x(x^*(\ep);\ep; o(\ep); 1,0) - p^+_x(x^*(\ep);\ep; o(\ep) ; 1,0)\}
}{\{p^-_x(x^*(\ep);\ep; o(\ep); 0,1) - p^+_x(x^*(\ep);\ep; o(\ep); 0,1)\}}, \\[0.3cm]
\alpha_3 = \alpha_1, \quad  \alpha_4 = \alpha_2. 
\end{array} \right. 
\end{equation}
Hence, $\alpha_i \ (i=2,3,4)$ are parametrized by $\alpha_1$, and we can assume that  $\alpha_1$ is independent of $\ep$ without loss of generality. Moreover, we note that 
\begin{equation*} \left\{ 
\begin{array}{l}
p^-_x(x^*(\ep);\ep; o(\ep); 1,0) - p^+_x(x^*(\ep);\ep; o(\ep); 1,0) 
 = (a_{11} - c_{11})\Big{|}_{\kappa=0} + o(1)  = o(1), \\[0.4cm]
p^-_x(x^*(\ep);\ep; o(\ep); 0,1) - p^+_x(x^*(\ep);\ep; o(\ep); 0,1) \\[0.2cm]
 = \dis\frac{1}{\ep} (b_{10} - d_{10} + o(1) ) = - \displaystyle \frac{1}{\ep} \left( \int^{h^+(v^*)}_{h^-(v^*)} f_v(u,v^*) du / \dot{ W}(0) + o(1) \right) \ne 0,
\end{array} \right. 
\end{equation*}
which implies that $\alpha_2 = o(\ep)$ as $\ep \to 0$. 
Then, it follows from
\eqref{eq002} that  $(p(x;\ep), q(x;\ep))$, 
the solution of \eqref{4c1} for $\lambda = \ep \kappa(\ep) = o(\ep)$, is represented by 
\begin{equation*}
\begin{array}{l}
\left[
\begin{array}{l}
p(x;\ep) \\
q(x;\ep)
\end{array} \right] = 
\left\{
\begin{array}{l}
\alpha_1 \left[ 
\begin{array}{l}
p^-(x;\ep;o(\ep) ; 1,0) \\
q^-(x;\ep;o(\ep) ; 1,0)
\end{array} \right] + \alpha_2 \left[ 
\begin{array}{l}
p^-(x;\ep;o(\ep) ; 0,1) \\
q^-(x;\ep;o(\ep) ; 0,1)
\end{array} \right] \\[0.4cm]
\hspace*{5cm}  ,x\in [0,x^*(\ep)], \\[1ex]
\alpha_1 \left[ 
\begin{array}{l}
p^+(x;\ep;o(\ep) ; 1,0) \\
q^+(x;\ep;o(\ep) ; 1,0)
\end{array} \right] + \alpha_2 \left[ 
\begin{array}{l}
p^+(x;\ep;o(\ep) ; 0,1) \\
q^+(x;\ep;o(\ep) ; 0,1)
\end{array} \right] \\[0.4cm]
\hspace*{5cm}  ,x\in [x^*(\ep),1].
\end{array}
\right.
\end{array}
\end{equation*}
Here, we note that the following:
$$
\displaystyle\int_0^{x^*(\ep)} p^-(x;\ep;o(\ep) ; 0,1)dx 
= - \frac {f_v^-}{f_u^-} \int_0^{x^*(\ep)} 1 \, dx + O(\ep) = 
- \frac {f_v^-}{f_u^-} x_0 +  O(\ep).
$$
Similarly, we have 
$$
\begin{array}{l}
\displaystyle \int^1_{x^*(\ep)} p^+(x;\ep;o(\ep) ; 0,1)dx 
= - \frac {f_v^+}{f_u^+}(1-  x_0) +  O(\ep), \\[0.4cm]
\displaystyle \int_0^{x^*(\ep)} q^-(x;\ep;o(\ep) ; 0,1)dx = x_0 + O(\ep), \\[0.4cm]
\displaystyle \int^1_{x^*(\ep)} q^+(x;\ep;o(\ep) ; 0,1)dx = (1 - x_0) + O(\ep). 
\end{array}
$$
By using $(p(x;\ep), q(x;\ep))$, we can calculate the following: 
\begin{equation*}
\begin{array}{lcl}
 0  & = & \displaystyle \frac{1}{\ep} \int_0^1 \left\{p(x;\ep) + q(x;\ep)\right\} dx \\[1ex]
 & = & \displaystyle \alpha_1 (h^+(v^*) - h^-(v^*))/\dot{ W}(0)  +  \frac{\alpha_2}{\ep}\left\{  \int^1_0 \left( \frac {f_u^* - f_v^*}{f_u^*} \right) dx + O(\ep) \right\}  \\[2ex]
  & = & \displaystyle \alpha_1 (h^+(v^*) - h^-(v^*))/\dot{ W}(0) + o(1). 
\end{array}
\end{equation*}
Since $\alpha_1$ is independent of $\ep$, 
we have $\alpha_1 (h^+(v^*) - h^-(v^*))/\dot{ W}(0) = 0$ by $\ep \to 0$, which leads to $\alpha_1 = 0$.
Therefore, we obtain $\alpha_i = 0 \ (i=1,2,3,4)$ by \eqref{eq003}, 
which implies that 
$\lambda = \ep \kappa(\ep) = o(\ep)$ is not an eigenvalue of \eqref{4c1}
because $(p(x;\ep), q(x;\ep)) \equiv (0, 0) $.
\end{remark} 

Thus, we obtain the following result:
\begin{theo}\ \label{l3-7} 
Assume that (A1)-(A4), and $\lambda = \lambda(\varepsilon) = O(\ep)$ as $\ep \to 0$. The eigenvalue problem \eqref{4c1} has only one eigenvalue 
\begin{equation} \label{eq:TH3.1}
\begin{array}{l}
\lambda(\ep)  = - \frac { \displaystyle (h^+(v^*) - h^-(v^*))
 \int^{h^+(v^*)}_{h^-(v^*)} f_v(u,v^*) du} {\displaystyle \int^\infty_{-\infty} (\dot{ W}(z))^2 dz \int^1_0 \left( \frac {f_u^* - f_v^*}{f_u^*} \right) dx} \ep + o(\ep) \in {\mathbb C}
\end{array}
\end{equation}
and the sign of the real part of $\lambda(\ep)$ is determined by 
$$
sign \{{\rm Re}(\lambda(\ep)) \} \ = \ sign \left\{ - \int^{h^+(v^*)}_{h^-(v^*)} f_v(u,v^*) du \right\} = sign \{ - J'(v^*) \}.
$$
\end{theo}
%
%%%%%%%%%%%%%%%%%%%%%%%%%%%%%%%%%%%%%%%%%%%%%%%%%%%%%%%%%%%%%%%%%%%%%%%%%%% 3-1
\subsection{Case $(\mathrm{II}) \ \ep \omega(\ep) \to 0$ and $\omega(\ep) \to \infty$ as $\ep \to 0$} \label{S3.2} %Section 3.2
%%%%%%%%%%%%%%%%%%%%%%%%%%%%%%%%%%%%%%%%%%%%%%%%%%%%%%%%%%%%%%%%%%%%%%%%

In case (II), $\lambda(\ep)$ satisfies 
$$ \lambda(\ep) = \ep \omega(\ep) \hat{\lambda}(\ep) \ \to \ 0, \ \mbox{as} \ \ep \to 0 $$
and $\hat{\lambda}(0) \ne 0$.
%
%%%%%%%%%%%%%%%%%%%%%%%%%%%%%%%%%%%%%%%%%%%%%%%%%%%%%%%%%%%%%%%%%%%
\subsubsection{Construction of $\bar{V}_1$ and $\bar{V}_2$}
%%%%%%%%%%%%%%%%%%%%%%%%%%%%%%%%%%%%%%%%%%%%%%%%%%%%%%%%%%%%%%%%%%%%
% Section 3.2.1
The essential part of the problem \eqref{4c3} for our purpose is to consider the following problem: 
\begin{equation}\label{eq008}
\left \{
\begin{array}{l}
\begin{array}{l}
\ep^2 p_{xx} +  f_u^\ep p + f_v^\ep q = \ep \omega(\ep)\hat{\lambda}(0) p, \\[1ex]
Dq_{xx}  - f_u^\ep p - f_v^\ep q = \ep \omega(\ep)\hat{\lambda}(0) q, 
\end{array} \quad  x \in (0,x^*(\ep)) \\[0.4cm]
(p_x, q_x)(0) = (0, 0), \ (p, q)(x^*(\ep)) = (a,b). 
\end{array} 
\right.
\end{equation}
Using the transformation $x=x^*(\ep)y$ in (\ref{eq008}), we have 
\begin{equation} \label{eq009}
\left\{
\begin{array}{l}
\begin{array}{l}
\varepsilon^2p_{yy}+x^*(\varepsilon)^2(f_u^\ep - \ep \omega(\ep)\hat{\lambda}(0))p 
+ x^*(\varepsilon)^2 f_v^\ep q = 0, \\[1ex]
D q_{yy}- x^*(\ep)^2f_u^\ep p - x^*(\ep)^2(f_v^\ep + \ep \omega(\ep)\hat{\lambda}(0))q = 0,
\end{array} y \in (0,1)\\[0.4cm]
(p_y,q_y)(0)=(0,0), \ (p,q)(1)= (a, b). 
\end{array}
\right.
\end{equation}
We first consider outer approximations of the form 
\begin{eqnarray*}
\left\{
\begin{array}{ll}
p^-(y) = P_0^-(y) + {\varepsilon}\omega(\ep)P_1^-(y) + o(\ep \omega(\ep)),\\[1ex]
q^-(y) = Q_0^-(y) + {\varepsilon}\omega(\ep)Q_1^-(y) + o(\ep \omega(\ep)).
\end{array}
\right.
\end{eqnarray*}
Substituting this into (\ref{eq009}), we equate the coefficients of the same powers of 
$\varepsilon$. \\[0.2cm]
$O((\ep \omega(\ep))^0)$:
\begin{eqnarray*}
\left\{
\begin{array}{l}
\begin{array}{l}
f_u^- P_0^- + f_v^- Q_0^- = 0, \\[1ex]
D Q_{0,yy}^- = 0, 
\end{array} y \in (0,1) \\[0.4cm]
Q_{0,y}^-(0)=0,\  Q_0^-(1)=b.
\end{array}
\right.
\end{eqnarray*}
$Q_0^-(y;b) = b$ and then $P_0^-(y;b) = - b f_v^-/f_u^- $, where $f_u^- = f_u(U_0^-, V_0^-)$
and $f_v^- = f_v(U_0^-, V_0^-)$. \\[0.2cm]
$O(\ep \omega(\ep))$:
\begin{eqnarray*}
\left\{
\begin{array}{l}
\begin{array}{l}
f_u^-P_1^- + f_v^- Q_1^-  -  \hat{\lambda}(0) P_0^-(y;b) = 0, \\[1ex]
D Q_{1,yy}^- + x_0^2 \hat{\lambda}(0) \{ f_v^-/f_u^- - 1\}Q_0^-(y;b) = 0, 
\end{array} y \in (0,1)\\[0.4cm]
Q_{1,y}^-(0)=0,\  Q_1^-(1)= 0. 
\end{array}
\right.
\end{eqnarray*}
We have 
\begin{eqnarray*}
\left\{
\begin{array}{l}
Q_1^-(y;\hat{\lambda}(0);b)= b x_0^2 \hat{\lambda}(0) \{ 1 - f_v^-/f_u^- \}(y^2-1)/(2D), \\[1ex]
P_1^-(y;\hat{\lambda}(0);b) = - \ f_v^- Q_1^-(y;\hat{\lambda}(0),b)/f_u^- + \hat{\lambda}(0)  P_0^-(y;b)/f_u^-.
\end{array}
\right.
\end{eqnarray*}
\indent Since the $p$ component does not satisfy the boundary condition at $y=1$, 
we have to modify this by adding correction terms $\rho_i^-, \pi_i^- \ (i=0,1)$ of the 
form 
\begin{eqnarray*}
\left\{
\begin{array}{lcl}
p^-(y) &=& P_0^-(y;b) + {\varepsilon}\omega(\ep)P_1^-(y;\hat{\lambda}(0);b) + \\[1ex]
 & & + \rho_{0}^-(\frac{y-1}{\varepsilon})+  {\varepsilon\omega(\ep)}\rho_1^-(\frac{y-1}{\varepsilon})+ o(\varepsilon\omega(\ep)),\\[1ex] 
q^-(y) &=& Q_0^-(y;b) + {\varepsilon}\omega(\ep) Q_1^-(y;\hat{\lambda}(0);b) \\[1ex]
 & & + \ \ep^2 \pi_0^-(\frac{y-1}{\varepsilon}) + \varepsilon^3\omega(\ep)\pi_1^-(\frac{y-1}{\varepsilon})+ o(\varepsilon^3\omega(\ep)).
\end{array}
\right.
\end{eqnarray*}
Substituting this into (\ref{eq009}) and using $z=(y-1)/\varepsilon$, we equate the coefficients of 
the same powers of $\varepsilon $. \\[0.2cm]
$O((\ep \omega(\ep))^{0}):$
\begin{eqnarray*}
\left\{
\begin{array}{l}
\begin{array}{l}
\ddot{\rho}_0^- + x_0^2\tilde{f}_u^-\rho_0^- = b x_0^2  \{\tilde{f}_u^- f_v^- - \tilde{f}_v^- f_u^- \}/f_u^-  , \\[1ex]
D \ddot{\pi}_0^- + \ddot{\rho}_0^-  =0, \\[1ex]
\end{array} \qquad z \in (-\infty ,0) \\[0.4cm]
\rho_0^-(-\infty)=0,\ \rho_0^-(0) = a + b f_v^-/f_u^-, \\[1ex]
\pi_0^-(-\infty)=0,\ \dot{\pi}_0^-(-\infty)=0, 
\end{array}
\right.
\end{eqnarray*}
where 
$\tilde{f}_u^- = f_u(U_0^-+\phi_0^-, V_0^-)$ and $\tilde{f}_v^- = f_v(U_0^-+\phi_0^-, V_0^-)$. 
By the first and third equations, and the form of \eqref{b21}, we have 
\begin{eqnarray*}
\begin{array}{l}
\rho_{0}^-(z;a,b)=( a + b f_v^-/f_u^-)\dot{\phi}_0^-(z)/\dot{\phi}_0^-(0) \\[0.2cm]
\hspace{1cm}- b x_0^2  \phi_0^-(z) \int_z^0 (\dot{\phi}_0^-(\zeta))^{-2} 
\int_{-\infty}^\eta \dot{\phi}_0^-(\eta) \{\tilde{f}_u^- f_v^- - \tilde{f}_v^- f_u^- \}/f_u^- d \eta d \zeta, 
\end{array}
\end{eqnarray*}
and then $\pi_0^-(z;a,b)=- \rho_0^-(z;a,b)/D$.  \\[0.2cm]
$O(\varepsilon\omega(\ep))$:
\begin{eqnarray*}
\left\{
\begin{array}{l}
\begin{array}{l}
\ddot{\rho}_1^- + x_0^2\tilde{f}_u^-\rho_1^- =   x_0^2 \hat{\lambda}(0) \rho_0^-(z;a,b) \\[1ex]
\hspace*{2cm} + \ b x_0^2 \hat{\lambda}(0)  (\tilde{f}_u^- f_v^- - f_u^- f_v^- )/(f_u^-)^2, \\[1ex]
D \ddot{\pi}_1^- + \ddot{\rho}_1^- = x_0^2 \hat{\lambda}(0) \rho_0^-(z;a,b), 
\\[1ex]
\end{array} z \in (-\infty ,0)\\[0.4cm]
\rho_1^-(-\infty)=0,\ \rho_1^-(0)= -P^-_1(1;\hat{\lambda}(0);b), \\[1ex]
\pi_1^-(-\infty)=0,\ \dot{\pi}_1^-(-\infty)=0.
\end{array}
\right.
%\label{eqn:BB4}
\end{eqnarray*}
We obtain 
\begin{eqnarray*}
\left\{
\begin{array}{l}
\rho_1^-(z;\hat{\lambda}(0);a,b) = -P^-_1(1;\hat{\lambda}(0);b)\dot{\phi}_0^-(z)/\dot{\phi}_0^-(0) \\[1ex]
\hspace*{2cm} - \ \dot{\phi}_0^-(z)\int_{z}^{0}(\dot{\phi}_0^-(\eta))^{-2}  \int_{-\infty}^{\eta}
\dot{\phi}_0^-(\zeta) \\[1ex]
\hspace*{2cm} \times \ x_0^2 \hat{\lambda}(0) \{ \rho_0^-(\eta;a,b) +  b (\tilde{f}_u^- f_v^- - f_u^- f_v^- )/(f_u^-)^2
\} d{\zeta}d{\eta}, \\[1ex]
\pi_1^-(z;\hat{\lambda}(0);a,b)= \{ -\rho_1^-(z;\hat{\lambda}(0);a,b) \\[1ex]
\hspace*{2cm}
+ \ x_0^2 \hat{\lambda}(0)\int_{-\infty}^{z}\int_{-\infty}^{\eta} 
\rho_0^-(\zeta;a,b) 
d{\zeta}d{\eta} \} /D.
\end{array} \right.
%\label{EQN:ddd1}
\end{eqnarray*}
Using the approximate functions defined by 
\begin{eqnarray*}
\left\{
\begin{array}{lcl}
P^-(y;\ep;\hat{\lambda}(0);a,b) &:=& P_0^-(y;b) + {\varepsilon}\omega(\ep)P_1^-(y;\hat{\lambda}(0);b) + \\[1ex]
 & & + \theta(y) \{ \rho_{0}^-(\frac{y-1}{\varepsilon};a,b)+  {\varepsilon\omega(\ep)}\rho_1^-(\frac{y-1}{\varepsilon};\hat{\lambda}(0);a,b)\},\\[1ex] 
Q^-(y;\ep;\hat{\lambda}(0);a,b) &:=& Q_0^-(y;b) + {\varepsilon}\omega(\ep) Q_1^-(y;\hat{\lambda}(0);b) \\[1ex]
 & & + \ \theta(y) \{\ep^2 (\pi_0^-(\frac{y-1}{\varepsilon};a,b) - \pi_0^-(0;a,b))\\[1ex]
 & &  + \ \varepsilon^3 \omega(\ep)(\pi_1^-(\frac{y-1}{\varepsilon};\hat{\lambda}(0);a,b) - \pi_1^-(0;\hat{\lambda}(0);a,b))\}, 
\end{array}
\right.
\end{eqnarray*}
we can obtain a solution $(p^-,q^-)(y;\ep;\hat{\lambda}(0);a,b)$ of \eqref{eq009} for any $\hat{\lambda}(0)\in{\mathbb C}$, which satisfies
$$
 ||(p^-, \, q^-) - (P^-, \, Q^-) ||_{X_\ep^-} = o(\varepsilon \omega(\ep)) \quad \mbox{as}  \ \ep \to 0.
$$
The proof is the same as that discussed in Subsection \ref{S3.1.1}. 

%%%%%%%%%%%%%%%%%%%%%%%%%%%%%%%%%%%%%%%%%%%%%%%%%%%%%%%%%%%%%%%
\subsubsection{Construction of $\bar{V}_3$ and $\bar{V}_4$}
%%%%%%%%%%%%%%%%%%%%%%%%%%%%%%%%%%%%%%%%%%%%%%%%%%%%%%%%%%%%%%%%%%%
% Section 3.2.2
Next, we consider the problem (\ref{4c4}) with $\lambda(\ep) = \ep \omega(\ep) \hat{\lambda}(\ep)$, in which $\hat{\lambda}(\ep)$ is replaced by $\hat{\lambda}(0)$. 
By using the transformation $x=x^*({\varepsilon}) + (1-x^*({\varepsilon}))y$, 
we have 
\begin{equation} \label{eqn:322_1}
\left\{
\begin{array}{l}
\begin{array}{l}
{\varepsilon}^2p_{yy}+(1-x^*({\varepsilon}))^2(f_u^\ep -\ep \omega(\ep)\hat{\lambda}(0))p + 
(1-x^*({\varepsilon}))^2f_v^\ep q = 0, \\[1ex]
Dq_{yy}-(1-x^*({\varepsilon}))^2f_u^\ep p 
-(1-x^*({\varepsilon}))^2(f_v^\ep+\ep \omega(\ep)\hat{\lambda}(0))q = 0, 
\end{array}
 y \in (0,1)\\[3ex]
(p,q)(0) = (a,b), \ (p_y,q_y)(1)=(0,0).
\end{array}
\right.
\end{equation}
First, we consider outer approximations of the form 
\begin{eqnarray*}
\left\{
\begin{array}{l}
p^+(y) = P_0^+(y) + \ep \omega(\ep) P_1^+(y) + o(\ep \omega(\ep)),\\[1ex]
q^+(y) = Q^+_0(y) + \ep \omega(\ep) Q_1^+(y) + o(\ep \omega(\ep)).
\end{array}
\right.
\end{eqnarray*}
Substituting this into (\ref{eqn:322_1}), we equate the coefficients of the same powers of 
$\varepsilon$ \\[0.2cm]
$O((\ep \omega(\ep))^0)$:
\begin{eqnarray*}
\left\{
\begin{array}{l}
\begin{array}{l}
f_u^+ P_0^+ + f_v^+ Q_0^+ = 0, \\[0.1cm]
D Q_{0,yy}^+ = 0,
\end{array} y \in (0,1 ) \\[0.3cm]
Q_0^+(0)=b, \ Q_{0,y}^+(1)=0.
\end{array} \right.
\end{eqnarray*}
We easily find that $Q_0^+(y;b) = b $ and $P_0^+(y;b) = - b f_v^+ /f_u^+$, where $f_u^+ = f_u(U_0^+, V_0^+)$ and $ f_v^+ = f_v(U_0^+, V_0^+)$.\\[0.2cm]
$O(\ep \omega(\ep))$:
\begin{eqnarray*}
\left\{
\begin{array}{l}
\begin{array}{l}
f_u^+ P_1^+ + f_v^+ Q_1^+ - \hat{\lambda}(0)P_0^+(y;b) = 0, \\[0.1cm]
D Q_{1,yy}^+ + (1-x_0)^2 \hat{\lambda}(0) \{ f_v^+/f_u^+ - 1\}Q_0^+(y;b)=0,
\end{array} y \in (0,1) \\[0.4cm]
Q_1^+(0)=0,\  Q_{1,y}^+(1)=0. 
\end{array}
\right.
\end{eqnarray*}
We have 
\begin{eqnarray*}
\left\{
\begin{array}{l}
Q_1^+(y;\hat{\lambda}(0),b)= b (1-x_0)^2 \hat{\lambda}(0) \{ -f_v^+/f_u^+ + 1\}(y^2-2y)/(2D), \\[1ex]
P_1^+(y;\hat{\lambda}(0);b) = - \ f_v^+ Q_1^+(y;\hat{\lambda}(0);b)/f_u^+ + \hat{\lambda}(0) P_0^+(y;b)/f_u^+. \\[1ex]

\end{array}
\right.
\end{eqnarray*}
\indent Since the $p$ component does not satisfy the boundary condition at $y=0$, we have to modify this by 
adding correction terms $\rho_i^+, \pi_i^+ \ (i=0,1)$ of the form 
\begin{eqnarray*}
\left\{
\begin{array}{lcl}
p^+(y) &=& P_0^+(y;b) + \ep \omega(\ep) P_1^+(y;\hat{\lambda}(0);b)  \\[0.1cm]
 & & + \ \rho_{0}^+(\frac{y}{\varepsilon}) + \ep \omega(\ep)\rho_1^+(\frac{y}{\varepsilon}) + o(\ep \omega(\ep)),\\[0.1cm] 
q^+(y) &=& Q_0^+(y;b) + \ep \omega(\ep) Q_1^+(y;\hat{\lambda}(0);b) \\[0.1cm]
 & & + \ {\varepsilon}^2\pi_0^+(\frac{y}{\varepsilon}) + {\varepsilon}^3 \omega(\ep)\pi_1^+(\frac{y}{\varepsilon}) + o(\ep^3 \omega(\ep)).
\end{array}
\right.
\end{eqnarray*}
Substituting this into (\ref{eqn:322_1}) and using $z=y/{\varepsilon}$, we equate the coefficient 
of the same powers of $\varepsilon$. \\[0.2cm]
$O((\ep \omega(\ep))^0)$:
\begin{eqnarray*}
\left\{
\begin{array}{l}
\begin{array}{l}
\ddot{\rho}_0^+ + (1-x_0)^2 \tilde{f}_u^+\rho_0^+= b (1-x_0)^2 \{\tilde{f}_u^+ f_v^+ - \tilde{f}_v^+ f_u^+ \}/f_u^+ , \\[0.1cm]
D \ddot{\pi}_0^+ + \ddot{\rho}_{0}^+ = 0,
\end{array} z \in (0 ,\infty)\\[0.3cm]
\rho_0^+(0)=a + b f_v^+/f_u^+, \ \rho_0^+(\infty)=0, \\[0.1cm]
\pi_0^+(\infty)=0, \ \dot{\pi}_0^+(\infty)=0, 
\end{array}
\right.
\end{eqnarray*}
where $\tilde{f}_u^+ = f_u(U_0^++\phi_0^+, V_0^+)$ and $\tilde{f}_v^+ = f_v(U_0^++\phi_0^+, V_0^+)$. 
By the first and third equations, and the form of \eqref{b32}, we have 
\begin{eqnarray*}
\begin{array}{l}
\rho_{0}^+(z;a,b)=(a + b f_v^+/f_u^+)\dot{\phi}_0^+(z)/\dot{\phi}_0^+(0) \\[1ex]
\hspace{1cm} - b (1-x_0)^2  \dot{\phi}_0^+(z) \int_0^z (\dot{\phi}_0^+(\eta))^{-2} \int_\eta^\infty  \dot{\phi}_0^+(\zeta) \{\tilde{f}_u^+ f_v^+ - \tilde{f}_v^+ f_u^+ \}/f_u^+d\zeta d \eta, 
\end{array}
\end{eqnarray*}
and then 
$\pi_0^+(z;a,b)=- \rho_0^+(z;a,b)/D$. \\[0.2cm]
$O(\ep \omega(\ep))$:
\begin{eqnarray*}
\left\{
\begin{array}{l}
\begin{array}{l}
 \ddot{\rho}_1^+ + (1-x_0)^2 \tilde{f}_u^+ \rho_1^+ =  (1-x_0)^2 \hat{\lambda}(0) \rho_0^+(z;a,b) \\[0.1cm]
\hspace{2cm}+ \ b (1-x_0)^2\hat{\lambda}(0)  (\tilde{f}_u^+ f_v^+ - f_u^+ f_v^+)/(f_u^+)^2, \\[0.1cm]
D \ddot{\pi}_1^+ + \ddot{\rho}_1^+ = (1-x_0)^2 \hat{\lambda}(0) \rho_0^+(z;a,b), 
\end{array}  z \in (0 ,\infty)\\[0.7cm]
 \rho_1^+(0) = -P^+_1(0;\hat{\lambda}(0);b), \ \rho_1^+(\infty)=0, \\[0.1cm]
 \pi_1^+(\infty)=0,\ \dot{\pi}_1^+(\infty)=0.
\end{array}
\right.
\end{eqnarray*}
We obtain 
\begin{eqnarray*}
\left\{ 
\begin{array}{l}
\rho_1^+(z;\hat{\lambda}(0);a,b) = -P^+_1(0;\hat{\lambda}(0);b)\dot{\phi}_0^+(z)/\dot{\phi}_0^+(0) - 
\dot{\phi}_0^+(z) \int_{0}^{z}(\dot{\phi}_0^+(\eta))^{-2}\int_{\eta}^{\infty}
\dot{\phi}_0^+(\zeta) \\[0.2cm]
\hspace*{0.5cm}\times  (1-x_0)^2 \hat{\lambda}(0)\{\rho_0^+(\eta;a,b) + b (\tilde{f}_u^+ f_v^+ - f_u^+ f_v^+)/(f_u^+)^2 \}d{\zeta}d{\eta}, \\[0.3cm]
\pi_1^+(z;\hat{\lambda}(0);a,b) = \{ -\rho_1^+(z;\hat{\lambda}(0);a,b) \\[0.2cm]
\hspace*{2.5cm} + (1-x_0)^2 \hat{\lambda}(0) 
\int_{z}^{\infty}\int_{\eta}^{\infty} \rho_0^+(\zeta;a,b) 
d{\zeta}d{\eta} \}/D. 
\end{array}
\right.
%\label{eqn:ddd5}
\end{eqnarray*}
Using the approximate functions defined by 
\begin{eqnarray*}
\left\{
\begin{array}{lcl}
P^+(y;\ep;\hat{\lambda}(0);a,b) &:=& P_0^+(y;b) + {\varepsilon}\omega(\ep)P_1^+(y;\hat{\lambda}(0);b) + \\[1ex]
 & & + \theta(1-y) \{ \rho_{0}^+(\frac{y}{\varepsilon};a,b)+  {\varepsilon\omega(\ep)}\rho_1^+(\frac{y}{\varepsilon};\hat{\lambda}(0);a,b)\},\\[1ex] 
Q^+(y;\ep;\hat{\lambda}(0);a,b) &:=& Q_0^+(y;b) + {\varepsilon}\omega(\ep) Q_1^+(y;\hat{\lambda}(0);b) \\[1ex]
 & & + \ \theta(1-y) \{\ep^2 (\pi_0^*(\frac{y}{\varepsilon};a,b) - \pi_0^+(0;a,b))\\[1ex]
 & &  + \ \varepsilon^3 \omega(\ep)(\pi_1^+(\frac{y}{\varepsilon};\hat{\lambda}(0);a,b) - \pi_1^+(0;\hat{\lambda}(0);a,b))\}, 
\end{array}
\right.
\end{eqnarray*}
we can obtain a solution $(p^+,q^+)(y;\ep;\hat{\lambda}(0);a,b)$ of \eqref{eqn:322_1} for any $\hat{\lambda}(0)\in {\mathbb C}$, which satisfies
$$
 ||(p^+, \, q^+) - (P^+,  \, Q^+) ||_{X_\ep^+} = o(\varepsilon \omega(\ep)) \quad \mbox{as}  \ \ep \to 0.
$$
The proof is the same as that discussed in Subsection \ref{S3.1.2}. 

%%%%%%%%%%%%%%%%%%%%%%%%%%%%%%%%%%%%%%%%%%%%%%%%%%%%%%%%%%%%%%%%%%%%%%%%%%%%%%%%%%%%%%%%%%%%%%%%%%%%%%%
\subsubsection{Evans function $g(\varepsilon;\ep \omega(\ep)\hat{\lambda}(\ep))$ corresponding to (\ref{4c1}) in  the case (II)}
%%%%%%%%%%%%%%%%%%%%%%%%%%%%%%%%%%%%%%%%%%%%%%%%%%%%%%%%%%%%%%%%%%%%%%%%%%%%%%%%%%%%%%%%%%%%%%%%%%%%%%%%
%Section 3.2.3

By using these $(p^-,q^-)(x;\varepsilon;\hat{\lambda}(0);a,b)$ and 
$(p^+,q^+)(x;\varepsilon;\hat{\lambda}(0);a,b)$, we can calculate 
$\bar{V}^{(i)}(\varepsilon;\hat{\lambda}(0)) := \bar{V}_i(x^*(\varepsilon);\varepsilon;\ep \omega(\ep)\hat{\lambda}(0))\ 
(i=1,2,3,4)$ (see \eqref{v12} and \eqref{v34}) as follows: 
$$
\begin{array}{l}
\bar{V}^{(1)}(\varepsilon;\hat{\lambda}(0))
=
\left[
\begin{array}{c}
1 \\
\tilde{a}_{10} + \ep \omega(\ep) \tilde{a}_{11}+o(\ep \omega(\ep))\\
0\\
o(\ep \omega(\ep))\\
\end{array}
\right], 
\end{array}
$$
$$
\begin{array}{l}
\bar{V}^{(2)}(\varepsilon;\hat{\lambda}(0))
=
\left[
\begin{array}{c}
0\\
\tilde{b}_{10} + \ep \omega(\ep) \tilde{b}_{11}+o(\ep \omega(\ep))\\
1\\
\ep \omega(\ep) \tilde{b}_{21} +o(\ep \omega(\ep))\\
\end{array}
\right],
\end{array}
$$
$$
\begin{array}{l}
\bar{V}^{(3)}(\varepsilon;\hat{\lambda}(0))
=
\left[
\begin{array}{c}
1 \\
\tilde{c}_{10} + \ep \omega(\ep) \tilde{c}_{11}+o(\ep \omega(\ep))\\
0\\
o(\ep \omega(\ep))\\
\end{array}
\right],
\end{array}
$$
$$
\begin{array}{l}
\bar{V}^{(4)}(\varepsilon;\hat{\lambda}(0))
=
\left[
\begin{array}{c}
0\\
\tilde{d}_{10} + \ep \omega(\ep) \tilde{d}_{11}+o(\ep \omega(\ep))\\
1\\
\ep \omega(\ep) \tilde{d}_{21} + o(\ep \omega(\ep))\\
\end{array}
\right], 
\end{array} 
$$
where 
\begin{eqnarray*}
\begin{array}{l}
\tilde{a}_{10} = \dot{\rho}^-_0(0;1,0)/x_0, \quad \tilde{a}_{11} = \dot{\rho}^-_1(0;\hat{\lambda}(0);1,0)/x_0, \quad \tilde{b}_{10} = \dot{\rho}^-_0(0;0,1)/x_0, \\[0.1cm]
\tilde{b}_{11} = \dot{\rho}^-_1(0;\hat{\lambda}(0);0,1)/x_0, \quad \tilde{b}_{21} = Q^-_{1,y}(1;\hat{\lambda}(0);1)/x_0, \\[0.1cm]
\tilde{c}_{10} = \dot{\rho}^+_0(0;1,0)/(1-x_0), \quad \tilde{c}_{11} = \dot{\rho}^+_1(0;\hat{\lambda}(0);1,0)/(1-x_0), \\[0.1cm]
\tilde{d}_{10} = \dot{\rho}^+_0(0;0,1)/(1-x_0), \quad \tilde{d}_{11} = \dot{\rho}^+_1(0;\hat{\lambda}(0);0,1)/(1-x_0), \\[0.1cm]
\tilde{d}_{21} = Q^+_{1,y}(0;\hat{\lambda}(0);1)/(1-x_0). 
\end{array}
\end{eqnarray*}
By a simple calculation, noting that $\tilde{a}_{10} \ = \ \tilde{c}_{10}$ (see Subsection \ref{S3.1.3}), 
we find that $g(\varepsilon;\ep \omega(\ep)\hat{\lambda}(0))$ is represented as 
follows:
\begin{eqnarray*}
\begin{array}{lcl}
g(\varepsilon;\ep \omega(\ep)\hat{\lambda}(0)) & = &
\det [\bar{V}^{(1)}(\varepsilon;\hat{\lambda}(0)), \bar{V}^{(2)}(\varepsilon;\hat{\lambda}(0)), \bar{V}^{(3)}(\varepsilon;\hat{\lambda}(0)), \bar{V}^{(4)}(\varepsilon;\hat{\lambda}(0)) ] \\[1ex]
 & = & (\ep \omega(\ep))^2 \left\{(\tilde{a}_{11}-\tilde{c}_{11})(\tilde{d}_{21}-\tilde{b}_{21}) + o(1)\right\} . 
\end{array}
\end{eqnarray*}
Moreover, a calculation  similar to that in Lemma \ref{l3-6} leads to 
\begin{eqnarray*} 
\begin{array}{l}
\tilde{a}_{11}-\tilde{c}_{11} = \hat{\lambda}(0) \displaystyle \int^\infty_{-\infty} (\dot{ W}(z))^2 dz/
(\dot{ W}(0))^2 , \\[0.3cm]
\tilde{d}_{21}-\tilde{b}_{21} = \displaystyle \frac { \hat{\lambda}(0)}D  \int^1_0 \left( 
\frac {f_u^* - f_v^*}{f_u^*} \right) dx. 
\end{array} 
\end{eqnarray*}
Hence, recalling $ \hat{\lambda}(0) \ne 0$, we see that 
$$
g(\varepsilon;\ep \omega(\ep)\hat{\lambda}(0)) \ne 0,
$$
which implies that $g(\varepsilon;\ep \omega(\ep)\hat{\lambda}(\ep)) \ne 0$ for small $\ep > 0$. Thus, $\lambda(\ep) = \ep \omega(\ep)\hat{\lambda}(\ep) $ of the case (II) is not an eigenvalue of \eqref{4c1}.

%%%%%%%%%%%%%%%%%%%%%%%%%%%%%%%%%%%%%%%%%%%%%%%%%%%%%%%%%%%%%%%%%%%%%%%%%%% 3-3 
\subsection{Case $(\mathrm{III}) \ \ep \omega(\ep) \to \omega_0$ as $\ep \to 0$ for some positive constant $\omega_0$} \label{S3.3} %Section 3.3.1 
For case (III), we already showed that there exist no eigenvalues of \eqref{4c1} in our previous paper \cite{KTI} using a different method. Only the following essential parts are described for self-completion of this paper. 
In case (III), $\lambda(\ep)$ satisfies 
$$ \lambda(\ep) = \ep \omega(\ep) \hat{\lambda}(\ep) \ \to \ \omega_0 \hat{\lambda}(0) \ne 0, \ \mbox{as} \ \ep \to 0 $$
and $\omega_0 > 0$.
We assume that ${\rm Re}\hat{\lambda}(0) \geq 0$  and $\hat{\lambda}(0) \ne 0$. 
%%%%%%%%%%%%%%%%%%%%%%%%%%%%%%%%%%%%%%%%%%%%%%%%%%%%%%%%%%%%%%%%%%%%%%%
\subsubsection{Construction of $\bar{V}_1, \bar{V}_2, \bar{V}_3$ and $\bar{V}_4$}
%%%%%%%%%%%%%%%%%%%%%%%%%%%%%%%%%%%%%%%%%%%%%%%%%%%%%%%%%%%%%%%%%%%%%%%
%
The problem \eqref{4c3} is approximately reduced to 
\begin{equation*}%\label{eq0008}
\left \{
\begin{array}{l}
\begin{array}{l}
\ep^2 p_{xx} +  f_u^\ep p + f_v^\ep q = \omega_0 \hat{\lambda}(0) p, \\[1ex]
Dq_{xx}  - f_u^\ep p - f_v^\ep q = \omega_0 \hat{\lambda}(0) q, 
\end{array} \quad  x \in (0,x^*(\ep)) \\[0.4cm]
(p_x, q_x)(0) = (0, 0), \ (p, q)(x^*(\ep)) = (a, b).
\end{array} 
\right.
\end{equation*}
Using $x = x^*(\ep)y$ and $z = (y-1)/\ep$, we see that 
the equations for the calculation of the Evans function consist of the following terms: 
\begin{eqnarray*}
\left\{
\begin{array}{ll}
p^-(x;\ep;\hat{\lambda}(0);a,b) = P_0^-(\frac x{x^*(\ep)};\hat{\lambda}(0);b) +  \theta(\frac x{x^*(\ep)}) \rho_{0}^-(\frac{x-x^*(\ep)}{\varepsilon x^*(\ep)};\hat{\lambda}(0);a,b) + O(\ep), \\[1ex] 
q^-(x;\ep;\hat{\lambda}(0);a,b) = Q_0^-(\frac x{x^*(\ep)};\hat{\lambda}(0);b) + O(\ep), 
\end{array}
\right.
\end{eqnarray*}
where $P_0^-, Q_0^-$ and $\rho_{0}^-$ satisfy the following equations:
\begin{equation}\label{eq1110}
\left\{
\begin{array}{l}
f_u^- P_0^- + f_v^- Q_0^- = \omega_0 \hat{\lambda}(0)P_0^- ,\\[0.2cm]
D Q_{0,yy}^- - x_0^2 g^- Q_0^- = 0, \ y \in (0,1) \\[0.2cm]
Q_{0,y}^-(0)=0,\  Q_0^-(1)=b, 
\end{array}
\right.
\end{equation}
\begin{equation} \label{eq1111}
\left\{
\begin{array}{l}
\ddot{\rho}_0^- + x_0^2(\tilde{f}_u^- - \omega_0 \hat{\lambda}(0))\rho_0^- \\[0.2cm]
= \ b x_0^2 \{ (\tilde{f}_u^- - f_u^-)f_v^-/(f_u^- - \omega_0 \hat{\lambda}(0)) 
- (\tilde{f}_v^- - f_v^-) \}, \ z \in (-\infty ,0) \\[0.2cm]
\rho_0^-(-\infty)=0,\ \rho_0^-(0) = a,
\end{array}
\right.
\end{equation}
respectively, and $g^- := \omega_0 \hat{\lambda}(0)\{f_u^- - f_v^- - \omega_0 \hat{\lambda}(0)\}/(f_u^- - \omega_0 \hat{\lambda}(0))$. The solutions of \eqref{eq1110} are solved as 
\begin{equation} \label{eqa122}
Q_0^-(y;\hat{\lambda}(0);b) = b \cdot \frac {e^{x_0 \sqrt{g^-/D}y} + e^{-x_0 \sqrt{g^-/D}y}}{e^{x_0 \sqrt{g^-/D}}+e^{-x_0 \sqrt{g^-/D}}}
\end{equation}
and 
$$
P_0^-(y;\hat{\lambda}(0);b) = - f_v^- Q_0^-(y;\hat{\lambda}(0);b)/(f_u^- - \omega_0 \hat{\lambda}(0)).
$$
For the solution $\rho_0^-(z;\hat{\lambda}(0);a,b)$ of \eqref{eq1111}, we can show the 
following lemma: 
\begin{lemma}\label{lem1c2} 
For any $\hat{\lambda}(0) \in {\bf \bar{C}}_0 := \{ \nu \in {\mathbb C} \ | \ {\rm Re}\{\nu\} \geq 0\}$, \eqref{eq1111} has a unique solution $\rho_0^-(z;\hat{\lambda}(0);a,b)$. 
\end{lemma}

Since this lemma is proved by the argument similar to \cite[Lemma 4.1]{IM}, we omit it.

Next, the problem \eqref{4c4} is approximately reduced to 
\begin{equation*}%\label{eq0009}
\left \{
\begin{array}{l}
\begin{array}{l}
\ep^2 p_{xx} +  f_u^\ep p + f_v^\ep q = \omega_0 \hat{\lambda}(0) p, \\[1ex]
Dq_{xx}  - f_u^\ep p - f_v^\ep q = \omega_0 \hat{\lambda}(0) q, 
\end{array} \quad  x \in (x^*(\ep),1) \\[0.4cm]
(p, q)(x^*(\ep)) = (a ,b), \  (p_x, q_x)(1) = (0, 0).
\end{array} 
\right.
\end{equation*}
Using $x =  x^*(\ep) + (1-x^*(\ep))y$ and $z = y/\ep$, we see that 
the equations for the calculation of the Evans function consist of the following terms: 
\begin{eqnarray*}
\left\{
\begin{array}{ll}
p^+(x;\ep;\hat{\lambda}(0);a,b) = P_0^+(\frac {x-x^*(\ep)}{1-x^*(\ep)};\hat{\lambda}(0);b) + \theta(\frac {1-x }{1-x^*(\ep)}) \rho_{0}^+(\frac{x-x^*(\ep)}{\varepsilon (1-x^*(\ep))};\hat{\lambda}(0);a,b) \\[1ex]
\hspace{3cm} + \ O(\ep), \\[1ex] 
q^+(x;\ep;\hat{\lambda}(0);a,b) = Q_0^+(\frac {x-x^*(\ep)}{1-x^*(\ep)};\hat{\lambda}(0);b) + O(\ep), 
\end{array}
\right.
\end{eqnarray*}
where 
$P_0^+, Q_0^+$ and $\rho_{0}^+$ satisfy the following equations:
\begin{equation}\label{eq1112}
\left\{
\begin{array}{l}
f_u^+ P_0^+ + f_v^+ Q_0^+ = \omega_0 \hat{\lambda}(0)P_0^+ ,\\[0.2cm]
D Q_{0,yy}^+ - (1-x_0)^2 g^+ Q_0^+ = 0, \ y \in (0,1) \\[0.2cm]
Q_0^+(0)=b,\  Q_{0,y}^+(1)=0 
\end{array}
\right.
\end{equation}
and
\begin{equation} \label{eq1113}
\left\{
\begin{array}{l}
\ddot{\rho}_0^+ + (1-x_0)^2(\tilde{f}_u^+ - \omega_0 \hat{\lambda}(0)) \rho_0^+ \\[0.2cm]
= \ b (1-x_0)^2 \{ (\tilde{f}_u^+ - f_u^+)f_v^+/(f_u^+ - \omega_0 \hat{\lambda}(0)) 
- (\tilde{f}_v^+ - f_v^+) \}, \ z \in (0,\infty) \\[0.2cm]
\rho_0^+(0) = a, \ \rho_0^+(\infty)=0, 
\end{array}
\right.
\end{equation}
respectively, and $g^+ := \omega_0 \hat{\lambda}(0)\{f_u^+ - f_v^+ - \omega_0 \hat{\lambda}(0)\}/(f_u^+ - \omega_0 \hat{\lambda}(0))$. The solutions of \eqref{eq1112} are solved as 
\begin{equation} \label{eqa123}
\begin{array}{l}
Q_0^+(y;\hat{\lambda}(0);b) = b \cdot \displaystyle 
\frac {e^{(1-x_0) \sqrt{g^+/D}(1-y)} + e^{-(1-x_0) \sqrt{g^+/D}(1-y)}}{e^{(1-x_0) \sqrt{g^+/D}}+e^{-(1-x_0) \sqrt{g^+/D}}}
\end{array}
\end{equation}
and 
$$
P_0^+(y;\hat{\lambda}(0);b) = - f_v^+ Q_0^+(y;\hat{\lambda}(0);b)/(f_u^+ - \omega_0 \hat{\lambda}(0)).
$$
Similarly to \eqref{eq1111}, we can show \eqref{eq1113} has a unique solution $\rho_0^+(z;\hat{\lambda}(0);a)$ for any $\hat{\lambda}(0) \in {\bf \bar{C}}_0$. 
%
%%%%%%%%%%%%%%%%%%%%%%%%%%%%%%%%%%%%%%%%%%%%%%%%%%%%%%%%%%%%%%%%%%%%%%%%%%%%%%%%%%%%%%%%%%%%%%%%%%%%%%%
\subsubsection{Evans function $g(\varepsilon;\ep \omega(\ep)\hat{\lambda}(\ep))$ corresponding to (\ref{4c1}) in  the case (III)}
%******************************************************************************
%Section 3.3.2

The Evans function $g(\varepsilon;\ep \omega(\ep)\hat{\lambda}(\ep))$ is calculated in the following form: 
\begin{eqnarray*}
\begin{array}{l}
g(\varepsilon;\ep \omega(\ep)\hat{\lambda}(\ep)) = \{Q_{0,y}^-(1;\hat{\lambda}(0);1)/x_0 - Q_{0,y}^+(0;\hat{\lambda}(0);1)/(1-x_0)\} \\[1ex]
\hspace{2cm} \ \times \ \{\dot{\rho}_0^+(0;\hat{\lambda}(0);1,0)/(1-x_0) - \dot{\rho}_0^-(0;\hat{\lambda}(0);1,0)/x_0\} + O(\ep)
\end{array}
\end{eqnarray*}
for small $\ep > 0$.

First, we show that for any $\hat{\lambda}(0) \in {\bf \bar{C}}_0 \backslash \{0\}$ 
\begin{equation} \label{eq222}
\begin{array}{l}
H_1(\hat{\lambda}(0)) := Q_{0,y}^-(1;\hat{\lambda}(0);1)/x_0 - Q_{0,y}^+(0;\hat{\lambda}(0);1)/(1-x_0) \ne 0.
\end{array}
\end{equation}
We define $A(y) := {\rm Re}\{Q_0^-(y;\hat{\lambda}(0);1)\}, B(y) := {\rm Im}\{Q_0^-(y;\hat{\lambda}(0);1)\}$, $\alpha := {\rm Re}\{g^-\}$ and $\beta := {\rm Im}\{g^-\}$, and rewrite the second and third equations of \eqref{eq1110} with $ b=1$ as 
\begin{equation}\label{eq1116}
\left\{
\begin{array}{l}
D A_{yy} - x_0^2 \alpha A + x_0^2 \beta B = 0, \ y \in (0,1) \\[0.2cm]
A_y(0)=0,\  A(1)=1
\end{array}
\right.
\end{equation}
and 
\begin{equation}\label{eq1117}
\left\{
\begin{array}{l}
D B_{yy} - x_0^2 \alpha B - x_0^2 \beta A = 0, \ y \in (0,1) \\[0.2cm]
B_y(0)=0,\  B(1)=0. 
\end{array}
\right.
\end{equation}
Multiplying $B(y)$ (resp. $A(y)$) to the first equation of \eqref{eq1116} (resp. \eqref{eq1117}) and integrating them on $y \in [0,1]$, we have $B_y(1) = x_0^2 \beta 
\int_0^1 \{A(y)^2 + B(y)^2\} dy/D$. That is, 
$$
{\rm Im}\{Q_{0,y}^-(1;\hat{\lambda}(0);1)\} = x_0^2 \, {\rm Im}\{g^-\} I^-, 
$$
where $I^- := \int_0^1 \{({\rm Re}\{Q_0^-(y;\hat{\lambda}(0);1)\})^2 + ({\rm Im}\{Q_0^-(y;\hat{\lambda}(0);1)\})^2 \} dy/D > 0 $. 
Similarly, by the second and third equations of \eqref{eq1112} with $b=1$, we have 
$$
{\rm Im}\{Q_{0,y}^+(0;\hat{\lambda}(0);1)\}  = - (1-x_0)^2 {\rm Im}\{g^+\} I^+,
$$
where $I^+ := \int_0^1 \{({\rm Re}\{Q_0^+(y;\hat{\lambda}(0);1)\})^2 + ({\rm Im}\{Q_0^+(y;\hat{\lambda}(0);1)\})^2\} dy/D > 0 $. 
Furthermore, we note that 
$$
{\rm Im}\{g^{\pm} \} = {\rm Im}\{\hat{\lambda}(0)\} S^{\pm}(\hat{\lambda}(0)), 
$$
where 
$$
S^{\pm}(\hat{\lambda}(0)) = \omega_0 \cdot 
\frac {f_u^{\pm}(f_u^{\pm}-f_v^{\pm})-2\omega_0{\rm Re}\{\hat{\lambda}(0)\}f_u^{\pm}
+(\omega_0|\hat{\lambda}(0)|)^2}{(f_u^{\pm}-\omega_0{\rm Re}\{\hat{\lambda}(0)\})^2 + (\omega_0{\rm Im}\{\hat{\lambda}(0)\})^2}. 
$$
From the assumptions (A1) and (A3), we find that 
$S^{\pm}(\hat{\lambda}(0)) > 0$ for any $\hat{\lambda}(0) \in {\bf \bar{C}}_0$. Thus, we obtain 
$$
{\rm Im}H_1(\hat{\lambda}(0))  = {\rm Im}\{\hat{\lambda}(0)\} 
  \ [x_0 S^-(\hat{\lambda}(0)) I^- 
 + \ (1-x_0) S^+(\hat{\lambda}(0)) I^+ ], 
$$
which implies that ${\rm Im}H_1(\hat{\lambda}(0)) \ne 0$ for any $\hat{\lambda}(0) \in {\bf \bar{C}}_0$ when ${\rm Im}\{\hat{\lambda}(0)\} \ne 0$. 
Then, we suppose that ${\rm Im}\{\hat{\lambda}(0)\} = 0$ and $\hat{\lambda}(0) > 0$. Noting $g^{\pm} > 0$, from the representation of \eqref{eqa122} and \eqref{eqa123}, 
we have 
$$
H_1(\hat{\lambda}(0)) = \  \sqrt{\frac{g^-}{D}}\tanh \left(x_0 \sqrt{\frac{g^-}{D}}\right) 
\ + \ \sqrt{\frac{g^+}{D}}\tanh \left((1-x_0) \sqrt{\frac{g^+}{D}}\right) > 0
$$
when $\hat{\lambda}(0) > 0$. Therefore, we see that $H_1(\hat{\lambda}(0)) \ne 0$. 

Next, we show that $H_2(\hat{\lambda}(0)) := \dot{\rho}_0^+(0;\hat{\lambda}(0);1,0)/(1-x_0) - \dot{\rho}_0^-(0;\hat{\lambda}(0);1,0)/x_0 \ne 0$. Applying an argument as applied to \eqref{eq222}, we can 
show that ${\rm Im}H_2(\hat{\lambda}(0)) \ne 0$ when ${\rm Im}\{\hat{\lambda}(0)\} \ne 0$. Then, we suppose that $\hat{\lambda}(0)$ is real. 
Setting $\hat{\lambda}(0)=\mu, a =1$ and $b=0$ in \eqref{eq1111}, we obtain 
\begin{equation} \label{eq1122}
\left\{
\begin{array}{l}
\ddot{\rho}_0^- + x_0^2(\tilde{f}_u^- - \omega_0 \mu)\rho_0^- = 0, \ z \in (-\infty ,0) \\[0.2cm]
\rho_0^-(-\infty)=0,\ \rho_0^-(0) = 1. 
\end{array}
\right.
\end{equation}
Differentiating \eqref{eq1122} by $\mu$ and putting $w(z) = {\partial \rho_0^-}/{\partial \mu} (z;\mu;1,0)$, we have 
\begin{equation} \label{eq1123}
\left\{
\begin{array}{l}
\ddot{w} + x_0^2(\tilde{f}_u^- - \omega_0 \mu)w = x_0^2\omega_0\rho_0^- , \ z \in (-\infty ,0) \\[0.2cm]
w(-\infty)=0,\ w(0) = 0. 
\end{array}
\right.
\end{equation}
Multiplying the first equation of \eqref{eq1123} by $\rho_0^-$ and integrating it on $z \in (-\infty,0)$,
we obtain 
$$
\frac \partial{\partial \mu} \dot{\rho}_0^-(0;\mu;1,0) = \dot{w}(0) = x_0^2 \omega_0
\int_{-\infty}^0 (\rho_0^-(z))^2 dz
$$
by using integration by parts.
Similarly, we have 
$$
\frac \partial{\partial \mu} \dot{\rho}_0^+(0;\mu;1,0) = -(1-x_0)^2 \omega_0
\int^{\infty}_0 (\rho_0^+(z))^2 dz.
$$
These equations imply that $\partial H_2(\hat{\lambda}(0))/{\partial \hat{\lambda}(0)} < 0$.
On the other hand, noting $\rho_0^\pm(z;0;1,0) = \dot{\phi}_0^\pm(z)/\dot{\phi}_0^\pm(0)$ by
\eqref{b21} and \eqref{b32}, we have 
$$
H_2(0) = \ddot{\phi}^+_0(0)/((1-x_0) \dot{\phi}^+_0(0)) -  \ddot{\phi}^-_0(0)/(x_0 \dot{\phi}^-_0(0))  =0,
$$
where we used 
$(1-x_0)\dot{\phi}_0^-(0) = x_0\dot{\phi}_0^+(0)$ by \eqref{b00009}, $\ddot{\phi}_0^-(0) =- x_0^2 f(\alpha,v^*)$
by \eqref{b18}, and $\ddot{\phi}_0^+(0) = -(1-x_0)^2 f(\alpha,v^*)$ by \eqref{b29}.
Hence, we see that $H_2(\hat{\lambda}(0)) \ne 0$ for any $\hat{\lambda}(0) \in {\bf \bar{C}}_0 \backslash \{0\}$. 
Therefore, $g(\varepsilon;\ep \omega(\ep)\hat{\lambda}(\ep))  \ne 0$ for any $\lambda(\ep) = \ep \omega(\ep) \hat{\lambda}(\ep)$ of the Case (III). That is, $\lambda(\ep) = \ep \omega(\ep) \hat{\lambda}(\ep)$ of the Case (III) is not an eigenvalue of \eqref{4c1}.

%%%%%%%%%%%%%%%%%%%%%%%%%%%%%%%%%%%%%%%%%%%%%%%%%%%%%%%%%%%%%%%%%%%%%%%%%%%%%%%%%%% Section 3.4
\subsection{Distribution of eigenvalues of (\ref{4c1})}
%%%%%%%%%%%%%%%%%%%%%%%%%%%%%%%%%%%%%%%%%%%%%%%%%%%%%%%%%%%%%%%%%%%%%%%%%%%
In Sections \ref{S3.1}, \ref{S3.2} and \ref{S3.3}, we find that eigenvalues $\lambda \in \mathbb{C}_d :=  
\{\lambda \in {\mathbb C}\ | \ {\rm Re}\lambda > -d \ep \} $ for any fixed $d>0$ are determined by \eqref{eq:TH3.1}. 
Then, we provide the following result on the stability of a single transition layer solution of \eqref{b1} with \eqref{b1_1}: 

\begin{theo}\ \label{l3-8}
Under the assumptions (A1)-(A4), for any fixed $d>0$  
the eigenvalue problem \eqref{4c1} has only one eigenvalue 
$$
\lambda(\ep)  = - \frac { \displaystyle (h^+(v^*) - h^-(v^*))
 \int^{h^+(v^*)}_{h^-(v^*)} f_v(u,v^*) du} {\displaystyle \int^\infty_{-\infty} (\dot{ W}(z))^2 dz \int^1_0 \left( \frac {f_u^* - f_v^*}{f_u^*} \right) dx} \ep + o(\ep) 
$$
in $\mathbb{C}_d$ and the sign of the real part of $\lambda(\ep)$ is determined by 
$$
sign \{{\rm Re}(\lambda(\ep)) \} \ = \ sign \left\{ - \int^{h^+(v^*)}_{h^-(v^*)} f_v(u,v^*) du \right\} = sign \{ - J'(v^*) \}.
$$
Then, the single transition layer solution $(u,v)(x;\ep)$ is stable when $J'(v^*) > 0$, conversely it is unstable when  $J'(v^*) < 0$. 
\end{theo}

\begin{remark}\ \label{rem3.2}
For examples of both cases $J'(v^*) > 0$ and $J'(v^*) < 0$, one can refer to our previous paper \cite[Section 4]{KTI}. There, we numerically showed that the single transition layer solution is stable (resp. unstable) when $J'(v^*) > 0$ (resp. $J'(v^*) < 0$). 
\end{remark}

\vspace{0.5cm}
\noindent
{\large \bf Acknowledgments.} 
We would like to thank anonymous referees for their valuable comments and suggestions to improve the original manuscript. The first author was supported in part by the JSPS KAKENHI Grant Numbers JP19K03618 and JP24K06864. The second author was also supported in part by the JSPS KAKENHI Grant Number JP24K06845.

\appendix
%%%%%%%%%%%%%%%%%%%%%%%%%%%%%%%%%%%%%%%%%%%%%%%%%%%%%%%%%%%%%%%%%%
%\section{Appendix}
%%%%%%%%%%%%%%%%%%%%%%%%%%%%%%%%%%%%%%%%%%%%%%%%%%%%%%%%%%%%%%%%%%%%
%\def\thesection{Appendix \Alph{section}}
\appendix
\numberwithin{equation}{section}
\makeatletter 
% "activate" the preparatory code, but for section-level headers only
\newcommand{\section@cntformat}{Appendix \thesection:\ }
\makeatother
\section{ Proof of Lemma \ref{l3-6}}
Differentiating the first and third equations of \eqref{b19} with $U_1^- = 0 = V_1^-$ by $z$, $Y(z) := \dot{\phi}_1^-(z)$ satisfies 
\begin{eqnarray*} \left\{ 
\begin{array}{l}
\ddot{Y} + x_0^2 \tilde{f}_u^- Y = -\{2 x_0 x_1  \tilde{f}_u^- + x_0^2  \tilde{f}_{uu}^- \phi_1^- \} \dot{\phi}_0^-, \ z \in (-\infty,0) \\[0.2cm]
Y(-\infty) = 0, \ Y(0) = \dot{\phi}_1^-(0).
\end{array} \right.
\end{eqnarray*}
In the same way as the solution of \eqref{b19} is represented by \eqref{b21}, noting \eqref{b18}, we have 
\begin{eqnarray*} 
\begin{array}{lll}
 Y(z) &=& \dot{\phi}_1^-(0) \dot{\phi}_0^-(z)/\dot{\phi}_0^-(0) \\[0.2cm]
 & & + \ \dot{\phi}_0^-(z) \int_z^0 (\dot{\phi}_0^-(\eta))^{-2} \int_{-\infty}^\eta 
 (\dot{\phi}_0^-(\zeta))^2 \{2 x_0 x_1  \tilde{f}_u^- + x_0^2  \tilde{f}_{uu}^- \phi_1^-\} d\zeta d \eta
\end{array}
\end{eqnarray*}
and 
\begin{eqnarray*} %\label{eqap1}
\begin{array}{lll}
\ddot{\phi}_1^-(0) = \dot{Y}(0) = \dot{\phi}_1^-(0) \ddot{\phi}_0^-(0)/\dot{\phi}_0^-(0) \\[0.2cm]
  - \ \int_{-\infty}^0 
 (\dot{\phi}_0^-(\zeta))^2 \{2 x_0 x_1  \tilde{f}_u^- + x_0^2  \tilde{f}_{uu}^- \phi_1^-\}d\zeta/\dot{\phi}_0^-(0). 
\end{array}
\end{eqnarray*}
Similarly, it follows from \eqref{b30} that 
\begin{eqnarray*} %\label{eqap2}
\begin{array}{l}
\ddot{\phi}_1^+(0)  =  \dot{\phi}_1^+(0) \ddot{\phi}_0^+(0)/\dot{\phi}_0^+(0) \\[0.2cm]
  - \ \int^{\infty}_0 
 (\dot{\phi}_0^+(\zeta))^2 \{2 (1-x_0) x_1  \tilde{f}_u^+ - (1-x_0)^2  \tilde{f}_{uu}^+ \phi_1^+\}d\zeta/\dot{\phi}_0^+(0). 
\end{array}
\end{eqnarray*}
From \eqref{eqn:C4}, \eqref{C6} and $\dot{\rho}_0^{\pm}(z;1,0) = \dot{\phi}_0^{\pm}(z)/\dot{\phi}_0^{\pm}(0)$ (see also \eqref{com1} and \eqref{c22}), we have 
\begin{eqnarray*} 
\begin{array}{lll}
\dot{\rho}_1^-(0;\kappa;1,0) &=& - \int_{-\infty}^0
(\dot{\phi}_0^-(\zeta) )^2 \{2x_0x_1\tilde{f}_u^- \\[0.2cm]
 & &+ \ x_0^2 \tilde{f}_{uu}^-\phi_1^-  - x_0^2 \kappa \}d{\zeta}/(\dot{\phi}_0^-(0))^2 \\[0.2cm]
 &=&  \ddot{\phi}_1^-(0)/\dot{\phi}_0^-(0) - \dot{\phi}_1^-(0) \ddot{\phi}_0^-(0)/(\dot{\phi}_0^-(0))^2 \\[0.2cm]
& &  + \ x_0^2 \kappa\int_{-\infty}^0 (\dot{\phi}_0^-(\zeta) )^2 d\zeta/(\dot{\phi}_0^-(0))^2
\end{array}
\end{eqnarray*}
and 
\begin{eqnarray*} 
\begin{array}{lll}
\dot{\rho}_1^+(0;\kappa;1,0) &=& - \int^{\infty}_0
(\dot{\phi}_0^+(\zeta) )^2 \{2(1-x_0)x_1\tilde{f}_u^+ \\[0.2cm]
& & - \ (1-x_0)^2 \tilde{f}_{uu}^+\phi_1^+  + (1-x_0)^2 \kappa \}d{\zeta}/(\dot{\phi}_0^+(0))^2 \\[0.2cm]
& = & \ddot{\phi}_1^+(0)/\dot{\phi}_0^+(0) - \dot{\phi}_1^+(0) \ddot{\phi}_0^+(0)/(\dot{\phi}_0^+(0))^2 \\[0.2cm]
& & - \ (1-x_0)^2 \kappa\int^{\infty}_0 (\dot{\phi}_0^+(\zeta) )^2 d\zeta/(\dot{\phi}_0^+(0))^2. 
\end{array}
\end{eqnarray*}
Then, noting that $(1-x_0)\dot{\phi}_0^-(0) = x_0\dot{\phi}_0^+(0)$ by $\Phi_0=0$, and $(1-x_0)^2\ddot{\phi}_0^-(0) =- x_0^2(1-x_0)^2f(\alpha,v^*) = x_0^2\ddot{\phi}_0^+(0)$ from \eqref{b18} and \eqref{b29}, and 
substituting these equations into the definitions of $a_{11}$ and $c_{11}$, we observe 
\begin{eqnarray} \label{eqapp3}
\begin{array}{l}
\hspace{0.5cm} a_{11} - c_{11} = \dot{\rho}^-_1(0;\kappa;1,0)/x_0 - \dot{\rho}^+_1(0;\kappa;1,0)/(1-x_0) \\[0.2cm]
\hspace{2cm}  - \ x_1\{\dot{\rho}^-_0(0;1,0)/(x_0)^2 +  \dot{\rho}^+_0(0;1,0)/(1-x_0)^2\} \\[0.2cm]
\hspace{0.5cm}   = \  \ddot{\phi}_1^-(0)/(x_0\dot{\phi}_0^-(0))  - \ddot{\phi}_1^+(0)/((1-x_0)\dot{\phi}_0^+(0)) \\[0.2cm]
\hspace{1cm}  - \ \dot{\phi}_1^-(0) \ddot{\phi}_0^-(0)/(x_0(\dot{\phi}_0^-(0))^2) + \dot{\phi}_1^+(0) \ddot{\phi}_0^+(0)/((1-x_0)(\dot{\phi}_0^+(0))^2) \\[0.2cm]
\hspace{0.5cm}  + \ \kappa \{ x_0 \int_{-\infty}^0 (\dot{\phi}_0^-(\zeta) )^2 d\zeta/(\dot{\phi}_0^-(0))^2   +  (1-x_0) \int^{\infty}_0 (\dot{\phi}_0^+(\zeta) )^2 d\zeta/(\dot{\phi}_0^+(0))^2\} \\[0.2cm]
\hspace{1cm} - \ x_1 \{\ddot{\phi}_0^-(0)/(x_0^2\dot{\phi}_0^-(0))  + \ddot{\phi}_0^+(0)/((1-x_0)^2\dot{\phi}_0^+(0))\}.
\end{array}
\end{eqnarray}
\par Setting $z=0$ in \eqref{b19} with $U_1^- = 0 = V_1^- $, we have $\ddot{\phi}_1^-(0) = -2x_0x_1f(\alpha,v^*)$. Hence, noting  $\ddot{\phi}_0^-(0) = -x_0^2f(\alpha,v^*)$ by \eqref{b18}, 
we obtain 
\begin{equation}  \label{eqappa3}
\ddot{\phi}_1^-(0) = 2x_1 \ddot{\phi}_0^-(0)/x_0. 
\end{equation}
Similarly, we obtain
\begin{equation}  \label{eqappa4}
\ddot{\phi}_1^+(0) = -2x_1 \ddot{\phi}_0^+(0)/(1-x_0).
\end{equation}
On the other hand, it follows from \eqref{b21} and \eqref{b32} that
$$ 
\dot{\phi}^-_1(0) = -2 x_0 x_1 \int_{h^-(v^*)}^\alpha f(u,v^*)du / \dot{\phi}^-_0(0)
$$
and 
$$ 
\dot{\phi}^+_1(0) = -2 (1-x_0) x_1 \int^{h^+(v^*)}_\alpha f(u,v^*)du / \dot{\phi}^+_0(0).
$$
Moreover, owing to the relations \eqref{b36} and \eqref{b37}, we have 
\begin{equation} \label{eqapp7}
\begin{array}{l}
\dot{\phi}^-_1(0) =  x_1 \dot{\phi}^-_0(0)/x_0 \quad \mbox{and}  \quad \dot{\phi}^+_1(0) = - x_1 \dot{\phi}^+_0(0)/(1-x_0).
\end{array} 
\end{equation}
Substituting \eqref{eqappa3}, \eqref{eqappa4} and \eqref{eqapp7} into \eqref{eqapp3}, we obtain 
\begin{eqnarray*}  
\begin{array}{l}
 a_{11} - c_{11} =  \kappa \{ x_0 \int_{-\infty}^0 (\dot{\phi}_0^-(\zeta) )^2 d\zeta/(\dot{\phi}_0^-(0))^2   +  (1-x_0) \int^{\infty}_0 (\dot{\phi}_0^+(\zeta) )^2 d\zeta/(\dot{\phi}_0^+(0))^2\} \\[0.2cm]
\hspace{2cm} + \ x_1 \left[ 2 \{ \ddot{\phi}_0^-(0)/(x_0^2\dot{\phi}_0^-(0)) + \ddot{\phi}_0^+(0)/((1-x_0)^2\dot{\phi}_0^+(0)) \} \right.\\[0.2cm]
\hspace{2.5cm} - \ \{\ddot{\phi}_0^-(0)/(x_0^2\dot{\phi}_0^-(0))  + \ddot{\phi}_0^+(0)/((1-x_0)^2\dot{\phi}_0^+(0))\} \\[0.2cm]
\hspace{2.5cm} \left.-  \{\ddot{\phi}_0^-(0)/(x_0^2\dot{\phi}_0^-(0))  + \ddot{\phi}_0^+(0)/((1-x_0)^2\dot{\phi}_0^+(0))\}\right] \\[0.2cm]
\hspace{1cm} = \ \kappa \{ x_0 \int_{-\infty}^0 (\dot{\phi}_0^-(\zeta) )^2 d\zeta/(\dot{\phi}_0^-(0))^2   +  (1-x_0) \int^{\infty}_0 (\dot{\phi}_0^+(\zeta) )^2 d\zeta/(\dot{\phi}_0^+(0))^2\} \\[0.2cm]
 \hspace{1cm} = \ \displaystyle \kappa \int^\infty_{-\infty} (\dot{ W}(z))^2 dz/(\dot{ W}(0))^2. 
\end{array}
\end{eqnarray*}

Using the definition of $a_{22}$, \eqref{eqn:C4} and $\pi_0^-(z;1,0) = -\rho_0^-(z;1,0)/D$, we have
\begin{eqnarray*}
\begin{array}{lll}
a_{22} & = & \dot{\pi}^-_1(0;\kappa;1,0)/x_0 -  x_1 \dot{\pi}^-_0(0;1,0)/(x_0)^2 \\[0.2cm]
 & = & - \dot{\rho}^-_1(0;\kappa;1,0)/(x_0 D) + x_0 \kappa \int_{-\infty}^0 \rho^-_0(\zeta;1,0)d\zeta/D \\[0.2cm]
 & & + \ x_1 \dot{\rho}_0^-(0;1,0)/(x_0^2 D) \\[0.2cm]
 & = & - \dot{\rho}^-_1(0;\kappa;1,0)/(x_0 D) + x_0 \kappa (\alpha - h^-(v^*))/(D \dot{\phi}_0^-(0)) \\[0.2cm]
  & & + \ x_1 \dot{\rho}_0^-(0;1,0)/(x_0^2 D).
\end{array} 
\end{eqnarray*}  
Similarly, using the definition of $c_{22}$, \eqref{C6} and $\pi_0^+(z;1,0) = -\rho_0^+(z;1,0)/D$, we have
\begin{eqnarray*}
\begin{array}{lll}
c_{22} & = & \dot{\pi}^+_1(0;\kappa;1,0)/(1-x_0) +  x_1 \dot{\pi}^+_0(0;1,0)/(1-x_0)^2 \\[0.2cm]
 & = & - \dot{\rho}^+_1(0;\kappa;1,0)/\{(1-x_0)D\} - (1-x_0) \kappa \int^{\infty}_0 \rho^+_0(\zeta;1,0)d\zeta/D \\[0.2cm]
   &  & - \ x_1 \dot{\rho}_0^+(0;1,0)/\{(1-x_0)^2 D\} \\[0.2cm]
 & = &  - \dot{\rho}^+_1(0;\kappa;1,0)/\{(1-x_0)D\} + (1-x_0) \kappa (\alpha - h^+(v^*))/(D \dot{\phi}_0^+(0)) \\[0.2cm]
   &  & - \ x_1 \dot{\rho}_0^+(0;1,0)/\{(1-x_0)^2 D\}.
\end{array} 
\end{eqnarray*}
Hence, noting $\dot{ W}(0) = \dot{\phi}_0^-(0)/x_0 = \dot{\phi}_0^+(0)/(1-x_0)$ by \eqref{b00009},
we obtain
\begin{eqnarray*}
\begin{array}{lll}
a_{22}-c_{22} & = & \left[ -(a_{11}-c_{11}) + \kappa \left\{x_0 (\alpha - h^-(v^*))/\dot{\phi}_0^-(0) \right. \right. \\[0.2cm]
  & & \left. \left. - \ (1-x_0)(\alpha - h^+(v^*))/\dot{\phi}_0^+(0) \right\} \right]/D \\[0.2cm]
 & = & \left[ -(a_{11}-c_{11}) + \kappa (h^+(v^*) - h^-(v^*))/\dot{ W}(0)
\right]/D.
\end{array} 
\end{eqnarray*}

Using \eqref{com1} and 
the relations 
\begin{eqnarray*}
\begin{array}{l}
\displaystyle 
\int_{-\infty}^0 \dot{\phi}_0^-(\zeta) \tilde{f}_u^- d\zeta = f(\alpha,v^*) \ \
\text{and} \ \ 
\int_{-\infty}^0 \dot{\phi}_0^-(\zeta) \tilde{f}_v^- d\zeta = \int_{h^-(v^*)}^\alpha 
f_v(u,v^*) du, 
\end{array} 
\end{eqnarray*}
we have 
\begin{eqnarray*}
\begin{array}{lll}
b_{10} & = & \dot{\rho}^-_0(0;0,1)/x_0 \\[0.2cm]
 & = & \displaystyle \frac 1{x_0} \frac {f_v^-}{f_u^-}\frac {\ddot{\phi}_0^-(0)}{\dot{\phi}_0^-(0)} - x_0 \int_{-\infty}^0 \dot{\phi}_0^-(\zeta)\left\{ \tilde{f}_u^- \left( - \frac {f_v^-}{f_u^-}\right) + \tilde{f}_v^- \right\} d\zeta/\dot{\phi}_0^-(0) \\[0.4cm]
 & = & \displaystyle \frac 1{x_0}  \frac {f_v^-}{f_u^-}\frac {\ddot{\phi}_0^-(0)}{\dot{\phi}_0^-(0)} + x_0 \left\{\frac {f_v^-}{f_u^-}f(\alpha,v^*) - \int_{h^-(v^*)}^\alpha 
f_v(u,v^*) du\right\}/\dot{\phi}_0^-(0) \\[0.4cm]
 & = & \displaystyle - x_0 \int_{h^-(v^*)}^\alpha f_v(u,v^*) du/\dot{\phi}_0^-(0), 
\end{array} 
\end{eqnarray*}
where we used $\ddot{\phi}_0^-(0) + x_0^2f(\alpha,v^*) = 0$ by \eqref{b18}. 
Similarly, we obtain 
\begin{eqnarray*}
\begin{array}{lll}
d_{10} & = & \dot{\rho}^+_0(0;0,1)/(1-x_0) \\[0.2cm]
 & = & \displaystyle \frac 1{(1-x_0)}  \frac {f_v^+}{f_u^+}\frac {\ddot{\phi}_0^+(0)}{\dot{\phi}_0^+(0)} \\[0.4cm]
  & & \displaystyle + \ (1-x_0) \left\{\frac {f_v^+}{f_u^+}f(\alpha,v^*) + \int^{h^+(v^*)}_\alpha f_v(u,v^*) du\right\}/\dot{\phi}_0^+(0) \\[0.4cm] 
 & = & \displaystyle (1-x_0) \int^{h^+(v^*)}_\alpha f_v(u,v^*) du/\dot{\phi}_0^+(0). 
\end{array} 
\end{eqnarray*}
Then, noting $(1-x_0)\dot{\phi}_0^-(0) = x_0\dot{\phi}_0^+(0)$ and $\dot{ W}(0) =\dot{\phi}_0^-(0)/x_0 =  \dot{\phi}_0^+(0)/(1-x_0)$ we have 
$$
b_{10}-d_{10} =  - \displaystyle \int^{h^+(v^*)}_{h^-(v^*)} f_v(u,v^*) du / \dot{ W}(0). 
$$

Finally, since
\begin{eqnarray*}
\begin{array}{lll}
b_{21} & = & (Q^-_{1,y}(1;\kappa;1) + \dot{\pi}^-_0(0;0,1))/x_0 \\[0.2cm]
 & = & - \dot{\rho}^-_0(0;0,1)/(x_0 D) + \dis\frac {x_0 \kappa}D \left\{1 - \frac {f_v^-}{f_u^-} \right\}, \\[0.4cm]
d_{21} & = & (Q^+_{1,y}(0;\kappa;1) + \dot{\pi}^+_0(0;0,1))/(1-x_0) \\[0.2cm]
 & = & - \dot{\rho}^+_0(0;0,1)/((1-x_0) D) - \dis\frac {(1-x_0) \kappa}D \left\{1 - \frac {f_v^+}{f_u^+} \right\}
\end{array} 
\end{eqnarray*}
by $\pi_0^\pm(z;a,b) = -\rho_0^\pm(z;a,b)/D$, a direct calculation shows that 
$$
b_{21}-d_{21} = \displaystyle  \left[ -(b_{10}-d_{10}) + \kappa \int^1_0 \left( 
\frac {f_u^* - f_v^*}{f_u^*} \right) dx \right]/D.
$$
Thus the proof of Lemma \ref{l3-6} is completed. 
%\vspace{0.5cm}

\end{document}